\documentclass[times,sort&compress,3p]{elsarticle}
\journal{Journal of Multivariate Analysis}
\usepackage[labelfont=bf]{caption}

\usepackage{amsfonts}
\usepackage{amsmath}
\usepackage{amsthm}
\usepackage{amssymb}
\usepackage{graphicx}
\usepackage{multirow}
\usepackage{longtable}
\usepackage[colorlinks=true,linkcolor=RawSienna,citecolor=RawSienna,urlcolor=RawSienna,bookmarksopen=true,pdfstartview=FitB]{hyperref}
\usepackage{natbib}
\usepackage[dvipsnames]{xcolor}
\usepackage{float}

\newtheorem{theorem}{Theorem}

\newtheorem{corollary}{Corollary}

\newtheorem{definition}{Definition}
\newtheorem{example}{Example}
\newtheorem{lemma}{Lemma}

\begin{document}

\begin{frontmatter}

\title{Uniformly consistently estimating the proportion of false null hypotheses via Lebesgue-Stieltjes integral equations}

\author[A1]{Xiongzhi Chen\corref{mycorrespondingauthor}}

\address[A1]{Department of Mathematics and Statistics, Washington State University, Pullman, WA 99164, USA}

\cortext[mycorrespondingauthor]{Corresponding author. Email address: \url{xiongzhi.chen@wsu.edu}}

\begin{abstract}
The proportion of false null hypotheses is a very important quantity in
statistical modelling and inference based on the two-component mixture model
and its extensions, and in control and estimation of the false discovery rate
and false non-discovery rate. Most existing estimators of this proportion
threshold p-values, deconvolve the mixture model under constraints on its
components, or depend heavily on the location-shift property of distributions.
Hence, they usually are not consistent, applicable to non-location-shift
distributions, or applicable to discrete statistics or p-values. To eliminate
these shortcomings, we construct uniformly consistent estimators of the
proportion as solutions to Lebesgue-Stieltjes integral equations. In
particular, we provide such estimators respectively for random variables whose
distributions have Riemann-Lebesgue type characteristic functions, form
discrete natural exponential families with infinite supports, and form natural
exponential families with separable moment sequences. We provide the speed of
convergence and uniform consistency class for each such estimator under
independence. In addition, we provide example distribution families for which
a consistent estimator of the proportion cannot be constructed using our
techniques.
\end{abstract}

\begin{keyword}
Analytic functions \sep
Bessel functions \sep
concentration inequalities \sep
Fourier transform \sep
Lambert W functions \sep
Lebesgue-Stieltjes integral equations \sep
Mellin transform \sep
natural exponential family \sep
proportion of false null hypotheses.
\MSC[2010] Primary 35C05, 60E05 \sep
Secondary 60F10
\end{keyword}

\end{frontmatter}

\section{Introduction}

\label{sec1intro}

The proportion of false null hypotheses and its dual, the proportion of true
null hypotheses, play important roles in statistical modelling and multiple
hypotheses testing. For example, they are components of the two-component
mixture model of \cite{Efron:2001}, its extensions by
\cite{Cai:2009}, \cite{Liu:2016} and \cite{Ploner:2006} and their induced statistics including
the \textquotedblleft local false discovery rate (local FDR)\textquotedblright%
\ of \cite{Efron:2001} and \textquotedblleft positive FDR
(pFDR)\textquotedblright\ and q-value of \cite{Storey:2003c}. They also form
the optimal discovery procedure (ODP) of \cite{Storey:2007}. However, without
information on the proportions, decision rules based on the local FDR cannot
be implemented, and none of the pFDR, q-value and ODP decision rule can be
computed in practice. On the other hand, these proportions form upper bounds
on the FDRs and false non-discovery rates (FNRs) of all FDR procedures
including those of \cite{Benjamini:2001}, \cite{Blanchard:2009}, \cite{Genovese:2004},
\cite{Sarkar:2006} and \cite{Storey:2004}. Therefore,
accurate information on either proportion helps better control and estimate
the FDR and FNR, and thus potentially enables an adaptive procedure to be more
powerful or have smaller FNR than its non-adaptive counterpart. Since neither
proportion is known in practice, it is very important to accurately estimate
the proportions.

In this work, we focus on consistently estimating the proportions without employing the two-component
mixture model, requiring p-values be identically distributed under the alternative hypotheses, or assuming that statistics or p-values have absolutely
continuous cumulative distribution functions (CDFs). The main motivation for
dealing with discrete p-values or statistics is the wide practice of FDR
control in multiple testing based on discrete data in genomics
\citep{Auer:2010,Robinson:2008}, genetics \citep{Gilbert:2005}, vaccine
efficacy studies \citep{Mehrotra:2004}, drug safety monitoring
\citep{Chen:2015discretefdr} and other areas, where Binomial Test (BT),
Fisher's Exact Test (FET) and Exact Negative Binomial Test (NBT) have been
routinely used to test individual hypotheses. In the sequel, we refer to an
estimator of the proportion of true (or false) null hypotheses as a ``null (or
alternative) proportion estimator''.

\subsection{A brief review on existing proportion estimators}

There are many estimators of the proportions, and their constructions can be
roughly categorized into 4 classes: (1) thresholding p-values
\citep{Storey:2002,Storey:2004}; (2) deconvolving the two-component mixture
model for the distribution of p-values or statistics, modulo identifiability
conditions \citep{Genovese:2004,Kumar:2016,swanepoel1999}; (3) bounding the
proportions via the use of uniform empirical process \citep{Meinshausen:2006};
(4) Fourier transform for Gaussian family or mixtures with a Gaussian
component \citep{Jin:2008,Jin:2007,Jin:2010a}. Further, perhaps the only
existing consistent proportion estimators are those mentioned right above, and only those in class (3) and (4) are consistent when
the proportion tends to zero (see (\ref{defConsistency}) for the definition of consistency).
However, they have the following disadvantages: (i) the construction in class
(4) based on Fourier transform does not work if the family of distributions of
statistics does not consist at least one component from a location-shift
family; (ii) consistency of estimators in class (1) and (2) requires the
two-component mixture model and various regularity conditions such as
concavity, smoothness, purity (defined by \cite{Genovese:2004}) of the
alternative component, or the conditions in Lemma 3 or Lemma 4 of
\cite{Kumar:2016}; (iii) consistency of the estimator in class (3) requires
p-values to have continuous distributions, be independent, and be identically
distributed under the alternative hypotheses. An extended comparison for popular proportion
estimators is provided in \autoref{tb_compEst}.

In practice, the requirements needed by the consistent proportion estimators
mentioned above are often not met. For example, for multiple testing where
each statistic follows a Chi-square or Binomial distribution the involved
distributions do not form a location-shift family, the two-component mixture
model is inappropriate for p-values of BTs, FETs or NBTs with different
marginal counts, and p-values under the alternative hypotheses cannot be
identically distributed when the signal levels are different across individual
hypotheses under the alternative. Further, almost all existing proportion
estimators were initially designed for p-values or statistics that have
continuous distributions, and it is not clear yet whether the null proportion
estimators in classes (2) to (4) are reciprocally conservative even though the
estimator of \cite{Storey:2004} is so for independent p-values whose null
distributions are uniform on $[0,1]$ (see Corollary 13 of
\cite{Blanchard:2009}). Therefore, it is important to develop new,
consistent proportion estimators that mitigate or eliminate the shortcomings
of existing ones.

{\footnotesize
\LTcapwidth=\textwidth
\begin{longtable}{|p{2cm}|p{2cm}|p{10cm}|}

    \hline
   \multirow{5}{2cm}{Thresholding estimators}
    & References & \cite{Chen:2016,Storey:2002,Storey:2004} \\
   \cline{2-3}
   & Method &  Thresholding p-values  \\
   \cline{2-3}
   &Disadvantages &  Only applicable to p-values \\
   \cline{2-3}
  &Consistency &  Consistency when p-values are independent, identically distributed, follow an identifiable, two-component mixture model, and have an absolutely continuous CDF with some regularity properties, which has been proved by  \cite{Genovese:2004} for Storey's estimator\\
  \cline{2-3}
  & Scope & Applicable to both discrete and continuous data\\

  \hline
  \hline
   \multirow{5}{2cm}{Deconvolution estimators: I}
   & References & \cite{Genovese:2004,swanepoel1999} \\
   \cline{2-3}
   & Method &  Deconvolving p-value densities  \\
   \cline{2-3}
  &Modelling assumptions &   p-values need to be identically distributed, follow a two-component mixture model, and have an absolutely continuous CDF\\
     \cline{2-3}
  &Consistency &  Joint assumptions for consistency: independence, an identifiable model, regularity conditions on the alternative component in the model, and the proportion being fixed\\
  \cline{2-3}
  & Scope & Applicable to continuous data only\\

   \hline
   \hline
   \multirow{4}{2cm}{Deconvolution estimators: II}
   & References & \cite{Kumar:2016} \\
   \cline{2-3}
   & Method &  Deconvolving p-value densities  \\
   \cline{2-3}
  &Modelling assumptions &   p-values need to be identically distributed and follow a two-component mixture model\\
     \cline{2-3}
  &Consistency &  Joint assumptions for consistency: independence, an identifiable model, regularity conditions on the alternative component in the model, and the proportion being fixed\\
  \cline{2-3}
   & Scope & Applicable to both discrete and continuous data\\

   \hline
   \hline

   \multirow{4}{2cm}{Deconvolution estimators: III}
   & References & \cite{Langaas:2005} \\
   \cline{2-3}
   & Method &  Deconvolving p-value densities  \\
   \cline{2-3}
  &Modelling assumptions &   p-values need to be identically distributed, follow a two-component mixture model, and have an absolutely continuous, concave CDF\\
     \cline{2-3}
  &Consistency &  Consistency unknown\\
  \cline{2-3}
  & Scope & Applicable to continuous data only\\
   \hline
   \hline

   \multirow{4}{2cm}{Fourier transform based estimators}
   & References & \cite{Chen:2018,Jin:2008,Jin:2007,Jin:2010a} \\
   \cline{2-3}
   & Method &  Fourier transform and convolution\\
   \cline{2-3}
  &Modelling assumptions; disadvantages&  The family of distributions of statistics should consist at least
one component from a location-shift family; implementation needs information on the scale parameter of the location-shift component\\
   \cline{2-3}
   &Consistency &  Uniform consistency under weak dependence when the proportion is fixed or converges to $0$ at certain rate\\
  \cline{2-3}
  & Scope & Applicable to continuous data only\\

  \hline
  \hline
  \multirow{4}{2cm}{Goodness-of-fit estimators}
   & References & \cite{Meinshausen:2006} \\
   \cline{2-3}
   & Method &  Measuring the relative excess of the uniform empirical process\\
   \cline{2-3}
  &Modelling assumptions &   p-values need to have absolutely continuous CDFs, be independent,
and be identically distributed under the alternative hypotheses\\
   \cline{2-3}
   &Consistency &  Consistency under independence when the proportion is fixed or converges to $0$ at certain rate\\
  \cline{2-3}
  & Scope & Applicable to continuous data only\\

  \hline
  \hline
  \multirow{4}{2cm}{Integro-equation based estimators}
   & References & This paper \\
   \cline{2-3}
   & Method &  Approximating the status of each hypothesis via solutions to Lebesgue-Stieltjes integral equations\\
   \cline{2-3}
  &Modelling assumptions and disadvantages &   Statistics need to have CDFs that have Riemann-Lebesgue type characteristic functions, that are members of an NEF with support $\mathbb{N}$, or that are members of an NEF with separable moment sequence; implementation needs information on the scale parameters of the CDFs, on the supremum norm of the parameter vector, or on the infimum of a transform of this vector\\
   \cline{2-3}
   &Consistency &  Uniform consistency under independence when the proportion is fixed or converges to $0$ at certain rate\\
  \cline{2-3}
  & Scope & Applicable to both continuous and discrete data\\
\hline
\hline
\caption{A comparison of some popular estimators of the proportion of false null hypotheses.}
\label{tb_compEst}
\end{longtable}
}

\subsection{Main contributions and summary of results}

We generalize \textquotedblleft Jin's Strategy\textquotedblright\ provided in
Section 2.1 of \cite{Jin:2008} that estimates the proportions by approximating
the indicator function of the true status of a null hypothesis and has only
been implemented for Gaussian family or mixtures with a Gaussian
component, and in much more general settings construct proportion estimators
as solutions to a specific type of Lebesgue-Stieltjes integral equation in the
complex domain.
Generalizing Jin's Strategy to non-location-shift families
is highly nontrivial since this requires solving Lebesgue-Stieltjes integral equations and outside location-shift families solutions to these equations are very hard to find or do not exist. Even after a solution is found, to prove the uniform consistency of the proposed estimator we have to derive concentration inequalities for non-Lipschitz transforms of independent, unbounded random variables, and this is highly nontrivial and an open area in probability theory.

Our methodology (referred to as \textquotedblleft the
Strategy\textquotedblright) lifts Jin's strategy to its full generality. In
addition to the key advantage of Jin's strategy, i.e., being independent from
but partially applicable to the two-component mixture model, the Strategy
applies to random variables that have continuous or discrete distributions or
form non-location-shift families. Further, it produces uniformly consistent
proportion estimators in the dense and moderately sparse regimes (that are defined at the beginning of \autoref{secModel}).

Our main contributions are summarized as follows:

\begin{itemize}
\item \textquotedblleft Construction I\textquotedblright: Construction of
uniformly consistent proportion estimators for random variables whose
distributions have Riemann-Lebesgue type characteristic functions (RL type
CFs); see \autoref{Def} in \autoref{secEstimator}. In particular, this covers the construction
for random variables from several location-shift families, including that in
\cite{Jin:2008} for Gaussian and Laplace distributions.

\item \textquotedblleft Construction II\textquotedblright: Construction of
uniformly consistent proportion estimators for random variables whose
distributions form natural exponential families with infinite supports; see
\autoref{SecConstructionII} and \autoref{SecConsistent2}. In particular, this
covers six of the twelve natural exponential families with cubic variance
functions (NEF-CVFs) proposed by \cite{Letac:1990} such as Poisson and
Negative Binomial families. Estimators from Construction II are perhaps the
first consistent estimators that are applicable to discrete random variables.

\item \textquotedblleft Construction III\textquotedblright: Construction of
uniformly consistent proportion estimators for random variables whose
distributions form natural exponential families with separable moment
sequences; see \autoref{SecConstructionII} and \autoref{secConsistent3}. In
particular, this covers Gamma family which includes Exponential and central Chi-square
families as special cases.

\item For each constructed estimator, its speed of convergence and uniform
consistency class is provided under independence; see \autoref{secEstimator},
\autoref{SecConsistent2} and \autoref{secConsistent3}.

\item Concentration inequalities for sums of certain non-Lipschitz transforms of independent but not necessarily almost surely bounded random variables.
\end{itemize}

Specifically, Construction I employs Fourier transform and an extended
Riemann-Lebesgue Lemma of \cite{Costin:2016}. Since the set of distributions
with RL type CFs contains several location-shift families, we give a unified
treatment of proportion estimators for location-shift families and reveal the
intrinsic mechanism of Fourier transform based construction. In contrast,
Construction II mainly uses generating functions (GFs), and Construction III
Mellin transform which can be regarded as inducing \textquotedblleft
multiplication-convolution equivalence\textquotedblright, in contrast to
Fourier transform inducing \textquotedblleft translation-convolution
equivalence\textquotedblright. As negative results, we show that the Strategy
is not implementable for Inverse Gaussian family (whose densities have a very
special structure; see \cite{Letac:1990} for a definition of the family) and Binomial family (which is discrete but with a finite
support); see \autoref{SecNonexist}.

We provide upper bounds on the variances of the proportion estimators, show
their uniform consistency, and provide their speeds of convergence for
consistency under independence. Additionally, we have found the following. For
estimators given by Construction I, uniform consistency in frequency domain
(see \autoref{DefUniformConsistency}) can also be achieved due to the global
Lipschitz property of the construction, and their uniform consistency classes
(see \autoref{DefUniformConsistency}) can be ordered via set inclusion
according to the magnitudes of the moduli of the corresponding CFs; see
\autoref{ThmLocationShift} and \autoref{CorInclusion1}. In contrast, for
estimators given by Constructions II and III, uniform consistency in frequency
domain is very hard to obtain since the constructions are not Lipschitz
transforms of the involved random variables.

For an estimator given by Construction I or Construction II where GFs have
finite radii of convergence, its speed of convergence and uniform consistency
class do not depend on the supremum norm of the parameter vector (see
\autoref{ThmLocationShift}, \autoref{CorInclusion1} and \autoref{CoroII}),
enabling the estimator to be fully data-adaptive. In contrast, for estimators
given by other instances of Construction II and Construction III, their speeds
of convergence and uniform consistency classes may depend explicitly on the
supremum norm of the parameter vector or the infimum of a transform of this
vector (see \autoref{II-consistency} and \autoref{III-uniformConsistent}).
However, since the Strategy approximates indicator functions of the status of
individual null hypotheses, these speeds and classes always depend on the
minimal magnitude of the differences between parameters of interest and their
reference value regardless of which among the three constructions is used ---
a universality phenomenon for the Strategy.

\subsection{Notations and conventions}

Throughout the article, we use the following conventions and notations: $C$
denotes a generic, positive constant whose values may differ at different
occurrences; $O\left(  \cdot\right)  $ and $o\left(  \cdot\right)  $ are
respectively Landau's big O and small o notations; $\mathbb{E}$ and
$\mathbb{V}$ are respectively the expectation and variance with respect to the
probability measure $\Pr$; $\mathbb{R}$ and $\mathbb{C}$ are respectively the
set of real and complex numbers; $\Re$, $\Im$ and $\arg$ denote the real,
imaginary part and argument of a complex number, respectively; $\mathbb{N}$
denotes the set of non-negative integers, and $\mathbb{N}_{+}=\mathbb{N}%
\setminus\left\{  0\right\}  $; $\delta_{y}$ is the Dirac mass at
$y\in\mathbb{R}$; $\nu$ the Lebesgue measure, and when it is clear that an
integral is with respect to $\nu$, the usual notation $d\cdot$ for
differential will be used in place of $\nu\left(  d\cdot\right)  $; for $1\leq
p\leq\infty$ and $A\subseteq\mathbb{R}$, $\left\Vert f\right\Vert _{p}=\left\{
\int_{A}\left\vert f\left(  x\right)  \right\vert ^{p}\nu\left(  dx\right)
\right\}  ^{1/p}$ and $L^{p}\left(  A\right)  =\left\{  f:\left\Vert
f\right\Vert _{p}<\infty\right\}  $, for which $\left\Vert f\right\Vert
_{\infty}$ is the essential supremum of $f$; for a set $A$ in $\mathbb{R\ }%
$and a scalar $a$, $1_{A}$ is the indicator of $A$ and $A-a=\left\{  x-a:x\in
A\right\}  $; $\varnothing$ is the empty set; for $x\in\mathbb{R}$, $\lfloor
x\rfloor$ denotes the integer part of $x$; the extended complex plane
$\mathbb{C}\cup\left\{  \infty\right\}  $ is identified with the Riemann
sphere so that $\infty$ corresponds to its north pole; for $x\in\mathbb{R}$,
$x\rightarrow\infty$ means $x\rightarrow+\infty$, and $x\rightarrow-\infty$
means $x\rightarrow-\infty$; for positive functions $a^{\prime}\left(
x\right)  $ and $b^{\prime}\left(  x\right)  $, $a^{\prime}\left(  x\right)
\sim b^{\prime}\left(  x\right)  $ means that $\lim_{x\rightarrow\infty
}a^{\prime}\left(  x\right)  \left\{  b^{\prime}\left(  x\right)  \right\}
^{-1}=1$; $\mathbb{R}^{\mathbb{N}}$ is the $\aleph$-Cartesian product of $\mathbb{R}$, where
$\aleph$ is the cardinality of $\mathbb{N}$.

\subsection{Organization of paper}

The rest of the article is organized as follows. In \autoref{secModel} we
formulate the problem of proportion estimation and state the Strategy for
constructing proportion estimators. In \autoref{secEstimator} we develop
uniformly consistent proportion estimators when the CDFs of random variables
have Riemann-Lebesgue type characteristic functions. In
\autoref{SecConstructionII}, we construct proportion estimators when the
distributions of random variables form discrete NEFs with infinite supports or
form NEFs with separable moment sequences. In \autoref{SecNonexist}, we
provide two families of random variables for which the Strategy cannot be
applied. In \autoref{SecConsistent2} and \autoref{secConsistent3} we justify
the uniform consistency of the constructed proportion estimators for random
variables whose distributions form NEFs with infinite supports or with
separable moment sequences, respectively. In \autoref{secSim} we conduct a simulation study on the proposed estimators, with comparison to those of \cite{Meinshausen:2006} and \cite{Jin:2008}.
We end the article with a discussion in \autoref{SecConcAndDisc}, together with several topics worthy of future investigations.
The supplementary material contains proofs of all theoretical results,
provides evidence that Construction III does not apply to Ressel or Hyperbolic
Cosine families (defined by \cite{Letac:1990}), and a discussion on uniform consistency of proportion
estimators in frequency domain for Construction II and III.

\section{The estimation problem and strategy}

\label{secModel}

In this section, we formulate the estimation problem and generalize Jin's
strategy. Let $z_{i}$, $1\leq i\leq m$ be $m$ random variables each with mean
or median $\mu_{i}$, such that, for a fixed value $\mu_{0}$ for the mean or
median and some integer $m_{0}$ between $0$ and $m$, $\mu_{i}=\mu_{0}$ for
each $i=1,\ldots,m_{0}$ and $\mu_{i}\neq\mu_{0}$ for each $i=m_{0}+1,\ldots
,m$. Consider simultaneously testing the null hypothesis $H_{i0}:\mu_{i}%
=\mu_{0}$ versus the alternative hypothesis $H_{i1}:\mu_{i}\neq\mu_{0}$ for
$1\leq i\leq m$. Let $I_{0,m}$ $=\left\{  1\leq i\leq m:\mu_{i}=\mu
_{0}\right\}  $ and $I_{1,m}=\left\{  1\leq i\leq m:\mu_{i}\neq\mu
_{0}\right\}  $. Then the cardinality of $I_{0,m}$ is $m_{0}$, the proportion
of true null hypothesis (\textquotedblleft null proportion\textquotedblright%
\ for short) is defined as $\pi_{0,m}=m^{-1}m_{0}$, and the proportion of
false null hypotheses (\textquotedblleft alternative
proportion\textquotedblright\ for short) $\pi_{1,m}=1-\pi_{0,m}$. In other
words, $\pi_{0,m}$ is the proportion of random variables that have a
prespecified mean or median. Our target is to consistently estimate $\pi
_{1,m}$ as $m\rightarrow\infty$ when $\left\{  z_{i}\right\}  _{i=1}^{m}$ are independent.

We will adopt the following convention from \cite{Jin:2008}: the dense regime is represented by $\pi_{1,m}$ such that $\liminf_{m\rightarrow\infty}\pi_{1,m}>0$, the moderately sparse regime by $\pi_{1,m}=Cm^{-s}$ for $s \in \left(0,0.5\right)$, the critically sparse regime by $\pi_{1,m}=Cm^{-0.5}$, and the very sparse regime by $\pi_{1,m}=Cm^{-s}$ for $s \in \left(0.5,1\right)$, where $C>0$ is a constant.

\subsection{The Strategy for proportion estimation}

Let $\mathbf{z}=\left(  z_{1},\ldots,z_{m}\right)  ^{\top}$ and $\boldsymbol{\mu
}=\left(  \mu_{1},...,\mu_{m}\right)  ^{\top}$. Denote by $F_{\mu_{i}}\left(
\cdot\right)  $ the CDF of $z_{i}$ for $1\leq i\leq m$ and suppose each
$F_{\mu_{i}}$ with $0\leq i\leq m$ is a member of a set $\mathcal{F}$ of CDFs
such that $\mathcal{F}=\left\{  F_{\mu}:\mu\in U\right\}  $ for some non-empty
$U$ in $\mathbb{R}$. For the rest of the paper, we assume that each $F_{\mu}$
is uniquely determined by $\mu$ and that $U$ has a non-empty interior.

To illustrate the intuitions behind Jin's strategy, we consider the setting
where each $z_{i}$ is Normally distributed with mean $\mu_{i}$ and standard
deviation $1$ (denoted by $z_{i}\sim\mathsf{Normal}\left(  \mu_{i},1\right)
$) and $\mu_{0}=0$. So, $\pi_{1,m}=m^{-1}\sum_{i=1}^{m}\left(  1-1_{\left\{
\mu_{i}=0\right\}  }\right)  $. If we can construct a function $\hat{\psi
}\left(  t,\mu\right)  $ such that $\hat{\psi}\left(  t,0\right)  =1$ for all
$t\in\mathbb{R}$ and $\lim_{t\rightarrow\infty}\hat{\psi}\left(  t,\mu\right)
=0$ for any $\mu\neq0$, then, for any fixed $\boldsymbol{\mu}$ and $m$,
$1_{\left\{  \mu_{i}=0\right\}  }=\hat{\psi}\left(  t,0\right)  $ when
$\mu_{i}=0$ and $\,1_{\left\{  \mu_{i}\neq0\right\}  }=1-\lim_{t\rightarrow
\infty}\hat{\psi}\left(  t,\mu_{i}\right)  $ when $\mu_{i}\neq0$. In other
words,
\[
\pi_{1,m}=\lim_{t\rightarrow\infty}m^{-1}\sum_{i=1}^{m}\left\{  1-\hat{\psi
}\left(  t,\mu_{i}\right)  \right\}  \text{ for any fixed }\boldsymbol{\mu
}\text{ and }m.
\]
One way to construct $\hat{\psi}$ is to utilize the Riemann-Lebesgue lemma
(see, e.g., \cite{Costin:2016}) and represent $\hat{\psi}$ as an integral.
For example, given any probability density function $\tilde{\omega}$ on
$\left[  -1,1\right]  $, setting
\[
\hat{\psi}\left(  t,\mu\right)  =\int\tilde{\omega}\left(  s\right)
\cos\left(  t\mu s\right)  ds
\]
immediately gives the desired $\hat{\psi}$. If $\tilde{\omega}$ is also even,
then $\hat{\psi}$ is the Fourier transform of $\tilde{\omega}$ evaluated at
$t\mu$. Notice that $\hat{\psi}$ is a deterministic function and that
$\tilde{\omega}$ is referred to as an \textquotedblleft averaging
function\textquotedblright\ in \autoref{secEstimator}. Even though a $\hat{\psi}$ has been
found, it cannot be used to estimate $\pi_{1,m}$ since $\hat{\psi}$ is not a
function of any $z_{i}$. So, the next step is to connect $\hat{\psi}\left(
t,\mu_{i}\right)  $ with $z_{i}$ probabilistically via a function $\hat
{K}\left(  x,t\right)  $, so that $\mathbb{E}\left\{  \hat{K}\left(
z_{i},t\right)  \right\}  =\hat{\psi}\left(  t,\mu_{i}\right)  $ for each $t$
and $i$. Once such a $\hat{K}$ is found, then $\pi_{1,m}^{\circ}\left(
t\right)  =m^{-1}\sum_{i=1}^{m}\left\{  1-\hat{K}\left(  z_{i},t\right)
\right\}  $ serves as an estimate of $\pi_{1,m}$ and can be very accurate when
$t$ is large. Specifically, setting%
\[
\hat{K}\left(  x,t\right)  =\int\tilde{\omega}\left(  s\right)  \exp\left(
2^{-1}t^{2}s^{2}\right)  \cos\left(  txs\right)  ds
\]
gives $\mathbb{E}\left\{  \hat{K}\left(  z_{i},t\right)  \right\}  =\hat{\psi
}\left(  t,\mu_{i}\right)  $ whenever $z_{i}\sim\mathcal{N}\left(  \mu
_{i},1\right)  $, and $\hat{\psi}\left(  t,\mu\right)  =\left\{  \hat{K}\left(
\cdot,t\right)  \ast\phi\right\}  \left(  \mu\right)  $, where $\phi$ is the
standard Normal density and $\ast$ denotes convolution. This fact can be found
in \cite{Jin:2008} or \cite{Chen:2018}. Finally, to show the consistency of
$\pi_{1,m}^{\circ}\left(  t\right)  $, we only need to control the difference
$\pi_{1,m}^{\circ}\left(  t\right)  -\mathbb{E}\left\{  \pi_{1,m}^{\circ
}\left(  t\right)  \right\}  $ as both $t$ and $m$ tend to $\infty$.

Now we state the Strategy below. For each fixed $\mu\in U$, if we can approximate
the indicator function $1_{\left\{  \mu\neq\mu_{0}\right\}  }$ by a function
$\psi\left(  t,\mu;\mu_{0}\right)  $ with $t\in\mathbb{R}$ satisfying
$\lim_{t\rightarrow\infty}\psi\left(  t,\mu_{0};\mu_{0}\right)  =1$ and
$\lim_{t\rightarrow\infty}\psi\left(  t,\mu;\mu_{0}\right)  =0$ for $\mu
\neq\mu_{0}$, then the \textquotedblleft phase function\textquotedblright%
\begin{equation*}
\varphi_{m}\left(  t,\boldsymbol{\mu}\right)  =\dfrac{1}{m}\sum_{i=1}%
^{m}\left\{  1-\psi\left(  t,\mu_{i};\mu_{0}\right)  \right\}  \label{eqd1}%
\end{equation*}
satisfies $\lim_{t\rightarrow\infty}\varphi_{m}\left(  t,\boldsymbol{\mu
}\right)  =\pi_{1,m}$ for any fixed $m$ and $\boldsymbol{\mu}$ and provides
the \textquotedblleft Oracle\textquotedblright\ $\Lambda_{m}\left(
\boldsymbol{\mu}\right)  =\lim_{t\rightarrow\infty}\varphi_{m}\left(
t,\boldsymbol{\mu}\right)  $. Further, if we can find a function
$K:\mathbb{R}^{2}\rightarrow\mathbb{R}$ that does not depend on any $\mu
\neq\mu_{0}$ and satisfies the Lebesgue-Stieltjes integral equation%
\begin{equation}
\psi\left(  t,\mu;\mu_{0}\right)  =\int K\left(  t,x;\mu_{0}\right)  dF_{\mu
}\left(  x\right)  , \label{eq3}%
\end{equation}
then the \textquotedblleft empirical phase function\textquotedblright\
\begin{equation*}
\hat{\varphi}_{m}\left(  t,\mathbf{z}\right)  =\dfrac{1}{m}\sum_{i=1}%
^{m}\left\{  1-K\left(  t,z_{i};\mu_{0}\right)  \right\}  \label{eqd2}%
\end{equation*}
satisfies $\mathbb{E}\left\{  \hat{\varphi}_{m}\left(  t,\mathbf{z}\right)
\right\}  =\varphi_{m}\left(  t,\boldsymbol{\mu}\right)  $ for any fixed $m,t$
and $\boldsymbol{\mu}$. Namely, $\hat{\varphi}_{m}\left(  t,\mathbf{z}\right)
$ is an unbiased estimator of $\varphi_{m}\left(  t,\boldsymbol{\mu}\right)
$. By the laws of large numbers, $\hat{\varphi}_{m}\left(  t,\mathbf{z}%
\right)$ can be close to $\varphi_{m}\left(  t,\boldsymbol{\mu}\right)$ for a fixed $t$ when $m$ is large. When the difference%
\begin{equation}
e_{m}\left(  t\right)  =\left\vert \hat{\varphi}_{m}\left(  t,\mathbf{z}%
\right)  -\varphi_{m}\left(  t,\boldsymbol{\mu}\right)  \right\vert
\label{eq2d}%
\end{equation}
is suitably small for large $t$, $\hat{\varphi}_{m}\left(  t,\mathbf{z}%
\right)  $ will accurately estimate $\pi_{1,m}$. Since $\varphi_{m}\left(
t,\boldsymbol{\mu}\right)  =\pi_{1,m}$ or $\hat{\varphi}_{m}\left(
t,\mathbf{z}\right)  =\pi_{1,m}$ rarely happens, $\hat{\varphi}_{m}\left(
t,\mathbf{z}\right)  $ usually employs a monotone increasing sequence
$\left\{  t_{m}\right\}  _{m\geq1}$ such that $\lim_{m\rightarrow\infty}%
t_{m}=\infty$ in order to achieve consistency, i.e., to achieve%
\begin{equation}
\label{defConsistency}
\Pr\left\{  \left\vert \frac{\hat{\varphi}_{m}\left(  t_{m},\mathbf{z}\right)
}{\pi_{1,m}}-1\right\vert \rightarrow0\right\}  \rightarrow1 \text{\quad as \quad} m \to \infty.
\end{equation}
So, the intrinsic speed for $\hat{\varphi}_{m}\left(  t_{m},\mathbf{z}\right)
$ to achieve consistency is better represented by $t_{m}$, and we will use
$t_{m}$ as the \textquotedblleft speed of convergence\textquotedblright\ of
$\hat{\varphi}_{m}\left(  t_{m},\mathbf{z}\right)  $. Throughout the paper, consistency of a proportion estimator
is defined via (\ref{defConsistency}) to accommodate the scenario $\lim_{m \to \infty}\pi_{1,m}=0$.

By duality, $\varphi_{m}^{\ast}\left(  t,\boldsymbol{\mu}\right)
=1-\varphi_{m}\left(  t,\boldsymbol{\mu}\right)  $ is the oracle for which
$\pi_{0,m}=\lim_{t\rightarrow\infty}\varphi_{m}^{\ast}\left(
t,\boldsymbol{\mu}\right)  $ for any fixed $m$ and $\boldsymbol{\mu}$,
$\hat{\varphi}_{m}^{\ast}=1-\hat{\varphi}_{m}\left(  t,\mathbf{z}\right)  $
satisfies $\mathbb{E}\left\{  \hat{\varphi}_{m}^{\ast}\left(  t,\mathbf{z}%
\right)  \right\}  =\varphi_{m}^{\ast}\left(  t,\boldsymbol{\mu}\right)  $ for
any fixed $m,t$ and $\boldsymbol{\mu}$, and $\hat{\varphi}_{m}^{\ast}\left(
t,\mathbf{z}\right)  $ will accurately estimate $\pi_{0,m}$ when $e_{m}\left(
t\right)  $ is suitably small for large $t$. Further, the stochastic
oscillations of $\hat{\varphi}_{m}^{\ast}\left(  t,\mathbf{z}\right)  $ and
$\hat{\varphi}_{m}\left(  t,\mathbf{z}\right)  $ are the same and is
quantified by $e_{m}\left(  t\right)  $.

We remark on the differences between the Strategy and Jin's Strategy. The
latter in our notations sets $\mu_{0}=0$, requires $\psi\left(  t,0;0\right)
=1$ for all $t$, requires location-shift families and uses Fourier transform
to construct $K$ and $\psi$, deals with distributions whose means are equal to
their medians, and intends to have%
\begin{equation}
1\geq\psi\left(  t,\mu;0\right)  \geq0\text{ \ for all }\mu\text{ and }t.
\label{eqd3}%
\end{equation}
It is not hard to see from the proof of Lemma 7.1 of \cite{Jin:2008} that
(\ref{eqd3}) may not be achievable for non-location-shift families. Further,
Jin's construction of $\psi$ is a special case of solving (\ref{eq3}).
Finally, it is easier to solve (\ref{eq3}) for $K:\mathbb{R}^{2}%
\rightarrow\mathbb{C}$ in the complex domain. However, a real-valued $K$ is
preferred for applications in statistics.

\section{Construction I}

\label{secEstimator}

In this section, we present the construction, referred to as \textquotedblleft
Construction I\textquotedblright, when the set of characteristic functions
(CFs) of the CDFs are of Riemann-Lebesgue type (see \autoref{Def}). This
construction essentially depends on generalizations of the Riemann-Lebesgue
Lemma (\textquotedblleft RL Lemma\textquotedblright) and subsumes that by
\cite{Jin:2008} for Gaussian family.

Recall the family of CDFs $\mathcal{F}=\left\{  F_{\mu}:\mu\in U\right\}  $
and let $\hat{F}_{\mu}\left(  t\right)  =\int e^{\iota tx}dF_{\mu}\left(
x\right)  $ be the CF of $F_{\mu}$ where $\iota=\sqrt{-1}$. Let $r_{\mu}$ be
the modulus of $\hat{F}_{\mu}$. Then $\hat{F}_{\mu}=r_{\mu}e^{\iota h_{\mu}}$,
where $h_{\mu}$ is the principal value of the argument that ranges in $(-\pi,\pi]$.
Further, $h_{\mu}$ restricted to the non-empty open interval
$\left(  -\tau_{\mu},\tau_{\mu}\right)  $ with $h_{\mu}\left(  0\right)  =0$
is uniquely defined, continuous and odd, where $\tau_{\mu}=\inf\left\{
t>0:\hat{F}_{\mu}\left(  t\right)  =0\right\}  >0$; see \cite{Luo:2004} for a
positive lower bound for $\tau_{\mu}$. Note that $\hat{F}_{\mu}$ has no real
zeros if and only if $\tau_{\mu}=\infty$.

\begin{definition}
\label{Def}If%
\begin{equation}
\left\{  t\in\mathbb{R}:\hat{F}_{\mu_{0}}\left(  t\right)  =0\right\}
=\varnothing\label{eq6b}%
\end{equation}
and, for each $\mu\in U\backslash\left\{  \mu_{0}\right\}  $%
\begin{equation}
\sup_{t\in\mathbb{R}}\frac{r_{\mu}\left(  t\right)  }{r_{\mu_{0}}\left(
t\right)  }<\infty\label{eq6}%
\end{equation}
and%
\begin{equation}
\lim_{t\rightarrow\infty}\frac{1}{t}\int_{\left[  -t,t\right]  }\frac{\hat
{F}_{\mu}\left(  y\right)  }{\hat{F}_{\mu_{0}}\left(  y\right)  }dy=0,
\label{eq6d}%
\end{equation}
then $\mathcal{\hat{F}}=\left\{  \hat{F}_{\mu}:\mu\in U\right\}  $ is said to
be of \textquotedblleft Riemann-Lebesgue type (RL type)\textquotedblright\ (at
$\mu_{0}$ on $U$).
\end{definition}

In \autoref{Def}, condition (\ref{eq6b}) requires that $\hat{F}_{\mu_{0}}$
have no real zeros, (\ref{eq6}) that $r_{\mu}$ with $\mu\in U\backslash
\left\{  \mu_{0}\right\}  $ approximately be of the same order as $r_{\mu_{0}%
}$, and (\ref{eq6d}) forces the \textquotedblleft mean value\textquotedblright%
\ of $\hat{F}_{\mu}\hat{F}_{\mu_{0}}^{-1}$ to converge to zero. Condition
(\ref{eq6b}) excludes the P\'{o}lya-type CF $\hat{F}\left(  t\right)  =\left\{
1-\left\vert t\right\vert \right\}  1_{\left\{  \left\vert t\right\vert
\leq1\right\}  }$ but holds for an infinitely divisible CF, (\ref{eq6})
usually cannot be relaxed to be $r_{\mu}r_{\mu_{0}}^{-1}\in L^{1}\left(
\mathbb{R}\right)  $ as seen from \autoref{lmLoc} for location-shift families,
and (\ref{eq6d}) induces a generalization of the RL Lemma as given by
\cite{Costin:2016} (which enables a construction via Fourier transform).

With a Fourier transform based construction comes the subtle issue of
determining an \textquotedblleft averaging function\textquotedblright%
\ $\omega$ that helps invoke the RL Lemma and facilitates easy numerical
implementation of the resulting proportion estimator. Indeed, a carefully
chosen $\omega$ will greatly simplify the construction and induce agreeable
proportion estimators. We adapt from \cite{Jin:2008} the concept of a
\textquotedblleft good\textquotedblright\ $\omega$:

\begin{definition}
\label{Def4}If a function $\omega:\left[  -1,1\right]  \rightarrow\mathbb{R}$
is non-negative and bounded such that $\int_{\left[  -1,1\right]  }%
\omega\left(  s\right)  ds=1$, then it is called \textquotedblleft
admissible\textquotedblright. If additionally $\omega$ is even on $\left[
-1,1\right]  $ and continuous on $\left(  -1,1\right)  $, then it is called
\textquotedblleft eligible\textquotedblright. If $\omega$ is eligible and
$\omega\left(  t\right)  \leq\tilde{\omega}\left(  1-t\right)  $ for all
$t\in\left(  0,1\right)  $ for some convex, super-additive function
$\tilde{\omega}$ over $\left(  0,1\right)  $, then it is called
\textquotedblleft good\textquotedblright.
\end{definition}

The definition above includes the end points $\left\{  -1,1\right\}  $ of the
compact interval $\left[  -1,1\right]  $ to tame $\omega$ at these points. For
example, the triangular density $\omega\left(  s\right)  =\max\left\{
1-\left\vert s\right\vert ,0\right\}  $ is good and its CF is%
\[
\hat{\omega}\left(  t\right)  =\frac{2\left\{  1-\cos\left(  t\right)  \right\}
}{t^{2}}1_{\left\{  t\neq0\right\}  }+1_{\left\{  t=0\right\}  }%
\]
as discussed by \cite{Jin:2008}. Since $r_{\mu}\left(  t\right)  $ is even and
uniformly continuous in $t$, an even $\omega$ matches $r_{\mu}$ and can induce
an even $K$ as a function of $t$. Unless otherwise noted, $\omega$ in this
work is always admissible.

\begin{theorem}
\label{ThmSeparable}Let $\mathcal{\hat{F}}$ be of RL type and define
$K:\mathbb{R}^{2}\rightarrow\mathbb{R}$ as%
\begin{equation}
K\left(  t,x;\mu_{0}\right)  =\int_{\left[  -1,1\right]  }\frac{\omega\left(
s\right)  \cos\left\{  tsx-h_{\mu_{0}}\left(  ts\right)  \right\}  }{r_{\mu_{0}%
}\left(  ts\right)  }ds. \label{eq2}%
\end{equation}
Then $\psi\left(  t,\mu;\mu_{0}\right)  $ in (\ref{eq3}) satisfies the following:

\begin{enumerate}
\item $\psi:\mathbb{R}\times U\rightarrow\mathbb{R}$ with%
\begin{equation}
\psi\left(  t,\mu;\mu_{0}\right)  =\int_{\left[  -1,1\right]  }\omega\left(
s\right)  \frac{r_{\mu}\left(  ts\right)  }{r_{\mu_{0}}\left(  ts\right)
}\cos\left\{  h_{\mu}\left(  ts\right)  -h_{\mu_{0}}\left(  ts\right)  \right\}
ds. \label{eq13}%
\end{equation}

\item $\psi\left(  t,\mu_{0};\mu_{0}\right)  =1$ for all $t$, and
$\lim_{t\rightarrow\infty}\psi\left(  t,\mu;\mu_{0}\right)  =0$ for each
$\mu\in U$ such that $\mu\neq\mu_{0}$.
\end{enumerate}
\end{theorem}

For random variables with RL type CFs in general, it is hard to ensure
(\ref{eqd3}), i.e., $1\geq\psi\left(  t,\mu;0\right)  \geq0$ for all $\mu$ and
$t$, since we do not have sufficient information on the phase $h_{\mu}$.
However, for certain location-shift families, (\ref{eqd3}) holds when $\omega$
is good; see \autoref{CorLocationShift}. Under slightly stronger conditions,
we have:

\begin{corollary}
\label{MethodICor}Assume that (\ref{eq6b}) and (\ref{eq6}) hold. If
$r_{\mu}/{r_{\mu_{0}}}\in L^{1}\left(  \mathbb{R}\right)  $, then the
conclusions of \autoref{ThmSeparable} hold.
\end{corollary}

When each $F_{\mu}$ has a density with respect to $\nu$, the condition
$r_{\mu}r_{\mu_{0}}^{-1}\in L^{1}\left(  \mathbb{R}\right)  \bigcap L^{\infty
}\left(  \mathbb{R}\right)  $ in \autoref{MethodICor} can be fairly strong
since it forces $\lim_{t\rightarrow\infty} r_{\mu}\left(  t\right)
/r_{\mu_{0}}\left(  t\right)  =0$. Unfortunately, \autoref{MethodICor} is
not applicable to location-shift families since for these families
${r_{\mu}\left(  t\right)  }/{r_{\mu_{0}}\left(  t\right)  }\equiv1$ for
all $\mu\in U$ whenever $r_{\mu_{0}}\left(  t\right)  \neq0$; see
\autoref{lmLoc}. In contrast, \autoref{ThmSeparable} is as we explain next.

Recall the definition of location-shift family, i.e., $\mathcal{F}=\left\{
F_{\mu}:\mu\in U\right\}  $ is a location-shift family if and only if
$z+\mu^{\prime}$ has CDF $F_{\mu+\mu^{\prime}}$ whenever $z$ has CDF $F_{\mu}$
for $\mu,\mu+\mu^{\prime}\in U$.

\begin{lemma}
\label{lmLoc}Suppose $\mathcal{F}$ is a location-shift family. Then, $r_{\mu}$
does not depend on $\mu$ and $h_{\mu}\left(  t\right)  =t\mu^{\prime}$, where
$\mu^{\prime}=\mu-\mu_{0}$. If in addition $\hat{F}_{\mu_{0}}\left(  t\right)
\neq0$ for all $t\in\mathbb{R}$, then $\mathcal{\hat{F}}$ is of RL type.
\end{lemma}

\autoref{lmLoc} implies that very likely the set of location-shift family
distributions is a subset of the set of distributions with RL type CFs.
However, verifying if $\hat{F}_{\mu_{0}}$ for a location-shift family has any
real zeros is quite difficult in general. Nonetheless, with \autoref{lmLoc},
\autoref{ThmSeparable} implies:

\begin{corollary}
\label{CorLocationShift}If $\mathcal{F}$ is a location-shift family for which
(\ref{eq6b}) holds, then%
\begin{equation}
\psi\left(  t,\mu;\mu_{0}\right)  =\int K\left(  t,y+  \mu-\mu
_{0}  ;\mu_{0}\right)  dF_{\mu_{0}}\left(  y\right)  =\int_{\left[
-1,1\right]  }\omega\left(  s\right)  \cos\left\{  ts\left(  \mu-\mu
_{0}\right)  \right\}  ds. \label{eq7a}%
\end{equation}
If additionally $\omega$ is good, then (\ref{eqd3}) holds, i.e., $1\geq
\psi\left(  t,\mu;\mu_{0}\right)  \geq0$ for all $\mu$ and $t$.
\end{corollary}

The identity (\ref{eq7a}) in \autoref{CorLocationShift} asserts that, for a
location-shift family, $\psi$ in (\ref{eq3}) reduces to a \textquotedblleft
convolution\textquotedblright, and it is the intrinsic mechanism behind the
Fourier transform based construction. When $\mu_{0}=0$, from
\autoref{CorLocationShift} we directly have%
\begin{equation*}
\psi\left(  t,\mu;0\right)  =\int K\left(  t,x+\mu;0\right)  dF_{0}\left(
x\right)  =\int_{\left[  -1,1\right]  }\omega\left(  s\right)  \cos\left(
ts\mu\right)  ds \label{eq7}%
\end{equation*}
recovering the construction by \cite{Jin:2008} for Gaussian and Laplace
families. Examples for which (\ref{eq7a}) holds are given in
\autoref{Sec:EgMeth1}.

\subsection{Construction I: some examples}

\label{Sec:EgMeth1}

We provide some examples from Construction I. They are all formed out of
location-shift families by allowing the location parameter $\mu$ to vary but
fixing the scale parameter $\sigma$, and are infinitely divisible so that none
of their CFs has any real zero; see \cite{Fischer:2014}, \cite{Lukacs:1970} and
\cite{Pitman:2003} for details on this. Note however that the Poisson family
as an example given in \autoref{SecEG2} is infinitely divisible but does not
have RL type CFs. There are certainly other examples that Construction I
applies to. But we will not attempt to exhaust them.

\begin{example}
Gaussian family $\mathcal{N}\left(  \mu,\sigma^{2}\right)  $ with mean
$\mu$ and standard deviation $\sigma>0$, for which%
\[
\frac{dF_{\mu}}{d\nu}\left(  x\right)  =f_{\mu}\left(  x\right)  =\left(
\sqrt{2\pi}\sigma\right)  ^{-1}\exp\left\{  -2^{-1}\sigma^{-2}\left(
x-\mu\right)  ^{2}\right\}  .
\]
The CF of $f_{\mu}$ is $\hat{f}_{\mu}\left(  t\right)  =\exp\left(  \iota
t\mu\right)  \exp\left(  -2^{-1}t^{2}\sigma^{2}\right)  $. Here $r_{\mu}%
^{-1}\left(  t\right)  =\exp\left(  2^{-1}t^{2}\sigma^{2}\right)  $.
Therefore,
\begin{equation}
K\left(  t,x;\mu_{0}\right)  =\int_{\left[  -1,1\right]  }r_{\mu}^{-1}\left(
ts\right)  \omega\left(  s\right)  \exp\left\{  \iota ts\left(  x-\mu
_{0}\right)  \right\}  ds \label{eqb14}%
\end{equation}
and (\ref{eq7a}) holds. When $\mu_{0}=0$ and $\sigma=1$, this construction
reduces to that in Section 2.1 of \cite{Jin:2008}.
\end{example}

\begin{example}
Laplace family $\mathsf{Laplace}\left(  \mu,2\sigma^{2}\right)  $ with mean
$\mu$ and standard deviation $\sqrt{2}\sigma>0$ for which%
\[
\frac{dF_{\mu}}{d\nu}\left(  x\right)  =f_{\mu}\left(  x\right)  =\frac
{1}{2\sigma}\exp\left(  -\sigma^{-1}\left\vert x-\mu\right\vert \right)
\]
and the CF is $\hat{f}_{\mu}\left(  t\right)  =\left(  1+\sigma^{2}%
t^{2}\right)  ^{-1}\exp\left(  \iota t\mu\right)  $. Therefore, (\ref{eqb14})
and (\ref{eq7a}) hold with $r_{\mu}^{-1}\left(  t\right)  =1+\sigma^{2}t^{2}$.
When $\mu_{0}=0$ and $\sigma=1$, this construction reduces to that in Section
7.1 of \cite{Jin:2008}.
\end{example}

\begin{example}
Logistic family $\mathsf{Logistic}\left(  \mu,\sigma\right)  $ with mean $\mu$
and scale parameter $\sigma>0$, for which%
\[
\frac{dF_{\mu}}{d\nu}\left(  x\right)  =f_{\mu}\left(  x\right)  =\frac
{1}{4\sigma}\operatorname{sech}^{2}\left(  \frac{x-\mu}{2\sigma}\right)
\]
and the CF of $f_{\mu}$ is $\hat{f}_{\mu}\left(  t\right)  =\pi\sigma t\left\{
\sinh\left(  \pi\sigma t\right)  \right\}  ^{-1}\exp\left(  \iota t\mu\right)
$. Therefore, (\ref{eqb14}) and (\ref{eq7a}) hold with $r_{\mu}^{-1}\left(
t\right)  =\left(  \pi\sigma t\right)  ^{-1}\sinh\left(  \pi\sigma t\right)
\sim\left(  2\pi\sigma t\right)  ^{-1}e^{\pi\sigma t}$ as $t\rightarrow\infty$.
\end{example}

\begin{example}
Cauchy family $\mathsf{Cauchy}\left(  \mu,\sigma\right)  $ with median $\mu$
and scale parameter $\sigma>0$, for which%
\[
\frac{dF_{\mu}}{d\nu}\left(  x\right)  =f_{\mu}\left(  x\right)  =\frac{1}%
{\pi\sigma}\frac{\sigma^{2}}{\left(  x-\mu\right)  ^{2}+\sigma^{2}}%
\]
and the CF is $\hat{f}_{\mu}\left(  t\right)  =\exp\left(  -\sigma\left\vert
t\right\vert \right)  \exp\left(  \iota t\mu\right)  $. Therefore,
(\ref{eqb14}) and (\ref{eq7a}) hold with $r_{\mu}^{-1}\left(  t\right)
=\exp\left(  \sigma\left\vert t\right\vert \right)  $. Note that the Cauchy
family does not have first-order absolute moment.
\end{example}

\begin{example}
The Hyperbolic Secant family $\mathsf{HSecant}\left(  \mu,\sigma\right)  $
with mean $\mu$ and scale parameter $\sigma>0$, for which%
\[
\frac{dF_{\mu}}{d\nu}\left(  x\right)  =f_{\mu}\left(  x\right)  =\frac
{1}{2\sigma}\frac{1}{\cosh\left(  \pi\frac{x-\mu}{\sigma}\right)  };
\]
see, e.g., Chapter 1 of \cite{Fischer:2014}. The identity%
\[
\int_{-\infty}^{+\infty}e^{\iota tx}\frac{dx}{\pi\cosh\left(  x\right)
}=\cosh\left(  2^{-1}\pi t\right)  ,
\]
implies $\hat{F}_{\mu}\left(  t\right)  =\sigma^{-1}\exp\left(  -\iota
t\mu\sigma^{-1}\right)  \operatorname{sech}\left(  t\sigma^{-1}\right)  $.
Therefore, (\ref{eqb14}) and (\ref{eq7a}) hold with $r_{\mu}^{-1}\left(
t\right)  =\sigma\cosh\left(  t\sigma^{-1}\right)  \sim2^{-1}\sigma\exp\left(
\sigma^{-1}t\right)  $ as $t\rightarrow\infty$. The Hyperbolic Secant family has been
used to model stock indices and exchange rates \citep{Fischer:2014} or status of coronary heart disease
\citep{Vaughan:2002}.
\end{example}

\subsection{Construction I: uniform consistency and speed of convergence}

\label{SecConsistencyI}

The performance of the estimator $\hat{\varphi}_{m}$ depends on how accurately
it approximates $\varphi_{m}$, the oracle that knows the true value $\pi
_{1,m}$ as $t\rightarrow\infty$. Specifically, the smaller $e_{m}\left(
t\right)  $ defined by (\ref{eq2d}) is when $t$ is large, the more accurately
$\hat{\varphi}_{m}\left(  t;\mathbf{z}\right)  $ estimates $\pi_{1,m}$. Two
key factors that affect $e_{m}$ are: (i) the magnitude of the reciprocal of
the modulus, $r_{\mu_{0}}^{-1}\left(  t\right)  $, which appears as a scaling
factor in the integrand in the definition of $K$ in (\ref{eq2}), and (ii) the
magnitudes of the $\mu_{i}$'s and the variabilities of $z_{i}$'s. As will be
shown by \autoref{ThmConsistency}, for independent $\left\{  z_{i}\right\}
_{i=1}^{m}$, the oscillation of $e_{m}\left(  t\right)  $ depends mainly on
$r_{\mu_{0}}^{-1}\left(  t\right)  $ due to concentration of measure for
independent, uniformly bounded random variables and their transforms by
Lipschitz functions, whereas the consistency of $\hat{\varphi}_{m}\left(
t;\mathbf{z}\right)  $ depends also on the magnitudes of $\pi_{1,m}$ and
$\mu_{i}$'s that affect how accurate the oracle is when $t$ is large.

\begin{theorem}
\label{ThmConsistency}If $\mathcal{\hat{F}}$ is of RL type and $\left\{
z_{i}\right\}  _{i=1}^{m}$ are independent. Let $a\left(  t;\mu_{0}\right)
=\int_{\left[  -1,1\right]  }{r_{\mu_{0}}^{-1}\left(  ts\right)  }ds$ for
$t\in\mathbb{R}$. Then%
\begin{equation}
\mathbb{V}\left\{  \left\vert \hat{\varphi}_{m}\left(  t,\mathbf{z}\right)
-\varphi_{m}\left(  t,\boldsymbol{\mu}\right)  \right\vert \right\}  \leq
\frac{1}{m}\left\Vert \omega\right\Vert _{\infty}^{2}a^{2}\left(  t;\mu
_{0}\right)  . \label{eq2g}%
\end{equation}
Further, for any fixed $\lambda>0$, with probability at least $1-2\exp\left(
-2^{-1}\lambda^{2}\right)  $,%
\begin{equation}
\left\vert \hat{\varphi}_{m}\left(  t,\mathbf{z}\right)  -\varphi_{m}\left(
t,\boldsymbol{\mu}\right)  \right\vert \leq\frac{\lambda\left\Vert
\omega\right\Vert _{\infty}a\left(  t;\mu_{0}\right)  }{\sqrt{m}}.
\label{eq2f}%
\end{equation}
If there are positive sequences $\left\{  u_{m}\right\}  _{m\geq1}$, $\left\{
\lambda_{m}\right\}  _{m\geq1}$ and $\left\{  t_{m}\right\}  _{m\geq1}$ such
that%
\begin{equation}
\lim_{m\rightarrow\infty}\sup\left\{  \psi\left(  t,\mu;\mu_{0}\right)
:\left(  t,\left\vert \mu\right\vert \right)  \in\lbrack t_{m},\infty
)\times\lbrack u_{m},\infty)\right\}  =0 \label{eq4}%
\end{equation}
and%
\begin{equation}
\lim_{m\rightarrow\infty}\frac{\lambda_{m}a\left(  t_{m};\mu_{0}\right)  }%
{\pi_{1,m}\sqrt{m}}=0\text{ \ and \ }\lim_{m\rightarrow\infty}\exp\left(
-2^{-1}\lambda_{m}^{2}\right)  =0, \label{eq2h}%
\end{equation}
then%
\begin{equation}
\Pr\left\{  \left\vert \pi_{1,m}^{-1}\hat{\varphi}_{m}\left(  t_{m}%
,\mathbf{z}\right)  -1\right\vert \rightarrow0\right\}  \rightarrow1
\label{eq4a}%
\end{equation}
whenever $\left\{  \mu_{i}:i\in I_{1,m}\right\}  \subseteq\lbrack u_{m}%
,\infty)$.
\end{theorem}

\autoref{ThmConsistency} bounds the variance of $\hat{\varphi}_{m}\left(
t,\mathbf{z}\right)  $, captures the key ingredients needed for and the
essence of proving the consistency of a proportion estimator based on
Construction I, and shows that such a proportion estimator is consistent as
long as $\lambda_{m}\rightarrow\infty$, $\lambda_{m}\pi_{1,m}^{-1}a\left(
t_{m};\mu_{0}\right)  $ is of smaller order than $\sqrt{m}$ and each
$\psi\left(  t_{m},\mu_{i};\mu_{0}\right)  $ with $i\in I_{1,m}$ is negligible
for a sequence $t_{m}\rightarrow\infty$. Using \autoref{ThmConsistency}, we
can characterize the consistency of the estimator $\hat{\varphi}_{m}\left(
t,\mathbf{z}\right)  $ for each example given in \autoref{Sec:EgMeth1} as follows:

\begin{corollary}
\label{CorConI}Consider $\hat{\varphi}\left(  t_{m},\mathbf{z}\right)  $ from
Construction I. Let $\left\{  z_{i}\right\}  _{i=1}^{m}$ be independent and
$\left\{  t_{m}:m\geq1\right\}  $ a positive sequence such that $\lim
_{m\rightarrow\infty}t_{m}=\infty$. Assume $\min_{i\in I_{1,m}}\left\vert
\mu_{i}-\mu_{0}\right\vert \geq\left(  t_{m}\right)  ^{-1}\ln\ln m$. Then
$\Pr\left\{  \left\vert \pi_{1,m}^{-1}\hat{\varphi}_{m}\left(  t_{m}%
,\mathbf{z}\right)  -1\right\vert \rightarrow0\right\}  \rightarrow1$ holds
\begin{itemize}
  \item for Gaussian, Hyperbolic Secant, Logistic and Cauchy family respectively when
$t_{m}=\sigma^{-1}\sqrt{2\gamma\ln m}$, $t_{m}=\sigma\gamma\ln m$,
$t_{m}=\left(  \sigma\pi\right)  ^{-1}\gamma\ln m$ and $t_{m}=\sigma
^{-1}\gamma\ln m$, $\lambda
_{m}=o\left(  t_{m}\right)  $, $\lambda_{m}\rightarrow\infty$ and $\pi
_{1,m}\geq Cm^{\gamma-0.5}$ with $\gamma\in\left(  0,0.5\right]  $;

 \item for Laplace family when $t_{m}=\ln m$,
$\lambda_{m}=O\left(  t_{m}\right)  $, $\lambda_{m}\rightarrow\infty$ and
$\pi_{1,m}\geq Cm^{-\gamma}$ with $\gamma\in\left[  0,0.5\right)  $,
\end{itemize}
 where
$C>0$ can be any constant for which $\pi_{1,m}\in\left(  0,1\right]  $ as
$\gamma$ varies in its designated range respectively.

\end{corollary}

To characterize if an estimator $\hat{\varphi}_{m}\left(  t,\mathbf{z}\right)
$ is uniformly consistent with respect to $\pi_{1,m}$ and $t$, we introduce
the following definition:

\begin{definition}
\label{DefUniformConsistency}Given a family $\mathcal{F}$, the sequence of
sets $\mathcal{Q}_{m}\left(  \boldsymbol{\mu},t;\mathcal{F}\right)
\subseteq\mathbb{R}^{m}\times\mathbb{R}\ $for each $m\in\mathbb{N}_{+}$ is
called a \textquotedblleft uniform consistency class\textquotedblright\ for
the estimator $\hat{\varphi}_{m}\left(  t,\mathbf{z}\right)  $ if%
\begin{equation}
\Pr\left\{  \sup\nolimits_{\boldsymbol{\mu}\mathcal{\in Q}_{m}\left(
\boldsymbol{\mu},t;\mathcal{F}\right)  }\left\vert \pi_{1,m}^{-1}%
\sup\nolimits_{t\in\mathcal{Q}_{m}\left(  \boldsymbol{\mu},t;\mathcal{F}%
\right)  }\hat{\varphi}_{m}\left(  t,\mathbf{z}\right)  -1\right\vert
\rightarrow0\right\}  \rightarrow1. \label{defUCC}%
\end{equation}
If (\ref{defUCC}) holds and the $t$-section of $\mathcal{Q}_{m}\left(  \boldsymbol{\mu
},t;\mathcal{F}\right)  $ (that is a subset of $\mathbb{R}^m$ containing $\boldsymbol{\mu}$) does not converge to the empty set in $\mathbb{R}^{\mathbb{N}}$ as $m \to \infty$, then $\hat{\varphi}_{m}\left(  t,\mathbf{z}\right)  $
is said to be \textquotedblleft uniformly consistent\textquotedblright. If
further the $\boldsymbol{\mu}$-section of $\mathcal{Q}_{m}\left(  \boldsymbol{\mu},t;\mathcal{F}\right)  $
(that is a subset of $\mathbb{R}$ containing $t$) contains a connected subset $G_{m}\subseteq$ $\mathbb{R}$ such that
$\lim_{m\rightarrow\infty}\nu\left(  G_{m}\right)  =\infty$, then
$\hat{\varphi}_{m}\left(  t,\mathbf{z}\right)  $ is said to be
\textquotedblleft uniformly consistent in frequency domain\textquotedblright.
\end{definition}

Now we discuss uniform consistency in frequency domain of an estimator from
Construction I. Define%
\[
\mathcal{B}_{m}\left(  \rho\right)  =\left\{  \boldsymbol{\mu}\in
\mathbb{R}^{m}:m^{-1}\sum\nolimits_{i=1}^{m}\left\vert \mu_{i}-\mu
_{0}\right\vert \leq\rho\right\}  \text{ \ for some }\rho>0
\]
and $u_{m}=\min\left\{  \left\vert \mu_{j}-\mu_{0}\right\vert :\mu_{j}\neq
\mu_{0}\right\}  $.

\begin{theorem}
\label{ThmLocationShift}Assume that $\mathcal{F}$ is a location-shift family
for which (\ref{eq6b}) holds and $\int\left\vert x\right\vert ^{2}dF_{\mu
}\left(  x\right)  <\infty$ for each $\mu\in U$. If $\sup_{y\in\mathbb{R}%
}\left\vert \frac{d}{dy}h_{\mu_{0}}\left(  y\right)  \right\vert =C_{\mu_{0}%
}<\infty$, then for the estimator $\hat{\varphi}_{m}\left(  t,\mathbf{z}%
\right)  $ from Construction I, a uniform consistency class is
\begin{equation*}
\mathcal{Q}_{m}\left(  \boldsymbol{\mu},t;\mathcal{F}\right)  =\left\{
\begin{array}
[c]{c}%
q\gamma^{\prime}>\vartheta>2^{-1},\gamma^{\prime}>0,\gamma^{\prime\prime
}>0,0\leq\vartheta^{\prime}<\vartheta-1/2,\\
R_{m}\left(  \rho\right)  =O\left(  m^{\vartheta^{\prime}}\right)  ,\tau
_{m}\leq\gamma_{m},u_{m}\geq\frac{\ln\ln m}{\gamma^{\prime\prime}\tau_{m}%
},\\
t\in\left[  0,\tau_{m}\right]  ,\lim\limits_{m\rightarrow\infty}\pi_{1,m}%
^{-1}\Upsilon\left(  q,\tau_{m},\gamma_{m},r_{\mu_{0}}\right)  =0
\end{array}
\right\}  \label{MainClass}%
\end{equation*}
where $q$, $\gamma^{\prime}$, $\gamma^{\prime\prime}$, $\vartheta$ and
$\vartheta^{\prime}$ are constants, $R_{m}\left(  \rho\right)  =2\int%
{\left\vert x\right\vert dF_{\mu_{0}}\left(  x\right)  }+2\rho+2C_{\mu_{0}}$,
$\gamma_{m}=\gamma^{\prime}\ln m$ and%
\begin{equation*}
\Upsilon\left(  q,\tau_{m},\gamma_{m},r_{\mu_{0}}\right)  =\frac{2\left\Vert
\omega\right\Vert _{\infty}\sqrt{2q\gamma_{m}}}{\sqrt{m}}\sup_{t\in\left[
0,\tau_{m}\right]  }\int_{\left[  0,1\right]  }\frac{ds}{r_{\mu_{0}}\left(
ts\right)  }. \label{eq12g}%
\end{equation*}
Moreover, for all sufficiently large $m$,
\begin{equation}
\sup_{\boldsymbol{\mu}\in\mathcal{B}_{m}\left(  \rho\right)  }\sup
_{t\in\left[  0,\tau_{m}\right]  }\left\vert \hat{\varphi}_{m}\left(
t,\mathbf{z}\right)  -\varphi_{m}\left(  t,\boldsymbol{\mu}\right)
\right\vert \leq\Upsilon\left(  q,\tau_{m},\gamma_{m},r_{\mu_{0}}\right)
\label{eq11d}%
\end{equation}
holds with probability at least $1-o\left(  1\right)  $.
\end{theorem}

Several remarks on \autoref{ThmLocationShift} are ready to be stated. Firstly,
since \autoref{ThmLocationShift} requires the random variables to have finite
absolute second-order moments, it may not apply to location-shift families
that do not have first-order absolute moments. Secondly, even though
\autoref{ThmLocationShift} potentially allows for many possible choices of $t$
for $\hat{\varphi}_{m}\left(  t,\mathbf{z}\right)  $, we should choose
$\tau_{m}$ such that $\tau_{m}\rightarrow\infty$ as fast as possible so that
$\pi_{1,m}^{-1}\hat{\varphi}_{m}\left(  t_{m},\mathbf{z}\right)  \rightarrow1$
as fast as possible. Thirdly, compared to Theorems 1.4 and 1.5 of
\cite{Jin:2008} where $\mathcal{B}_{m}\left(  \rho\right)  $ is for a fixed
$\rho$ for Gaussian family, we allow $\rho\rightarrow\infty$ for
location-shift families. Fourthly,\ the bound (\ref{eq11d}) together with
(\ref{eq12g}) imply that, when other things are kept fixed, the larger
$r_{\mu}^{-1}$ is, the slower $\pi_{1,m}^{-1}\hat{\varphi}_{m}\left(
t_{m},\mathbf{z}\right)  \rightarrow1$. This has been observed by
\cite{Jin:2008} for the Gaussian and Laplace families since $r_{\mu}%
^{-1}\left(  t\right)  $ for the former is much larger than the latter when
$t$ is large. Fifthly, for location-shift families the constants $q$,
$\gamma^{\prime}$, $\gamma^{\prime\prime}$ and $\vartheta$ can be specified by
a user. If additionally $\vartheta^{\prime}=0$, then $\mathcal{Q}_{m}\left(
\boldsymbol{\mu},t;\mathcal{F}\right)  $ is fully data-adaptive and depends
only on $u_{m}=\min_{i\in I_{1,m}}\left\vert \mu_{i}-\mu_{0}\right\vert $; see
also \autoref{CorInclusion1} on this for examples given in
\autoref{Sec:EgMeth1}.

Let $\mathcal{F}$ and $\mathcal{\tilde{F}}$ be two location-shift families
that have RL type CFs and are determined by the same set of parameters, and
$r_{\mu_{0}}$ and $\tilde{r}_{\mu_{0}}$ respectively the moduli of the CFs of
$F_{\mu_{0}}\in\mathcal{F}$ and $\tilde{F}_{\mu_{0}}\in\mathcal{\tilde{F}}$.
For the estimator $\hat{\varphi}_{m}\left(  t,\mathbf{z}\right)  $ under the
same settings and under the conditions stated in \autoref{ThmLocationShift},
consider the two uniform consistency classes\ $\mathcal{Q}_{m}\left(
\boldsymbol{\mu},t;\mathcal{F}\right)  $ and $\mathcal{Q}_{m}\left(
\boldsymbol{\mu},t;\mathcal{\tilde{F}}\right)  $ for $\hat{\varphi}_{m}\left(
t,\mathbf{z}\right)  $ that have the same constants $q$, $\gamma^{\prime}$,
$\vartheta$ and $\gamma^{\prime\prime}$. Then \autoref{ThmLocationShift}
implies that%
\begin{equation}
r_{\mu_{0}}\left(  t\right)  \leq\tilde{r}_{\mu_{0}}\left(  t\right)  \text{
for all large }t>0\text{ }\Longrightarrow\text{ }\mathcal{Q}_{m}\left(
\boldsymbol{\mu},t;\mathcal{F}\right)  \supseteq\mathcal{Q}_{m}\left(
\boldsymbol{\mu},t;\mathcal{\tilde{F}}\right)  \label{eqInc1}%
\end{equation}
In other words, for location-shift families that have RL type CFs, the uniform
consistency classes for estimators from Construction I\ can be reversely
ordered via set inclusion by the magnitudes of the associated moduli when
other things are kept fixed.

Another consequence of \autoref{ThmLocationShift} is as follows. If we fix
$q,\gamma^{\prime},\gamma^{\prime\prime},\vartheta$ and $\vartheta^{\prime}$,
then $r_{\mu_{0}}$ determines $\mathcal{Q}_{m}\left(  \boldsymbol{\mu
},t;\mathcal{F}\right)  $. In particular, if $F_{\mu}$ has a density with
respect to $\nu$, then $\lim_{t\rightarrow\infty}r_{\mu}\left(  t\right)  =0$
must hold, which forces (\ref{eq12g}) to give%
\begin{equation}
\Upsilon\left(  q,\tau_{m},\gamma^{\prime}\ln m,r_{\mu_{0}}\right)  \geq
\frac{C\sqrt{q\gamma^{\prime}\ln m}}{\sqrt{m}}\text{ \ as \ }m\rightarrow
\infty. \label{eq12g2}%
\end{equation}
Since Hoeffding inequality, used in the proof of \autoref{ThmLocationShift},
is asymptotically optimal for independent, uniformly almost surely bounded random variables (see \cite{Cohen:1999})
induced by Construction I, \autoref{ThmLocationShift} and (\ref{eq12g2})
together imply that, when Construction I applied to location-shift families
with absolutely continuous CDFs, a uniform consistency class for the
corresponding proportion estimator is unlikely able to contain any $\pi
_{1,m}\in\left(  0,m^{-0.5}\right]  $. Comparing this with the conclusions of
\autoref{ThmConsistency} and \autoref{CorConI}, we see that a sacrifice to
achieve uniform consistency in frequency domain is the reduction of the range
for $\pi_{1,m}$ for which $\hat{\varphi}_{m}\left(  t,\mathbf{z}\right)  $ can
be consistent.

Finally, \autoref{ThmLocationShift} gives

\begin{corollary}
\label{CorInclusion1}When $\left\{  z_{i}\right\}  _{i=1}^{m}$ are
independent, the following hold for $\hat{\varphi}_{m}\left(  t,\mathbf{z}%
\right)  $:

\begin{enumerate}
\item For Gaussian family: $q\sigma^{-1}>\vartheta>2^{-1},0\leq\vartheta
^{\prime}<\vartheta-1/2,\gamma_{m}=\sigma^{-1}\ln m$ and
\begin{equation}
\mathcal{Q}_{m}\left(  \boldsymbol{\mu},t;\mathcal{F}\right)  =\left\{
\begin{array}
[c]{c}%
R_{m}\left(  \rho\right)  =O\left(  m^{\vartheta^{\prime}}\right)  ,u_{m}%
\geq\frac{\ln\ln m}{\sqrt{2\ln m}},0<\gamma\leq0.5\\
\pi_{1,m}\geq Cm^{\gamma-0.5},t\in\left[  0,\sigma^{-1}\sqrt{2\gamma\ln
m}\right]
\end{array}
\right\}  \label{defB2}%
\end{equation}
The fastest speed of convergence is $\sqrt{\ln m}$, achieved when
$\liminf_{m\rightarrow\infty}\pi_{1,m}>0$.

\item For Laplace family: $q>\vartheta>2^{-1},0\leq\vartheta^{\prime
}<\vartheta-1/2,\gamma_{m}=\ln m$ and%
\[
\mathcal{Q}_{m}\left(  \boldsymbol{\mu},t;\mathcal{F}\right)  =\left\{
\begin{array}
[c]{c}%
R_{m}\left(  \rho\right)  =O\left(  m^{\vartheta^{\prime}}\right)  ,u_{m}%
\geq\frac{\ln\ln m}{\ln m},0\leq\gamma<1/2\\
\pi_{1,m}\geq Cm^{-\gamma},t\in\left[  0,\ln m\right]
\end{array}
\right\}  .
\]

\item Hyperbolic Secant family: $q\sigma>\vartheta>2^{-1}$, $\gamma_{m}%
=\sigma\ln m$, $0\leq\vartheta^{\prime}<\vartheta-1/2$ and
\[
\mathcal{Q}_{m}\left(  \boldsymbol{\mu},t;\mathcal{F}\right)  =\left\{
\begin{array}
[c]{c}%
R_{m}\left(  \rho\right)  =O\left(  m^{\vartheta^{\prime}}\right)  ,u_{m}%
\geq\frac{\ln\ln m}{\ln m},0<\gamma<1/2\\
\pi_{1,m}\geq Cm^{\gamma-0.5},t\in\left[  0,\sigma\gamma\ln m\right]
\end{array}
\right\}  .
\]

\item Logistic family: $q\left(  \sigma\pi\right)  ^{-1}>\vartheta>2^{-1}$,
$\gamma_{m}=\left(  \sigma\pi\right)  ^{-1}\ln m$, $0\leq\vartheta^{\prime
}<\vartheta-1/2$ and
\[
\mathcal{Q}_{m}\left(  \boldsymbol{\mu},t;\mathcal{F}\right)  =\left\{
\begin{array}
[c]{c}%
R_{m}\left(  \rho\right)  =O\left(  m^{\vartheta^{\prime}}\right)  ,u_{m}%
\geq\frac{\ln\ln m}{\ln m},0<\gamma<1/2\\
\pi_{1,m}\geq Cm^{\gamma-0.5},t\in\left[  0,\left(  \sigma\pi\right)
^{-1}\gamma\ln m\right]
\end{array}
\right\}  .
\]

\end{enumerate}

In each case above, $C>0$ can be any constant for which $\pi_{1,m}\in\left(
0,1\right]  $ as $\gamma$ varies in its designated range.
\end{corollary}

\autoref{CorInclusion1} provides uniform consistency classes for estimators
from Construction I when it is applied to five location-shift families with
general scale parameter $\sigma>0$. In particular, if $\sigma=1$ and $q>3/2$,
$\vartheta=q/3$ and $\vartheta^{\prime}=0$ is set in (\ref{defB2}), then we
recover the uniform consistency class given by Theorems 1.4 and 1.5 of
\cite{Jin:2008}.

\section{Construction II and Construction III}

\label{SecConstructionII}

When the CDFs of $\left\{  z_{i}\right\}  _{i=1}^{m}$ do not have RL type CFs,
Construction I in \autoref{secEstimator} cannot be used. In particular,
outside location-shift families, the translation-convolution equivalence does
not hold and Hoeffding inequality is no longer applicable. This makes the
construction of a proportion estimator using the Strategy much more
challenging. So, we will restrict our attention to $\left\{  z_{i}\right\}
_{i=1}^{m}$ whose CDFs do not have RL type CFs but belong to NEFs whose mean
and variance are functionally related. Specifically, we show that the Strategy
is implementable for discrete NEFs with infinite supports or continuous NEFs with \textquotedblleft
separable moment functions\textquotedblright\ (see \autoref{Def3}). These
include 8 of the total of 12 NEF-CVFs. The techniques of construction mainly
use generating functions (GFs) and Mellin transform.

\subsection{A brief review on natural exponential families}

\label{sec:onNEFs}

We provide a very brief review on NEF, whose details can be found in
\cite{Letac:1992}. Let $\beta$ be a positive Radon measure on $\mathbb{R}$
that is not concentrated on one point. Let $L\left(  \theta\right)  =\int
e^{x\theta}\beta\left(  dx\right)  $ for $\theta\in\mathbb{R}$ be its Laplace
transform and $\Theta$ be the maximal open set containing $\theta$ such that
$L\left(  \theta\right)  <\infty$. Suppose $\Theta$ is not empty and let
$\kappa\left(  \theta\right)  =\ln L\left(  \theta\right)  $ be the cumulant
function of $\beta$. Then%
\[
\mathcal{F}=\left\{  G_{\theta}:G_{\theta}\left(  dx\right)  =\exp\left\{
\theta x-\kappa\left(  \theta\right)  \right\}  \beta\left(  dx\right)
,\theta\in\Theta\right\}
\]
forms an NEF with respect to the basis $\beta$. Note that $\Theta$ has a
non-empty interior if it is not empty and that $L$ is analytic on the strip
$A_{\Theta}=\left\{  z\in\mathbb{C}:\Re\left(  z\right)  \in\Theta\right\}  $.

The NEF $\mathcal{F}$ can be equivalently characterized by its mean domain and
variance function. Specifically, the mean function $\mu:\Theta\rightarrow U$
with $U=\mu\left(  \Theta\right)  $ is given by $\mu\left(  \theta\right)
=\frac{d}{d\theta}\kappa\left(  \theta\right)  $, and the variance function is
$V\left(  \theta\right)  =\frac{d^{2}}{d\theta^{2}}\kappa\left(
\theta\right)  $ and can be parametrized by $\mu$ as
\[
V\left(  \mu\right)  =\int\left(  x-\mu\right)  ^{2}F_{\mu}\left(  dx\right)
\text{ for }\mu\in U,
\]
where $\theta=\theta\left(  \mu\right)  $ is the inverse function of $\mu$ and
$F_{\mu}=G_{\theta\left(  \mu\right)  }$. Namely, $\mathcal{F}=\left\{
F_{\mu}:\mu\in U\right\}  $. The pair $\left(  V,U\right)  $ is called the
variance function of $\mathcal{F}$, and it characterizes $\mathcal{F}$.

For the constructions of proportion estimators for NEFs, we will reuse the
notation $K$ but take $K$ as a function of $t$ and $\theta\in\Theta$. Note
that $K$ depends on $\theta_{0}$ but not on any $\theta\neq\theta_{0}$.
Further, we will reuse the notation $\psi$ but take it as a function of $t$
and $\theta\in\Theta$. For an NEF, $\psi$ defined by (\ref{eq3}) becomes%
\[
\psi\left(  t,\theta;\theta_{0}\right)  =\int K\left(  t,x;\theta_{0}\right)
dG_{\theta}\left(  x\right)  \text{ for }G_{\theta}\in\mathcal{F}.
\]
Let $\boldsymbol{\theta}=\left(  \theta_{1},\ldots,\theta_{m}\right)  $. Then
accordingly $\hat{\varphi}_{m}\left(  t,\mathbf{z}\right)  =\dfrac{1}{m}%
\sum_{i=1}^{m}\left\{  1-K\left(  t,z_{i};\theta_{0}\right)  \right\}  $ and
$\varphi_{m}\left(  t,\boldsymbol{\theta}\right)  =\dfrac{1}{m}\sum_{i=1}%
^{m}\left\{  1-\psi\left(  t,\theta_{i};\theta_{0}\right)  \right\}  $.

\subsection{Construction II: discrete NEFs with infinite supports}

\label{secConstII}

Suppose the basis $\beta$ for $\mathcal{F}$ is discrete with support
$\mathbb{N}$, i.e., there exists a positive sequence $\left\{  c_{k}\right\}
_{k\geq0}$ such that
\begin{equation}
\beta=\sum\nolimits_{k=0}^{\infty}c_{k}\delta_{k}. \label{eq12a}%
\end{equation}
Then the power series $H\left(  z\right)  =\sum_{k=0}^{\infty}c_{k}z^{k}$ with
$z\in\mathbb{C}$ must have a positive radius of convergence $R_{H}$, and $H$
is the generating function (GF) of $\beta$. Further, if $\beta$ is a
probability measure, then $(-\infty,0]\subseteq\Theta$ and $R_{H}\geq1$, and
vice versa. The following approach, which we refer to as \textquotedblleft
Construction II\textquotedblright, provides the construction for discrete NEFs
with support $\mathbb{N}$.

\begin{theorem}
\label{ThmDiscreteInfiniteSupport}Let $\mathcal{F}$ be the NEF generated by
$\beta$ in (\ref{eq12a}) and $\omega$ admissible. For $x\in\mathbb{N}$ and
$t\in\mathbb{R}$ set
\begin{equation}
K\left(  t,x;\theta_{0}\right)  =H\left(  e^{\theta_{0}}\right)  \int_{\left[
-1,1\right]  }\frac{\left(  ts\right)  ^{x}\cos\left(  \frac{\pi x}%
{2}-tse^{\theta_{0}}\right)  }{H^{\left(  x\right)  }\left(  0\right)  }%
\omega\left(  s\right)  ds. \label{II-c}%
\end{equation}
Then%
\[
\psi\left(  t,\theta;\theta_{0}\right)  =\int K\left(  t,x;\theta_{0}\right)
dG_{\theta}\left(  x\right)  =\frac{H\left(  e^{\theta_{0}}\right)  }{H\left(
e^{\theta}\right)  }\int_{\left[  -1,1\right]  }\cos\left\{  st\left(
e^{\theta}-e^{\theta_{0}}\right)  \right\}  \omega\left(  s\right)  ds,
\]
$\psi\left(  t,\theta_{0};\theta_{0}\right)  =1$ for any $t$, and
$\lim_{t\rightarrow\infty}$ $\psi\left(  t,\theta;\theta_{0}\right)  =0$ for
each $\theta\neq\theta_{0}$.
\end{theorem}

In \autoref{ThmDiscreteInfiniteSupport}, $H^{\left(  k\right)  }\left(
0\right)  =c_{k}k!$ for $k\in\mathbb{N}$. So, for Construction II, if $\beta$
is known, then we can use $c_{k}k!$ instead of $H^{\left(  k\right)  }\left(
0\right)  $, whereas if $H$ is known and $H^{\left(  k\right)  }$ is easy to
compute, we can use $H^{\left(  k\right)  }\left(  0\right)  $. This will
greatly aid the numerical implementation of Construction II.

\subsection{Construction II: some examples}

\label{SecEG2}

\autoref{ThmDiscreteInfiniteSupport} covers the construction for Abel,
Negative Binomial, Poisson, Strict Arcsine, Large Arcsine and Tak\'{a}cs
families, each of which is an NEF-CVF, has basis $\beta$ with support
$\mathbb{N}$, is infinite divisible such that $H\left(  e^{\theta}\right)
\neq0$ for each $\theta\in\Theta$, and has non-RL type CFs; see
\cite{Letac:1990} for details on these distributions. However, for each of
Abel and Large Arcsine families, the corresponding GF is a composition of two
analytic functions, and manually computing $H\left(  e^{\theta_{0}}\right)  $
in the construction of $K$ in the statement of
\autoref{ThmDiscreteInfiniteSupport} may be cumbersome.

\begin{example}
Poisson family $\mathcal{P}\left(  \mu\right)  $ with mean $\mu>0$, for
which $F_{\mu}\left(  \left\{  k\right\}  \right)  =\Pr\left(  X=k\right)
=\frac{\mu^{k}}{k!}e^{-\mu}$ for $k\in\mathbb{N}$. Clearly,
\[
\hat{F}_{\mu}\left(  t\right)  =\exp\left\{  \mu\left(  e^{\iota t}-1\right)
\right\}  =\exp\left\{  \mu\left(  \cos t-1\right)  \right\}  \exp\left(
\iota\mu\sin t\right)  .
\]
However, $\left\{  \hat{F}_{\mu}:\mu>0\right\}  $ are not of RL type since
when $\mu=\mu_{0}+1$ and $t>0$,%
\begin{align*}
\frac{1}{t}\Re\left(  \int_{\left[  -t,t\right]  }\frac{\hat{F}_{\mu}\left(
y\right)  }{\hat{F}_{\mu_{0}}\left(  y\right)  }dy\right)   &  =\frac{e^{-1}%
}{t}\int_{\left[  -t,t\right]  }\exp\left(  \cos y\right)  \cos\left(  \sin
y\right)  dy\\
&  \geq2e^{-1}e^{-1}\cos1>0.
\end{align*}
The basis is $\beta=\sum_{k=0}^{\infty}(k!)^{-1}\delta_{k}$, $L\left(
\theta\right)  =\exp\left(  e^{\theta}\right)  $ with $\theta\in\mathbb{R}$,
$\mu\left(  \theta\right)  =e^{\theta}$, $H\left(  z\right)  =e^{z}$ with
$R_{H}=\infty$ and $H^{\left(  k\right)  }\left(  0\right)  =1$ for all
$k\in\mathbb{N}$. The Poisson family has been used to model RNA-Seq data \citep{Chen:2014intra}.
\end{example}

\begin{example}
\label{EgNegBin}Negative Binomial family $\mathsf{NegBinomial}\left(
\theta,n\right)  $ with $\theta<0$ and $n\in\mathbb{N}_{+}$ such that%
\[
G_{\theta}\left(  \left\{  k\right\}  \right)  =\Pr\left(  X=k\right)
=\frac{c_{k}^{\ast}}{k!}e^{k\theta}\left(  1-e^{\theta}\right)  ^{n}%
\]
with $c_{k}^{\ast}={\left(  k+n-1\right)  !}/{\left(  n-1\right)  !}$
for\ $k\in\mathbb{N}$. The basis is $\beta=\sum_{k=0}^{\infty}
{c_{k}^{\ast}}(k!)^{-1}\delta_{k}$, $L\left(  \theta\right)  =\left(  1-e^{\theta
}\right)  ^{-n}$ with $\theta<0$, $\mu\left(  \theta\right)  =ne^{\theta
}\left(  1-e^{\theta}\right)  ^{-1}$, $H\left(  z\right)  =\left(  1-z\right)
^{-n}$ with $R_{H}=1$, and $H^{\left(  k\right)  }\left(  0\right)
=c_{k}^{\ast}$ for all $k\in\mathbb{N}$. The Negative Binomial family has also been used to model RNA-Seq data \citep{Robinson:2008,Di:2011}.
\end{example}

\begin{example}
Strict Arcsine family. Its VF is $V\left(  u\right)  =u\left(  1+u^{2}\right)
$ and $\beta=\sum_{n=0}^{\infty}
(n!)^{-1}{c_{n}^{\ast}\left(  1\right)  }\delta_{n}$, where $c_{0}^{\ast}\left(  1\right)  =c_{1}^{\ast}\left(
1\right)  =1,$%
\begin{equation}
c_{2n}^{\ast}\left(  \sigma\right)  =\prod_{k=0}^{n-1}\left(  \sigma
^{2}+4k^{2}\right)  \text{ \ and \ }c_{2n+1}^{\ast}\left(  \sigma\right)
=\sigma\prod_{k=0}^{n-1}\left(  \sigma^{2}+\left(  2k+1\right)  ^{2}\right)
\label{eqPcoef}%
\end{equation}
for $\sigma>0$ and $n\in\mathbb{N}_{+}$. Further, $H\left(  z\right)
=\exp\left(  \arcsin z\right)  $ with $R_{H}=1$. The Strict Arcsine family has been used to model insurance claims \citep{Kokonendji:2004}.
\end{example}

\begin{example}
Large Arcsine family. Its VF is $V\left(  u\right)  =u\left(  1+2u+2u^{2}%
\right)  $ and $\beta=\sum\nolimits_{n=0}^{\infty}c_{n}\delta_{n}$ with
$c_{n}={c_{n}^{\ast}\left(  1+n\right)  }/{\left(  n+1\right)  !}$ for
$n\in\mathbb{N}$, where $c_{n}^{\ast}\left(  \sigma\right)  $ is defined in
(\ref{eqPcoef}) and for which $H\left(  z\right)  =\exp\left\{  \arcsin\left(
h\left(  z\right)  \right)  \right\}  $ with $h\left(  z\right)  =\sum
_{k=0}^{\infty}c_{k}z^{k+1}$. It can be seen that $R_{H}$ must be finite;
otherwise, $\lim_{\left\vert z_{l}\right\vert \rightarrow\infty}\left\vert
h\left(  z_{l}\right)  \right\vert =\infty$ for a sequence $\left\{
z_{l}:l\geq1\right\}  $, and $H\left(  z_{l}\right)  $ cannot be expanded into
a convergent power series at $z_{l}$ for $l$ sufficiently large.
\end{example}

\begin{example}
Abel family. Its VF $V\left(  u\right)  =u\left(  1+u\right)  ^{2}$,
$\beta=\sum\nolimits_{k=0}^{\infty}c_{k}\delta_{k}$ with $c_{k}={\left(
1+k\right)  ^{k-1}}/{k!}$ for $k\in\mathbb{N}$, and $H\left(  z\right)
=e^{h\left(  z\right)  }$ with $h\left(  z\right)  =\sum_{k=0}^{\infty}%
c_{k}z^{k+1}$ with $R_{H}=e^{-1}$. The Abel family has been used to model birds' migration patterns and other phenomena \citep{Nandi:1994}.
\end{example}

\begin{example}
Tak\'{a}cs family. Its VF is $V\left(  u\right)  =u\left(  1+u\right)  \left(
1+2u\right)  $ and $\beta=\delta_{0}+\sum\nolimits_{k=1}^{\infty}c_{k}%
\delta_{k}$ with $c_{k}= {\left(  2k\right)  !}/({k!\left(  k+1\right)  !})$
for $k\in\mathbb{N}_{+}$, for which $H\left(  z\right)  =\left({1-\sqrt{1-4z}%
}\right)/{2z}$ with $R_{H}=4^{-1}$ and $z=0$ is a removable singularity of $H$.
\end{example}

The above calculations show that the GFs of Negative Binomial, Strict Arcsine,
Large Arcsine, Abel and Tak\'{a}cs families all have positive and finite radii
of convergence whereas that of Poisson family has infinite radius of
convergence. This will be very helpful in determining uniform consistency
classes for Construction II and its numerical implementation; see
\autoref{CoroII}.

\subsection{Construction III: continuous NEFs with separable moments}

\label{secConstructionIII}

In contrast to NEFs with support $\mathbb{N}$, we consider non-location-shift
NEFs whose members are continuous distributions. Assume $0\in\Theta$, so that
$\beta$ is a probability measure with finite moments of all orders. Let%
\begin{equation}
\tilde{c}_{n}\left(  \theta\right)  =\frac{1}{L\left(  \theta\right)  }\int
x^{n}e^{\theta x}\beta\left(  dx\right)  =\int x^{n}dG_{\theta}\left(
x\right)  \text{ \ for \ }n\in\mathbb{N} \label{eq11c}%
\end{equation}
be the moment sequence for $G_{\theta}\in\mathcal{F}$. Note that (\ref{eq11c})
is the Mellin transform of the measure $G_{\theta}$.

\begin{definition}
\label{Def3}If there exist two functions $\zeta,\xi:\Theta\rightarrow
\mathbb{R\ }$and a sequence $\left\{  \tilde{a}_{n}\right\}  _{n\geq0}$ that
satisfy the following:

\begin{itemize}
\item $\xi$ $\left(  \theta\right)  \neq\xi\left(  \theta_{0}\right)  $
whenever $\theta\neq\theta_{0}$, $\zeta\left(  \theta\right)  \neq0$ for all
$\theta\in U$, and $\zeta$ does not depend on any $n\in\mathbb{N}$,

\item $\tilde{c}_{n}\left(  \theta\right)  =\xi^{n}\left(  \theta\right)
\zeta\left(  \theta\right)  \tilde{a}_{n}$ for each $n\in\mathbb{N}$ and
$\theta\in\Theta$,

\item $\Psi\left(  t,\theta\right)  =\sum_{n=0}^{\infty}\frac{t^{n}\xi
^{n}\left(  \theta\right)  }{\tilde{a}_{n}n!}$ is absolutely convergent
pointwise in $\left(  t,\theta\right)  \in\mathbb{R}\times\Theta$,
\end{itemize}

\noindent then the moment sequence $\left\{  \tilde{c}_{n}\left(
\theta\right)  \right\}  _{n\geq0}$ is called \textquotedblleft
separable\textquotedblright\ (at $\theta_{0}$).
\end{definition}

The concept of separable moment sequence is an analogy to the structured
integrand used in (\ref{eq13g}) for Construction II, and the condition on
$\Psi\left(  t,\theta\right)  $ is usually satisfied since $\sum_{n=0}%
^{\infty}\frac{t^{n}\xi^{n}\left(  \theta\right)  }{n!}$ already is convergent
pointwise on $\mathbb{R}\times\Theta$. The next approach, which we refer to as
\textquotedblleft Construction III\textquotedblright, is based on Mellin
transform of $G_{\theta}$ and applies to NEFs with separable moment sequences.

\begin{theorem}
\label{ThmConstructionMoments}Assume that the NEF $\mathcal{F}$ has a
separable moment sequence $\left\{  \tilde{c}_{n}\left(  \theta\right)
\right\}  _{n\geq0}$ at $\theta_{0}$, and let $\omega$ be admissible. For
$t,x\in\mathbb{R}$ set%
\begin{equation*}
K\left(  t,x;\mu_{0}\right)  =\frac{1}{\zeta\left(  \theta_{0}\right)  }%
\int_{\left[  -1,1\right]  }\sum_{n=0}^{\infty}\frac{\left(  -tsx\right)
^{n}\cos\left\{  \frac{\pi}{2}n+ts\xi\left(  \theta_{0}\right)  \right\}
}{\tilde{a}_{n}n!}\omega\left(  s\right)  ds. \label{III-a}%
\end{equation*}
Then%
\begin{equation*}
\psi\left(  t,\mu;\mu_{0}\right)  =\int K\left(  t,x;\theta_{0}\right)
dG_{\theta}\left(  x\right)  =\frac{\zeta\left(  \theta\right)  }{\zeta\left(
\theta_{0}\right)  }\int_{\left[  -1,1\right]  }\cos\left[  ts\left\{
\xi\left(  \theta_{0}\right)  -\xi\left(  \theta\right)  \right\}  \right]
\omega\left(  s\right)  ds, \label{eq1g}%
\end{equation*}
$\psi\left(  t,\theta_{0};\theta_{0}\right)  =1$ for any $t$ and
$\lim_{t\rightarrow\infty}$ $\psi\left(  t,\theta;\theta_{0}\right)  =0$ for
each $\theta\neq\theta_{0}$.
\end{theorem}

Compared to Constructions I and II, Construction III involves the integral of
an infinite series and is more complicated. However, it deals with NEFs that
have more complicated structures than the former two.

\subsection{Construction III: two examples}

\label{SecEG3}

We provide two examples from Construction III for Exponential and Gamma
families, respectively. Note that Gamma family contains Exponential family and central
Chi-square family.

\begin{example}
Exponential family $\mathsf{Exponential}\left(  \mu\right)  $ with mean
$\mu>0$ and basis $\beta\left(  dx\right)  =e^{-x}\nu\left(  dx\right)  $, for
which $L\left(  \theta\right)  =\left(  1-\theta\right)  ^{-1}$ and
$\mu\left(  \theta\right)  =1-\theta$ for $\theta<1$. Further,%
\[
\frac{dF_{\mu}}{d\nu}\left(  x\right)  =f_{\mu}\left(  x\right)  =\mu e^{-\mu
x}1_{[0,\infty)}\left(  x\right)
\]
and $\tilde{c}_{n}\left(  \mu\right)  =\left(\mu{n!}\right)/{\mu^{n}}$. So, $L\left(
\mu\right)  =\mu^{-1}$, $\xi\left(  \mu\right)  =\mu^{-1}$, $\tilde{a}_{n}=n!$
and $\zeta\equiv1$. Setting%
\begin{equation*}
K\left(  t,x;\mu_{0}\right)  =\int_{\left[  -1,1\right]  }\omega\left(
s\right)  \sum_{n=0}^{\infty}\frac{\left(  -tsx\right)  ^{n}\cos\left(
\frac{\pi}{2}n+\frac{ts}{\mu_{0}}\right)  }{\left(  n!\right)  ^{2}}ds
\label{eq1a}%
\end{equation*}
gives%
\[
\psi\left(  t,\mu;\mu_{0}\right)  =\int K\left(  t,x;\mu_{0}\right)  dF_{\mu
}\left(  x\right)  =\int_{\left[  -1,1\right]  }\cos\left\{  ts\left(  \mu
^{-1}-\mu_{0}^{-1}\right)  \right\}  \omega\left(  s\right)  ds.
\]

\end{example}

\begin{example}
\label{EgGamma}Gamma family $\mathsf{Gamma}\left(  \theta,\sigma\right)  $
with basis $\beta$ such that%
\[
\frac{d\beta}{d\nu}\left(  x\right)  =\frac{x^{\sigma-1}e^{-x}}{\Gamma\left(
\sigma\right)  }1_{\left(  0,\infty\right)  }\left(  x\right)  dx\text{ \ with
}\sigma>0,
\]
where $\Gamma$ is the Euler's Gamma function. So, $L\left(  \theta\right)
=\left(  1-\theta\right)  ^{-\sigma}$,%
\begin{equation}
\frac{dG_{\theta}}{d\nu}\left(  x\right)  =f_{\theta}\left(  x\right)
=\left(  1-\theta\right)  ^{\sigma}\frac{e^{\theta x}x^{\sigma-1}e^{-x}%
}{\Gamma\left(  \sigma\right)  }1_{\left(  0,\infty\right)  }\left(  x\right)
\label{eq1h}%
\end{equation}
for $\theta<1$, and $\mu\left(  \theta\right)  ={\sigma}/({1-\theta})$.
Since%
\[
\tilde{c}_{n}\left(  \theta\right)  =\left(  1-\theta\right)  ^{\sigma}%
\int_{0}^{\infty}\frac{e^{-y}y^{n+\sigma-1}}{\Gamma\left(  \sigma\right)
}\frac{dy}{\left(  1-\theta\right)  ^{n+\sigma-1}}=\frac{\Gamma\left(
n+\sigma\right)  }{\Gamma\left(  \sigma\right)  }\frac{1}{\left(
1-\theta\right)  ^{n}},
\]
we see $\xi\left(  \theta\right)  =\left(  1-\theta\right)  ^{-1}$, $\tilde
{a}_{n}={\Gamma\left(  n+\sigma\right)  }/{\Gamma\left(  \sigma\right)  }$
and $\zeta\equiv1$. Setting%
\[
K\left(  t,x;\theta_{0}\right)  =\int_{\left[  -1,1\right]  }\omega\left(
s\right)  \sum_{n=0}^{\infty}\frac{\left(  -tsx\right)  ^{n}\Gamma\left(
\sigma\right)  \cos\left(  \frac{\pi}{2}n+\frac{ts}{1-\theta_{0}}\right)
}{n!\Gamma\left(  \sigma+n\right)  }ds,
\]
we obtain%
\[
\psi\left(  t,\theta;\theta_{0}\right)  =\int_{\left[  -1,1\right]  }%
\cos\left[  ts\left\{  \left(  1-\theta_{0}\right)  ^{-1}-\left(
1-\theta\right)  ^{-1}\right\}  \right]  \omega\left(  s\right)  ds,
\]
for which $\psi\left(  t,\theta_{0};\theta_{0}\right)  =1$ for all $t$ and
$\lim_{t\rightarrow\infty}$ $\psi\left(  t,\mu;\theta_{0}\right)  =0$ for each
$\theta\neq\theta_{0}$.

Recall (\ref{eq7a}) of Construction I based on Fourier transform, i.e.,%
\[
\psi\left(  t,\mu;\mu_{0}\right)  =\int K\left(  t,x;\mu_{0}\right)  dF_{\mu
}\left(  x\right)  =\int K\left(  t,y+   \mu-\mu_{0}   ;\mu
_{0}\right)  dF_{\mu_{0}}\left(  y\right)  ,
\]
where the action of Fourier transform is seen as the translation $K\left(
t,x;\mu_{0}\right)  \mapsto K\left(  t,y+   \mu-\mu_{0}   ;\mu
_{0}\right)  $. In contrast, the action of Mellin transform is seen via%
\begin{align}
\psi\left(  t,\theta;\theta_{0}\right)   &  =\int K\left(  t,x;\theta
_{0}\right)  dG_{\theta}\left(  x\right) \nonumber\\
&  =\int_{0}^{\infty}K\left(  t,x;\theta_{0}\right)  \left(  1-\theta\right)
^{\sigma}\frac{e^{-\left(  1-\theta\right)  x}x^{\sigma-1}}{\Gamma\left(
\sigma\right)  }dx\label{eq1d}\\
&  =\int K\left(  t,\frac{y}{1-\theta};\theta_{0}\right)  \beta\left(
dy\right)  , \label{eq1b}%
\end{align}
where from (\ref{eq1d}) to (\ref{eq1b}) scaling $K\left(  t,x;\theta
_{0}\right)  \mapsto K\left(  t,y \left(  1-\theta\right)^{-1}  ;\theta
_{0}\right)  $ is induced. This comparison clearly shows the action of Mellin
transform as the multiplication-convolution equivalence, in contrast to the
action of Fourier transform as the translation-convolution equivalence.
\end{example}

\section{Two non-existence results for the Strategy}

\label{SecNonexist}

In this section, we provide two example families for which the Strategy is not
implementable. To state them, we introduce

\begin{definition}
The proposition \textquotedblleft PropK\textquotedblright: there exists a
$K:\mathbb{R}^{2}\rightarrow\mathbb{R}$ such that $K$ does not depend on any
$\theta\in\Theta_{1}$ with $\theta\neq\theta_{0}$ and that%
\[
\psi\left(  t,\theta;\theta_{0}\right)  =\int K\left(  t,x;\theta_{0}\right)
dG_{\theta}\left(  x\right)
\]
satisfies $\lim_{t\rightarrow\infty}\psi\left(  t,\theta_{0};\theta
_{0}\right)  =1$ and $\lim_{t\rightarrow\infty}\psi\left(  t,\theta;\theta
_{0}\right)  =0$ for $\theta\in\Theta_{1}$ with $\theta\neq\theta_{0}$, where
$\Theta_{1}$ is a subset of $\Theta$ that has a non-empty interior and does
not contain $\theta_{0}$.
\end{definition}

\noindent For the following two families, i.e., Inverse Gaussian and Binomial
families, PropK does not hold. Note that the former family is discrete but has
a finite support, whereas the latter is continuous but does not have a
separable moment sequence.

\begin{example}
The Inverse Gaussian family $\mathsf{InvGaussian}\left(  \theta,\sigma\right)
$ with scale parameter $\sigma>0$ and basis%
\[
\beta\left(  dx\right)  =x^{-3/2}\exp\left(  -\frac{\sigma^{2}}{2x}\right)
\frac{\sigma}{\sqrt{2\pi}}1_{\left(  0,\infty\right)  }\left(  x\right)
\nu\left(  dx\right)  ,
\]
for which $L\left(  \theta\right)  =\exp\left(  -\sigma\sqrt{-2\theta}\right)
$ with $\theta<0$. So,%
\[
\frac{dG_{\theta}}{d\nu}\left(  x\right)  =f_{\theta}\left(  x\right)
=\frac{x^{-3/2}}{L\left(  \theta\right)  }\exp\left(  -\frac{\sigma^{2}}%
{2x}\right)  \frac{\sigma}{\sqrt{2\pi}}1_{\left(  0,\infty\right)  }\left(
x\right)  .
\]
By a change of variables $y=\frac{\sigma^{2}}{2x}$, we obtain%
\begin{align}
\int K\left(  t,x;\theta_{0}\right)  dG_{\theta}\left(  x\right)   &
=\frac{1}{L\left(  \theta\right)  }\frac{\sigma}{\sqrt{2\pi}}\int_{0}^{\infty
}K\left(  t,x;\theta_{0}\right)  x^{-3/2}\exp\left(  -\frac{\sigma^{2}}%
{2x}\right)  dx\nonumber\\
&  =\frac{1}{\sqrt{\pi}L\left(  \theta\right)  }\int_{0}^{\infty}K\left(
t,\frac{\sigma^{2}}{2y};\theta_{0}\right)  y^{-1/2}e^{-y}dy. \label{eq11f}%
\end{align}
Since the integral on the right hand side of (\ref{eq11f}) is not a function
of $\theta$ for $\theta\neq\theta_{0}$, PropK does not hold and $K$ does not
exist. Note that the Gaussian family and Inverse Gaussian family are
reciprocal pairs, called so by \cite{Letac:1990}. The Inverse Gaussian family has been used to model the shelf life of products \citep{Folks:1978}.
\end{example}

\begin{example}
\label{egBin}Binomial family $\mathsf{Binomial}\left(  \theta,n\right)  $ such
that%
\[
\mathbb{P}\left(  X=k\right)  =\binom{n}{k}\theta^{k}\left(  1-\theta\right)
^{n-k}\text{ \ for \ \ }\theta\in\left(  0,1\right)  .
\]
We will show that if PropK holds, then $K$ must be a function of $\theta
\in\Theta_{1}$, a contradiction. Assume PropK holds. Then,
\[
\psi\left(  t,\theta;\theta_{0}\right)  =\int K\left(  t,x;\theta_{0}\right)
dG_{\theta}\left(  x\right)  =\left(  1-\theta\right)  ^{n}\sum_{k=0}^{n}%
a_{k}\left(  t;\theta_{0}\right)  \binom{n}{k}\frac{\theta^{k}}{\left(
1-\theta\right)  ^{k}},
\]
where $a_{x}\left(  t;\theta_{0}\right)  =K\left(  t,x;\theta_{0}\right)  $
for $x=0,...,n$.

Let $d_{k}\left(  t;\theta_{0}\right)  =a_{k}\left(  t;\theta_{0}\right)
\binom{n}{k}$ and $q\left(  \theta\right)  ={\theta}/({1-\theta})$. Pick
$n+1$ distinct values $\varpi_{i},i=0,\ldots,n$ from $\left(  0,1\right)  $
such that $\varpi_{0}=\theta_{0}$ and $\varpi_{i}\in\Theta_{1}$ for
$i=1,\ldots,n$. Further, define the $\left(  n+1\right)  \times\left(
n+1\right)  $ Vandermonde matrix $\mathbf{V}_{n+1}$ whose $\left(  i,j\right)
$ entry is $q^{j-1}\left(  \varpi_{i-1}\right)  $, i.e., $\mathbf{V}%
_{n+1}=\left(  q^{j-1}\left(  \varpi_{i-1}\right)  \right)  $, $\mathbf{d}%
_{t}=\left(  d_{0}\left(  t;\theta_{0}\right)  ,\ldots,d_{n}\left(
t;\theta_{0}\right)  \right)  ^{\top}$ and $\mathbf{b}=\left(  \left(
1-\theta_{0}\right)  ^{-n},0,\ldots,0\right)  ^{\top}$. Then the determinant
$\left\vert \mathbf{V}_{n+1}\right\vert \neq0$, and the properties of $K$
imply%
\[
\lim_{t\rightarrow\infty}\mathbf{V}_{n+1}\mathbf{d}_{t}=\mathbf{b}.
\]
However, $\mathbf{V}_{n+1}:\mathbb{R}^{n+1}\rightarrow\mathbb{R}^{n+1}$ as a
bounded linear mapping is a homeomorphism with the bounded inverse
$\mathbf{V}_{n+1}^{-1}$. So,%
\begin{equation*}
\lim_{t\rightarrow\infty}\mathbf{d}_{t}=\lim_{t\rightarrow\infty}%
\mathbf{V}_{n+1}^{-1}\left(  \mathbf{V}_{n+1}\left(  \mathbf{d}_{t}\right)
\right)  =\mathbf{V}_{n+1}^{-1}\mathbf{b.} \label{eq11i}%
\end{equation*}
By Cramer's rule, we obtain
\begin{equation}
\lim_{t\rightarrow\infty}d_{n}\left(  t;\theta_{0}\right)  =\frac{\left(
-1\right)  ^{n}\left\vert \mathbf{V}_{n}\right\vert }{\left\vert
\mathbf{V}_{n+1}\right\vert \left(  1-\theta_{0}\right)  ^{n}}=\frac{\left(
-1\right)  ^{n}\left(  1-\theta_{0}\right)  ^{-n}}{\prod\nolimits_{k=1}%
^{n}\left(  q\left(  \varpi_{k}\right)  -q\left(  \theta_{0}\right)  \right)
}, \label{eq11h}%
\end{equation}
where $\mathbf{V}_{n}$ is the submatrix of $\mathbf{V}_{n+1}$ obtained by
removing the first row and last column of $\mathbf{V}_{n+1}\,$. But
(\ref{eq11h}) is a contradiction since $d_{n}\left(  t;\theta_{0}\right)  $
does not depend on any $\theta\in\Theta_{1}$ for all $t$. For the case $n=1$,
we easily see from (\ref{eq11h}) the contradiction%
\[
\lim_{t\rightarrow\infty}a_{1}\left(  t;\theta_{0}\right)  =\frac{-\left(
1-\theta_{0}\right)  ^{-1}}{q\left(  \theta\right)  -q\left(  \theta
_{0}\right)  }\text{ \ for each \ }\theta\in\Theta_{1}.
\]
To summarize, PropK does not hold and $K$ does not exist.
\end{example}

\section{Construction II: consistency and speed of convergence}

\label{SecConsistent2}

Recall $H^{\left(  k\right)  }\left(  0\right)  =c_{k}k!$ for $k\in\mathbb{N}%
$. Call the sequence $\left\{  \left(  c_{k}k!\right)  ^{-1},k\in
\mathbb{N}\right\}  $ the \textquotedblleft reciprocal derivative sequence (of
$H$ at $0$)\textquotedblright. The following lemma gives bounds on the
magnitudes of this sequence for the examples given in \autoref{SecEG2}. It
will help derive concentration inequalities for estimators from Construction II.

\begin{lemma}
\label{LmNEFDiscrete}Consider the examples given in \autoref{SecEG2}. Then
$c_{k}k!\equiv1$ for Poisson family, whereas for Negative Binomial, Abel and
Tak\'{a}cs families with a fixed $n$ and $\sigma>0$,%
\begin{equation}
\frac{1}{c_{k}k!}\leq\frac{C}{k!}\text{ \ for all }k\in\mathbb{N}.
\label{eqUPde}%
\end{equation}
However, for Strict Arcsine and Large Arcsine families, both with a fixed
$\sigma>0$, (\ref{eqUPde}) does not hold. On the other hand, for any
$\tilde{r}>0$ smaller than the radius $R_{H}$ of convergence of $H$,
\begin{equation}
\frac{1}{c_{k}k!}=\frac{1}{H^{\left(  k\right)  }\left(  0\right)  }\geq
\frac{1}{H\left(  \tilde{r}\right)  }\frac{\tilde{r}^{k}}{k!}\text{ \ \ \ for
all \ }k\in\mathbb{N}. \label{eqLBder}%
\end{equation}

\end{lemma}

\autoref{LmNEFDiscrete} shows that, among the six discrete NEFs with support
$\mathbb{N}$ given in \autoref{SecEG2}, the reciprocal derivative sequence for
Poisson family has the largest magnitude, whereas this sequence for Negative
Binomial, Abel and Tak\'{a}cs families with a fixed $n$ and $\sigma>0$ are all
dominated by the \textquotedblleft reciprocal factorial
sequence\textquotedblright\ $\left\{  {1}/{k!}:k\in\mathbb{N}\right\}  $
approximately. Further, \autoref{LmNEFDiscrete} asserts that the reciprocal
derivative sequence dominates the \textquotedblleft exponential
sequence\textquotedblright\ $\left\{  H\left(  \tilde{r}\right)
{k!}/{\tilde{r}^{k}}:k\in\mathbb{N}\right\}  $ for any $\tilde{r}>0$ smaller
than the radius of convergence of $H$.

Let $\eta=e^{\theta}$ for $\theta\in\Theta$, $\eta_{i}=e^{\theta_{i}}$ for
$0\leq i\leq m$ and $\boldsymbol{\eta}=\left(  \eta_{1},\ldots,\eta
_{m}\right)  $. First, we provide upper bounds on the variance and
oscillations of $\hat{\varphi}_{m}\left(  t,\mathbf{z}\right)  -\varphi
_{m}\left(  t,\boldsymbol{\theta}\right)  $ when $t$ is positive and
sufficiently large.

\begin{theorem}
\label{II-concentration}Let $\mathcal{F}$ be an NEF generated by $\beta$ in
(\ref{eq12a}), $\left\{  z_{i}\right\}  _{i=1}^{m}$ independent with CDFs
$\left\{  G_{\theta_{i}}\right\}  _{i=1}^{m}$ belonging to $\mathcal{F}$,
$\lambda$ a positive constant, and $t$ positive and sufficiently large.

\begin{enumerate}
\item Let $\phi_{m}\left(  L,\boldsymbol{\theta}\right)  =\min_{1\leq i\leq
m}\left\{  L\left(  \theta_{i}\right)  \eta_{i}^{1/4}\right\}  $. If
(\ref{eqUPde}) holds, then%
\[
\mathbb{V}\left\{  \hat{\varphi}_{m}\left(  t,\mathbf{z}\right)  -\varphi
_{m}\left(  t,\boldsymbol{\theta}\right)  \right\}  \leq V_{\mathrm{II}%
,1}^{\left(  m\right)  }=\frac{C}{m}\frac{\exp\left(  2t\left\Vert
\boldsymbol{\eta}\right\Vert _{\infty}^{1/2}\right)  }{\sqrt{t}\phi_{m}\left(
L,\boldsymbol{\theta}\right)  }%
\]
and%
\begin{equation}
\Pr\left\{  \left\vert \hat{\varphi}_{m}\left(  t,\mathbf{z}\right)
-\varphi_{m}\left(  t,\boldsymbol{\theta}\right)  \right\vert \geq
\lambda\right\}  \leq\lambda^{-2}V_{\mathrm{II},1}^{\left(  m\right)  }.
\label{eq15e4}%
\end{equation}

\item Let $L_{\min}^{\left(  m\right)  }=\min_{1\leq i\leq m}L\left(
\theta_{i}\right)  $. Then for Poisson family,%
\begin{equation}
\mathbb{V}\left\{  \hat{\varphi}_{m}\left(  t,\mathbf{z}\right)  -\varphi
_{m}\left(  t,\boldsymbol{\theta}\right)  \right\}  \leq V_{\mathrm{II}%
,2}^{\left(  m\right)  }=\frac{C}{m}\frac{\exp\left(  t^{2}\left\Vert
\boldsymbol{\eta}\right\Vert _{\infty}^{1/2}\right)  }{L_{\min}^{\left(
m\right)  }} \label{eqII3}%
\end{equation}
and%
\[
\Pr\left\{  \left\vert \hat{\varphi}_{m}\left(  t,\mathbf{z}\right)
-\varphi_{m}\left(  t,\boldsymbol{\theta}\right)  \right\vert \geq
\lambda\right\}  \leq\lambda^{-2}V_{\mathrm{II},2}^{\left(  m\right)  }.
\]

\end{enumerate}
\end{theorem}

We remark that the assertion in \autoref{II-concentration} on Poisson family
holds for any NEF with support $\mathbb{N}$ such that $c_{k}k!\leq C$ for all
$k\in\mathbb{N}$. With \autoref{II-concentration}, we derive the uniform
consistency classes and speeds of convergence for the estimators from
Construction II. Recall $\eta=e^{\theta}$ for $\theta\in\Theta$ and
$\boldsymbol{\eta}=\left(  e^{\theta_{1}},\ldots,e^{\theta_{m}}\right)  $.

\begin{theorem}
\label{II-consistency}Let $\mathcal{F}$ be the NEF generated by $\beta$ in
(\ref{eq12a}), $\left\{  z_{i}\right\}  _{i=1}^{m}$ independent with CDFs
$\left\{  G_{\theta_{i}}\right\}  _{i=1}^{m}$ belonging to $\mathcal{F}$, and
$\rho$ a finite, positive constant.

\begin{enumerate}
\item If (\ref{eqUPde}) holds, then a uniform consistency class is%
\[
\mathcal{Q}_{\mathrm{II},1}\left(  \boldsymbol{\theta},t,\pi_{1,m}%
;\gamma\right)  =\left\{
\begin{array}
[c]{c}%
\left\Vert \boldsymbol{\theta}\right\Vert _{\infty}\leq\rho,\pi_{1,m}\geq
m^{\left(  \gamma-1\right)  /2},t=2^{-1}\left\Vert \boldsymbol{\eta
}\right\Vert _{\infty}^{-1/2}\gamma\ln m,\\
\lim_{m\rightarrow\infty}t\min_{i\in I_{1,m}}\left\vert \eta_{0}-\eta
_{i}\right\vert =\infty
\end{array}
\right\}
\]
for any fixed $\gamma\in\left(  0,1\right]  .$ The speed of convergence is poly-log.

\item For Poisson family, a uniform consistency class is%
\[
\mathcal{Q}_{\mathrm{II},2}\left(  \boldsymbol{\theta},t,\pi_{1,m}%
;\gamma\right)  =\left\{
\begin{array}
[c]{c}%
\left\Vert \boldsymbol{\theta}\right\Vert _{\infty}\leq\rho,\pi_{1,m}\geq
m^{\left(  \gamma^{\prime}-1\right)  /2},t=\sqrt{\left\Vert \boldsymbol{\eta
}\right\Vert _{\infty}^{-1/2}\gamma\ln m},\\
\lim_{m\rightarrow\infty}t\min_{i\in I_{1,m}}\left\vert \eta_{0}-\eta
_{i}\right\vert =\infty
\end{array}
\right\}
\]
for any fixed $\gamma\in\left(  0,1\right)  $ and $\gamma^{\prime}>\gamma$.
The speed of convergence is poly-log.
\end{enumerate}
\end{theorem}

In \autoref{II-consistency}, the speed of convergence and uniform consistency
class depend on $\left\Vert \boldsymbol{\eta}\right\Vert _{\infty}$ and
$\min_{i\in I_{1,m}}\left\vert \eta_{0}-\eta_{i}\right\vert $, whereas they
only depend on $\min_{i\in I_{1,m}}\left\vert \mu_{i}-\mu_{0}\right\vert $ for
location-shift families. However, when the GF $H$ of the basis $\beta$ has
finite radius of convergence, their dependence on $\left\Vert \boldsymbol{\eta
}\right\Vert _{\infty}$ can be removed, as justified by:

\begin{corollary}
\label{CoroII}Let $\mathcal{F}$ be the NEF generated by $\beta$ in
(\ref{eq12a}), $\left\{  z_{i}\right\}  _{i=1}^{m}$ independent with CDFs
$\left\{  G_{\theta_{i}}\right\}  _{i=1}^{m}$ belonging to $\mathcal{F}$, and
$\rho$ a finite, positive constant. If (\ref{eqUPde}) holds and $H$ has a
finite radius of convergence $R_{H}$, then for positive and sufficiently large
$t$%
\[
\mathbb{V}\left\{  \hat{\varphi}_{m}\left(  t,\mathbf{z}\right)  -\varphi
_{m}\left(  t,\boldsymbol{\theta}\right)  \right\}  \leq\frac{C}{m}\frac
{\exp\left(  2tR_{H}^{1/2}\right)  }{\sqrt{t}\phi_{m}\left(
L,\boldsymbol{\theta}\right)  },
\]
and a uniform consistency class is%
\[
\mathcal{\tilde{Q}}_{\mathrm{II},1}\left(  \boldsymbol{\theta},t,\pi
_{1,m};\gamma\right)  =\left\{
\begin{array}
[c]{c}%
\left\Vert \boldsymbol{\theta}\right\Vert _{\infty}\leq\rho,t=2^{-1}%
R_{H}^{-1/2}\gamma\ln m,\\
\pi_{1,m}\geq m^{\left(  \gamma-1\right)  /2},\min_{i\in I_{1,m}}\left\vert
\eta_{0}-\eta_{i}\right\vert \geq\frac{\ln\ln m}{\sqrt{2\ln m}}%
\end{array}
\right\}
\]
for any $\gamma\in\left(  0,1\right]  $.
\end{corollary}

Since $\sup\Theta\leq\ln R_{H}$ and $\left\Vert \boldsymbol{\eta}\right\Vert
_{\infty}\leq R_{H}$, the proof \autoref{CoroII} follows easily from the first
claims of \autoref{II-concentration} and \autoref{II-consistency} and is
omitted. The uniform consistency class $\mathcal{\tilde{Q}}_{\mathrm{II}%
,3}\left(  \boldsymbol{\theta},t,\pi_{1,m},\gamma\right)  $ is fully
data-adaptive and only requires information on $\min_{i\in I_{1,m}}\left\vert
\eta_{0}-\eta_{i}\right\vert $, as do Jin's estimator of \cite{Jin:2008} for
Gaussian family and Construction I for location-shift families on $\min_{i\in
I_{1,m}}\left\vert \mu_{i}-\mu_{0}\right\vert $.

\section{Construction III: consistency and speed of convergence}

\label{secConsistent3}

We will focus on Gamma family and show that the corresponding estimators are
uniformly consistent. Recall $f_{\theta}$ defined by (\ref{eq1h}) for Gamma
family. Then $f_{\theta}\left(  x\right)  =O\left(  x^{\sigma-1}\right)  $ as
$x\rightarrow0+$, which tends to $0$ when $\sigma>1$. The next result provides
upper bounds on the variance and oscillations of $\hat{\varphi}_{m}\left(
t,\mathbf{z}\right)  -\varphi_{m}\left(  t,\boldsymbol{\theta}\right)  $ when
$t$ is positive and large.

\begin{theorem}
\label{ConcentrationIII}Consider Gamma family such that $\left\{
z_{j}\right\}  _{j=1}^{m}$ are independent with parameters $\left\{  \left(
\theta_{i},\sigma\right)  \right\}  _{i=1}^{m}$ and a fixed $\sigma>0$. Assume
$t$ is positive and sufficiently large and set $u_{3,m}=\min_{1\leq i\leq
m}\left\{  1-\theta_{i}\right\}  $. Then%
\[
\mathbb{V}\left\{  \hat{\varphi}_{m}\left(  t,\mathbf{z}\right)  -\varphi
_{m}\left(  t,\boldsymbol{\theta}\right)  \right\}  \leq V_{\mathrm{III}%
}=\frac{C}{m^2}\exp\left(  \frac{4t}{u_{3,m}}\right)
\sum_{i=1}^{m}\left(  \frac{t}{1-\theta_{i}}\right)  ^{3/4-\sigma}%
\]
and%
\begin{equation}
\Pr\left\{  \left\vert \hat{\varphi}_{m}\left(  t,\mathbf{z}\right)
-\varphi_{m}\left(  t,\boldsymbol{\theta}\right)  \right\vert \geq
\lambda\right\}  \leq\lambda^{-2}V_{\mathrm{III}}.
\label{eq18}%
\end{equation}

\end{theorem}

Recall
$\xi\left(  \theta\right)  =\left(  1-\theta\right)  ^{-1}$, such that
$u_{3,m}=\min_{1\leq i\leq m}\xi^{-1}\left(  \theta_{i}\right)  $. Using
\autoref{ConcentrationIII}, we show the uniform consistency and speed of
convergence of the estimator for Gamma family.

\begin{theorem}
\label{III-uniformConsistent}Consider Gamma family such that $\left\{
z_{j}\right\}  _{j=1}^{m}$ are independent with parameters $\left\{  \left(
\theta_{i},\sigma\right)  \right\}  _{i=1}^{m}$ and a fixed $\sigma>0$. Let
$\rho>0$ be a finite constant. If $\sigma>3/4$, then a uniform consistency
class%
\[
\mathcal{Q}_{\mathrm{III}}\left(  \boldsymbol{\theta},t,\pi_{1,m}%
;\gamma\right)  =\left\{
\begin{array}
[c]{c}%
\left\Vert \boldsymbol{\theta}\right\Vert_{\infty} \leq\rho,t=4^{-1}\gamma u_{3,m}\ln
m,\lim_{m\rightarrow\infty}u_{3,m}\ln m=\infty,\\
\pi_{1,m}\geq m^{\left(  \gamma-1\right)  /2},\lim_{m\rightarrow\infty}%
t\min_{i\in I_{1,m}}\left\vert \xi\left(  \theta_{0}\right)  -\xi\left(
\theta_{i}\right)  \right\vert =\infty
\end{array}
\right\}
\]
for any fixed $\gamma\in\left(  0,1\right]  $. On the other hand, if
$\sigma\leq3/4$, then%
\[
\mathcal{Q}_{\mathrm{III}}\left(  \boldsymbol{\theta},t,\pi_{1,m}%
;\gamma\right)  =\left\{
\begin{array}
[c]{c}%
\left\Vert \boldsymbol{\theta}\right\Vert_{\infty} \leq\rho,t=4^{-1}\gamma u_{3,m}\ln
m,\pi_{1,m}\geq m^{\left(  \gamma^{\prime}-1\right)  /2},\\
\lim_{m\rightarrow\infty}t\min_{i\in I_{1,m}}\left\vert \xi\left(  \theta
_{0}\right)  -\xi\left(  \theta_{i}\right)  \right\vert =\infty
\end{array}
\right\}
\]
for any fixed $\gamma\in\left(  0,1\right)  $ and $\gamma^{\prime}>\gamma$.
However, in either case, the speed of convergence is poly-log.
\end{theorem}

In \autoref{III-uniformConsistent}, the speed of convergence and uniform
consistency class depend on $u_{3,m}=\min_{1\leq i\leq m}\left\{  1-\theta
_{i}\right\}  $ and $\xi_{3,m}=\min_{i\in I_{1,m}}\left\vert \xi\left(
\theta_{0}\right)  -\xi\left(  \theta_{i}\right)  \right\vert $. Since
$\theta<1$ for Gamma family, $u_{3,m}$ measures how close a $G_{\theta_{i}}$
is to the singularity where a Gamma density is undefined, and it is sensible
to often assume $\liminf_{m\rightarrow\infty}u_{3,m}>0$. On the other hand,
$\sigma\xi\left(  \theta\right)  =\mu\left(  \theta\right)  $ for all
$\theta\in\Theta$ for Gamma family. So, $\xi_{3,m}$ measures the minimal
difference between the means of $G_{\theta_{i}}$ for $i\in I_{1,m}$ and
$G_{\theta_{0}}$, and $\xi_{3,m}$ cannot be too small relative to $t$ as
$t\rightarrow\infty$ in order for the estimator to achieve consistency.

\section{Simulation studies}
\label{secSim}

We present a simulation study on $\hat{\varphi}_{m}\left(t_{m},\mathbf{z}\right)  $, with comparison to the ``MR'' estimator of \cite{Meinshausen:2006} and the ``hybrid estimator" induced by ``Jin's estimator" of \cite{Jin:2008}.
Since the MR estimator is only applicable to p-values that have continuous distributions, its performance, when applied to discrete p-values such as those induced by Poisson and Negative Binomial distributions (to be considered hereunder), provides information on its robustness. Specifically, when $X_0$ is a realization of a Poisson or Negative Binomial random variable with CDF $F_{0}^{\ast}$, its p-value is $F_{0}^{\ast}\left(X_0\right)$.

\subsection{Simulation design}
\label{simDesign}

For $a<b$, let $\mathsf{Unif}\left(  a,b\right)  $ be the uniform random
variable or the uniform distribution on the closed interval $[a,b]$. We
consider $7$ values for $m$ as $10^3$, $5\times10^{3}$, $10^{4}$,
$5\times10^{4}$, $10^{5}$, $5\times10^{5}$ or $10^{6}$, and $4$ sparsity
levels for $\pi_{1,m}$, i.e., the dense regime $\pi_{1,m}=0.2$, moderately sparse
regime $\pi_{1,m}=m^{-0.2}$, critically sparse regime $\pi_{1,m}=m^{-0.5}$ and
very sparse regime $\pi_{1,m}=m^{-0.7}$, where we recall $\pi_{1,m}=1-m_{0}m^{-1}$ and
$m_{0}=\left\vert \left\{  1\leq i\leq m:\mu_{i}=\mu_{0}\right\}  \right\vert
$. Further, we consider $5$ distribution families $\mathcal{F}$, i.e.,
Laplace, Cauchy, Poisson, Negative Binomial and central Chi-square families, set as follows:

\begin{itemize}
\item For Laplace and Cauchy families, $\sigma=1$ and $\mu_{0}=0$ is set, and
the nonzero $\mu_{i}$'s are generated independently such that their absolute
values $\left\vert\mu_{i}\right\vert$ are from $\mathsf{Unif}\left(0.75,5\right)  $
but each $\mu_{i}$ has probability $0.5$ to be
negative or positive.

\item For Negative Binomial family, $n=5$ and $\theta_{0}=-4.5$, for
Poisson family $\theta_{0}=0.08$, and for Gamma family $\theta
_{0}=0.5$ and $\sigma=$ $6$. Note that for Gamma family, under the null
hypothesis the corresponding distribution is a central Chi-square distribution with
$2^{-1}\sigma$ degrees of freedom. For all three families, $\theta_{i}%
=\theta_{0}\rho_{i}$ for $i=m_{0}+1,\ldots,m_{0}+\lfloor2^{-1}m\pi_{m}\rfloor$
and $\theta_{i}=\theta_{0}\rho_{i}^{-1}$ for $i=m_{0}+\lfloor2^{-1}m\pi
_{m}\rfloor+1,\ldots,m$, where $\left\{  \rho_{i}\right\}  _{i=m_{0}+1}^{m}$
are independently generated from $\mathsf{Unif}\left(  10,13\right)  $ for
Poisson family, from $\mathsf{Unif}\left(  8,15\right)  $ for Negative Binomial family, and from
$\mathsf{Unif}\left(  1.2,1.5\right)  $ for Gamma family.
\end{itemize}
Each triple $\left(  m,\pi_{m},\mathcal{F}\right)  $ gives an experiment, and
there are a total of $140$ experiments. Each experiment is repeated independently
$250$ times so that summary statistics can be obtained. For an estimator $\hat{\pi}_{1,m}$ of $\pi_{1,m}$ for $m\geq1$ for each experiment,
the mean and standard deviation of $\tilde{\delta}_{m}=\hat{\pi}_{1,m}\pi_{1,m}^{-1}-1$ is estimated from the $250$ realizations.

Details on the implementations of the estimators to be compared are given below.
For the estimator $\hat{\varphi}_{m}\left(  t_{m},\mathbf{z}\right)  $, the
sequence $\left\{  t_{m}\right\}  _{m\geq1}$ is set by \autoref{CorConI} as
$t_{m}=\ln m$ for Laplace family and as $t_{m}=\sqrt{\ln m}$ for Cauchy
family, by \autoref{II-consistency} as $t_{m}=\left\Vert \boldsymbol{\eta
}\right\Vert _{\infty}^{-1/4}\sqrt{\ln m}$ for Poisson family, by
\autoref{CoroII} as $t_{m}=2^{-1}\ln m$ for Negative Binomial family, and by
\autoref{III-uniformConsistent} to be $t_{m}=4^{-1}\ln m\min_{1\leq i\leq
m}\left\{  1-\theta_{i}\right\}  $ for central Chi-square family, where the
tuning parameter $\gamma$ for $\hat{\varphi}_{m}\left(  t_{m},\mathbf{z}\right)  $ has been set to be the supremum of its feasible range to
enable the estimator to achieve a fast speed of convergence. The averaging
function $\omega$ is chosen to be the triangular density on $\left[
-1,1\right]  $. Further, from \autoref{Sec:EgMeth1},
\autoref{ThmDiscreteInfiniteSupport} and \autoref{SecEG3}, we see the
following: for Cauchy family with $\mu_{0}=0$ and $\sigma=1$,
\begin{equation}
K\left(  t,x;0\right)  =\int_{\left[  -1,1\right]  }\exp\left(  \left\vert
t\right\vert \right)  \omega\left(  s\right)  \cos\left(  tsx\right)
ds;\label{E1}%
\end{equation}
for Laplace family with $\mu_{0}=0$ and $\sigma=1$,%
\begin{equation}
K\left(  t,x;0\right)  =\int_{\left[  -1,1\right]  }\left(  1+t^{2}\right)
\omega\left(  s\right)  \cos\left(  tsx\right)  ds;\label{E2}%
\end{equation}
for Poisson family and Negative Binomial family%
\begin{equation}
K\left(  t,x;\theta_{0}\right)  =H\left(  e^{\theta_{0}}\right)  \int_{\left[
-1,1\right]  }\frac{\left(  ts\right)  ^{x}\cos\left(  \frac{\pi x}%
{2}-tse^{\theta_{0}}\right)  }{H^{\left(  x\right)  }\left(  0\right)  }%
\omega\left(  s\right)  ds;\label{E3}%
\end{equation}
for central Chi-square family%
\begin{equation}
K\left(  t,x;\theta_{0}\right)  =\int_{\left[  -1,1\right]  }\omega\left(
s\right)  \sum_{n=0}^{\infty}\frac{\left(  -tsx\right)  ^{n}\Gamma\left(
\sigma\right)  \cos\left(  \frac{\pi}{2}n+\frac{ts}{1-\theta_{0}}\right)
}{n!\Gamma\left(  \sigma+n\right)  }ds.\label{E4}%
\end{equation}
For the integrals in (\ref{E1}), (\ref{E2}), (\ref{E3}) and (\ref{E4}), each
integral is approximated by a Riemann sum for which the interval $\left[
-1,1\right]  $ is partitioned into $400$ equal subintervals, for which each integrand is evaluated at the end points of these subintervals. Further, for the
integral in (\ref{E4}), the power series in the integrand is replaced by the partial sum of its
first $21$ terms, i.e., it is truncated at $n=20$. The MR estimator (defined for continuous p-values) is
implemented as follows: let the ascendingly ordered p-values be $p_{\left(
1\right)  }<p_{\left(  2\right)  }<\cdots<p_{\left(  m\right)  }$ for $m>4$, set $b_{m}^{\ast}=m^{-1/2}\sqrt{2\ln\ln m}$,
define
\[
q_{i}^{\ast}=\left(  1-p_{\left(  i\right)  }\right)  ^{-1}\left\{
im^{-1}-p_{\left(  i\right)  }-b_{m}^{\ast}\sqrt{p_{\left(  i\right)  }\left(
1-p_{\left(  i\right)  }\right)  }\right\};
\]
then $\hat{\pi}_{1,m}^{mr}=\min\left\{  1,\max\left\{  0,\max_{2\leq i\leq
m-2}q_{i}^{\ast}\right\}  \right\}  $ is the MR estimator. Note that the MR estimator implicitly assumes that the probability for any tie between the p-values is zero and that when it is applied
to discrete p-values, it is ok to allow for such ties.
The hybrid estimator is implemented as
follows: first, each $z_{i}$ is transformed into $\tilde{z}_{i}=\Phi
^{-1}\left(  F_{i0}\left(  z_{i}\right)  \right)  $, where $F_{i0}$ is the CDF
of $z_{i}$ under the null hypothesis and $\Phi^{-1}$ the inverse of the
CDF\ of the standard Normal random variable; secondly, Jin's estimator of
\cite{Jin:2008} for Gaussian family is applied to $\left\{  \tilde{z}%
_{i}\right\}  _{i=1}^{m}$, for which the integral in (\ref{eq2}) is
approximated by a Riemann sum based on partitioning $\left[  -1,1\right]  $
into subintervals of equal length $0.01$ (the default setting in
\cite{Jin:2008}), $\omega$ is set as the triangular density on $\left[
-1,1\right]  $, and $\gamma=0.5$ is set in $t_{m}=\sqrt{2\gamma\ln m}$ to
provide the fastest possible speed of convergence.

\subsection{Simulation results}

For an estimator $\hat{\pi}_{1,m}$ of $\pi_{1,m}$, we will measure its stability by the standard deviation $\sigma_m^{\ast}$ of $\tilde{\delta}_{m}$ and its accuracy by the mean $\mu_{m}^{\ast}$ of $\tilde{\delta}_{m}$. Among two estimators for a fixed $m$, the one that has both smaller $\sigma_m^{\ast}$ and $\left\vert \mu_{m}^{\ast} \right\vert$ will be considered better. The supplementary material contains boxplots that summarize the performances of the three estimators under investigation.

The following five observations have been made from the comparison between the proposed estimator and the MR estimator: (1) For Laplace and Cauchy families, the new estimator is very accurate and much better than the MR estimator. In the dense and moderately sparse regimes, there is very strong evidence on the consistency of the new estimator since $\tilde{\delta}_{m}$ displays a strong trend to converge to $0$. These may hold true when the new estimator is applied to other location-shift families. (2) For Poisson and Negative Binomial families, the new estimator is accurate and much better than the MR estimator. In the dense and moderately sparse regimes, there is strong evidence on the convergence of $\tilde{\delta}_{m}$ (even though not necessarily to $0$) when the new estimator is applied to Negative Binomial family, whereas there is no strong evidence of convergence of $\tilde{\delta}_{m}$ when it is applied to Poisson family. In other words, we have not observed strong evidence on the consistency of the new estimator. This may be a consequence of non-adaptively choosing $t_{m}$ for the estimator, and is worth further investigation. In contrast, the MR estimator is almost always zero, i.e., it is rarely able to detect the existence of false null hypotheses. (3) For Gamma family, the new estimator is more accurate than the MR estimator, and it often severely underestimates $\pi_{1,m}$. When the new estimator is applied to the dense and moderately sparse regimes, there is strong evidence on the convergence of $\tilde{\delta}_{m}$ but there is no strong evidence on the consistency of the new estimator. This may be due to truncating the power series in Construction III when implementing the estimator for Gamma family and non-adaptively choosing $t_{m}$ for the estimator, and requires further investigation. In contrast, the MR estimator is almost always zero, often failing to detect the existence of false null hypotheses. Such an interesting behavior for the MR estimator has not been reported before. (4) The MR estimator, if not being $0$ for all almost all repetitions of an experiment, is less stable than the new estimator and can be much so when $m$ is large. In the critically sparse and very sparse regimes, the new estimator does not seem to be consistent and its $\tilde{\delta}_{m}$ does not display a trend of convergence as $m$ increases. Similarly, in these regimes, the MR estimator does not seem to be consistent unless it is identically zero. However, this does not contradict the theory for the MR estimator since its consistency requires that p-values under the alternative hypothesis be identically distributed and they are not so in the simulation study here. (5) In the dense and moderately sparse regimes, the new estimator usually underestimates $\pi_{1,m}$, i.e., the estimated proportion of true null hypotheses, $\hat{\pi}_{0,m}$, induced by the estimator is usually conservative. This is appealing in that a one-step adaptive FDR procedure that employs $\hat{\pi}_{0,m}$ is usually conservative.

On the other hand, the following three observations have been made from the comparison between the proposed estimator and the hybrid estimator: (1) when $z_{i}$'s have Laplace or Cauchy distributions, the proposed estimator is (much) more accurate but a bit less stable than the hybrid estimator in the dense and moderately sparse regimes, whereas they have competitive performances in the critically and very sparse regimes;
(2) when $z_{i}$'s have central Chi-square distributions, the hybrid estimator is more accurate but a bit less stable than the
proposed estimator in the dense and moderately sparse regimes, whereas they have competitive performances in the critically and very sparse regimes. However, for this scenario, the proposed estimator is implemented by truncating the power series in the integral in (\ref{E4}) at the $20$th
term, and its performance can be improved by better approximating the power series;
(3) when $z_{i}$'s have Poisson or Negative Binomial distributions, the proposed estimator is much more accurate than and as stable as the hybrid estimator across all sparsity regimes. This is reasonable since when $z_i$'s have discrete CDFs, the assumptions on the continuity and Normality of the transformed random variables $\tilde{z}_i$'s are violated, and Jin's method is not applicable. We remark that when $z_{i}$'s have central Chi-square distributions such that $\theta_{0}=0.05$ and $\sigma=9$ are set in \autoref{simDesign}, the proposed estimator is more accurate and stable than the hybrid estimator in the dense and moderately sparse regimes, whereas they have competitive performances in the critically and very sparse regimes.
In summary, the hybrid estimator is not able to serve as the universal estimator, and tailored ones such as the proposed are needed for specific scenarios.

\section{Discussion}

\label{SecConcAndDisc}

We have demonstrated that solutions of Lebesgue-Stieltjes integral equations
can serve as a universal construction for proportions estimators, provided
proportion estimators for random variables with three types of distributions,
and justified under independence the uniform consistency and speeds of
convergence of the estimators. For a proposed estimator to achieve uniform consistency, the tuning parameter
that determines its intrinsic speed of convergence needs to be determined adaptively based on data.
On the other hand, in applications we usually have information from domain scientists on a lower bound on
the proportion of false null hypotheses and on the minimal effect size. Further, for each estimator from each construction, an upper bound on its variance with respect to the oracle has been provided, and the minimal effect size to ensure its uniform consistency can converge to zero quite fast. So, inspired by the work of \cite{Jin:2008}, we can adaptively determine the tuning parameter so that the intrinsic speed forces the variance upper bound to converge to zero at certain rate and to be of smaller order than the lower bound on the proportion, thus achieving uniform consistency adaptively.
Specifically, this can be done for estimators from Construction I, for estimators from Construction II where generating functions have finite radii of convergence, and for other estimators of Construction II and those from Construction III after estimating the supremum norm of the parameter vector or the infimum of a known transform of the vector
itself.
In an accompanying article, we will deal with this estimation problem, discuss how
truncating the power series in Construction III affects the accuracy of the
induced estimators, and report via extensive simulation studies the adaptive,
non-asymptotic performances of the proposed estimators.

Our work induces three topics that are worthy of future investigations.
Firstly, we have only considered estimating the proportion of parameters that
are unequal to a fixed value, i.e., the proportion induced by the functional
that maps a parameter to a fixed value. It would be interesting to construct
uniformly consistent estimators of proportions induced by other functionals.
Further, we have only considered independent random
variables. Extending the consistency results provided here to dependent case
will greatly enlarge the scope of applications of the estimators, as did by
\cite{Chen:2018} and \cite{Jin:2007} to Jin's estimator of \cite{Jin:2008} for
Gaussian family and Gaussian mixtures. Moreover, Construction I, II and III
are applicable to random variables whose distributions are from different
sub-families of the same type of distributions, and results on the uniform
consistency of the corresponding proportion estimators can be extended to this
case. Finally, following the principles in Section 3 of \cite{Jin:2008},
Construction I, II and III can possibly be applied to consistently estimate
the mixing proportions for two-component mixture models at least one of whose
components follows a distribution discussed in this work.

Secondly, we have only been able to construct proportion estimators for three
types of distributions, and provide Gamma family as an example for
Construction III. It is worthwhile to explore other settings for which
solutions of Lebesgue-Stiejtjes integro-differential equations exit, can be
analytically expressed, and serve as consistent proportion estimators.
Further, we have not studied optimal properties of the proposed estimators, and with regard to this the techniques of \cite{Cai:2010} and \cite{Carpentier:2017} may
be useful.

Thirdly, we have introduced the concept of \textquotedblleft the family of
distributions with Riemann-Lebesgue type characteristic functions (RL type
CFs)\textquotedblright\ (see \autoref{Def}) for which%
\begin{equation}
\left\{  t\in\mathbb{R}:\hat{F}_{\mu_{0}}\left(  t\right)  =0\right\}
=\varnothing, \label{eq17b}%
\end{equation}
and shown that it contains several location-shift families. The requirement
(\ref{eq17b}) precludes the characteristic function $\hat{F}_{\mu_{0}}$ to
have any real zeros. We are aware that Poisson family is infinitely divisible
but not a location-shift family and does not have RL type CFs. However, it is
unclear to us the relationships (with respect to set inclusion) between
infinitely divisible distributions, location-shift families and distributions
with RL type CFs. So, a better understanding of such relationships will
contribute both to the theory of probability distributions and finding
examples different than those given here that Constructions I, II and III
apply to.

\section*{Acknowledgements}

Part of the research was funded by the New Faculty Seed Grant provided by
Washington State University. I would like to thank the Editor and Associate Editor for handling my submission and the reviewers
for their helpful comments.
I am very grateful to G\'{e}rard Letac for his
guidance on my research involving natural exponential families and constant
encouragements, Jiashun
Jin for providing the technical report \cite{Jin:2006b}, warm encouragements
and comments on the presentation of an earlier version of the manuscript, Mark
D. Ward for comments on exponential generating functions, Kevin Vixie,
Hong-Ming Yin and Charles N. Moore for discussions on solutions of
integro-differential equations, Sheng-Chi Liu for discussions on solutions of
algebraic equations, and Ovidiu Costin and Sergey Lapin for help with access
to two papers.

\bibliographystyle{myjmva}


\newpage

\appendix

\begin{center}
{\Large{Supplementary material for ``Uniformly consistently estimating the proportion of false null hypotheses via Lebesgue-Stieltjes integral equations"}}
\bigskip
\end{center}

We will discuss in \ref{SecIII_notworking} Construction III for Resell and
Hyperbolic Cosine families and in \ref{secUniformFrequence} uniform
consistency in frequency domain for Constructions II and III and its relation
to concentration inequalities for non-Lipschitz functions of independent
random variables. Proofs related to Construction I, II and III are provided
respectively in \ref{ProofsI}, \ref{ProofsII} and \ref{ProofsIII}. Additional simulation results are given in \ref{secSimAdd}. Here we will use $\log$ to denote
the natural logarithm in order to maintain consistency with the notation for the complex natural logarithm.

\section{On Construction III for Ressel and Hyperbolic Cosine families}

\label{SecIII_notworking}

For Ressel and Hyperbolic Cosine families, which are non-location-shift
NEF-CVFs, we suspect that Construction III cannot be implemented based on the
following initial results on their moment sequences.

\begin{example}
Ressel family with basis%
\begin{equation}
\frac{d\beta}{d\nu}\left(  x\right)  =f\left(  x\right)  =\frac{\sigma
x^{x+\sigma-1}e^{-x}}{\Gamma\left(  x+\sigma+1\right)  }1_{\left(
0,\infty\right)  }\left(  x\right)  \text{ \ for }\sigma>0 \label{eq3a}%
\end{equation}
and variance function $V\left(  \mu\right)  =\frac{\mu^{2}}{\sigma}\left(
1+\frac{\mu}{\sigma}\right)  $ for\ $\mu>0$. Note that $\int_{0}^{\infty}%
\beta\left(  dx\right)  =1$ by Proposition 5.5 of \cite{Letac:1990}. From
\cite{Lev:2016IJSP} and references therein, we know the following: the Laplace
transform $L_{\sigma}\left(  \theta\right)  $ of $\beta$ cannot be explicitly
expressed in $\theta$; $L_{1}\left(  \theta\right)  =\exp\left(  -\tilde{\eta
}\left(  -\theta\right)  \right)  $ where $\tilde{\eta}\left(  -\theta\right)
$ is the solution to the functional equation%
\[
\tilde{\eta}\left(  -\theta\right)  =\log\left(  1+\tilde{\eta}\left(
-\theta\right)  -\theta\right)  \text{ with }\theta\leq0;
\]
$L_{\sigma}\left(  \theta\right)  =\left(  L_{1}\left(  \theta\right)
\right)  ^{\sigma}$, $\theta\left(  \mu\right)  =\log\frac{1+\mu}{\mu}%
-\mu^{-1}$ and $L_{1}\left(  \theta\left(  \mu\right)  \right)  =\frac{\mu
}{1+\mu};$ $\lim_{x\rightarrow0+}\frac{d\beta}{d\nu}\left(  x\right)  =0$ when
$\sigma>1$.

Let us compute $\tilde{c}_{n}^{\ast}=\int x^{n}e^{x\theta}\beta\left(
dx\right)  $ for\ $n\in\mathbb{N}$. Recall Hankel's formula for the reciprocal
Gamma function, i.e.,%
\[
\frac{1}{\Gamma\left(  z\right)  }=\frac{\iota}{2\pi}\int_{\mathcal{C}}\left(
-t\right)  ^{-z}e^{-t}dt\text{ \ with \ }\Re\left(  z\right)  >0,
\]
where $\mathcal{C}$ is the
\href{https://en.wikipedia.org/wiki/Hankel_contour}{Hankel contour} that wraps
the non-negative real axis counterclockwise once and the logarithm function
$\log$ is such that $\log\left(  -t\right)  \in\mathbb{R}$ for $t<0$; see,
e.g., Section 12.22 of \cite{Whittaker:1940} for details on this. Then, for
$n\geq1$ we obtain%
\[
\frac{x^{x+\sigma}}{\Gamma\left(  x+\sigma+1\right)  }=\frac{\iota}{2\pi}%
\int_{\mathcal{C}}\left(  -t\right)  ^{-x-\sigma-1}e^{-xt}dt\text{ \ for
\ }x>0
\]
and
\begin{align}
\tilde{c}_{n}^{\ast}  &  =\int_{0}^{\infty}\frac{\sigma x^{x+\sigma+n-1}%
}{\Gamma\left(  x+\sigma+1\right)  }e^{-x\left(  1-\theta\right)
}dx\nonumber\\
&  =\int_{0}^{\infty}\sigma x^{n-1}\exp\left(  -x\left(  1-\theta
+t+\log\left(  -t\right)  \right)  \right)  dx\frac{\iota}{2\pi}%
\int_{\mathcal{C}}\left(  -t\right)  ^{-\sigma-1}dt\nonumber\\
&  =\frac{\iota\sigma\left(  n-1\right)  !}{2\pi}\int_{\mathcal{C}}%
\frac{\left(  -t\right)  ^{-\sigma-1}}{\left(  1-\theta+t+\log\left(
-t\right)  \right)  ^{n}}dt. \label{eq3b}%
\end{align}
Let $b\left(  \theta\right)  $ be the lower real branch of the solutions in
$t$ to the functional equation%
\[
1-\theta+t+\log\left(  -t\right)  =0.
\]
Then $b\left(  \theta\right)  $ is the Lambert W function $W_{-1}\left(
\tilde{z}\right)  $ with $\tilde{z}=-\exp\left(  \theta-1\right)  $ and domain
$\tilde{z}\in\lbrack-e^{-1},0)$ such that $W_{-1}\left(  \tilde{z}\right)  $
decreases from $W_{-1}\left(  -e^{-1}\right)  =-1$ to $W_{-1}\left(
0-\right)  =-\infty$; see, e.g., \cite{Corless1996} for details on this. When
$\theta<0$, $b\left(  \theta\right)  >-1$. Since $b\left(  \theta\right)  $ is
a pole of order $n$ for the integrand $\tilde{R}\left(  t\right)  $ in
(\ref{eq3b}), the residue theorem implies%
\[
\tilde{c}_{n}^{\ast}=-\sigma\lim_{t\rightarrow b\left(  \theta\right)  }%
\frac{d^{n-1}}{dt^{n-1}}\left(  \left(  t-b\left(  \theta\right)  \right)
^{n}\tilde{R}\left(  t\right)  \right)  .
\]
In particular,%
\[
\tilde{c}_{1}^{\ast}=-\sigma\frac{\left(  -b\left(  \theta\right)  \right)
^{-\sigma-1}}{1+b^{-1}\left(  \theta\right)  }=-\sigma\frac{\left(  -1\right)
^{-\sigma-1}W_{-1}^{-\sigma}\left(  -\exp\left(  \theta-1\right)  \right)
}{1+W_{-1}\left(  -\exp\left(  \theta-1\right)  \right)  },
\]
and $\tilde{c}_{n}^{\ast}$ is a complicated function in $W_{-1}\left(
-\exp\left(  \theta-1\right)  \right)  $ when $n$ is large. So, Ressel family
is unlikely to have a separable moment sequence.
\end{example}

\begin{example}
Hyperbolic Cosine family with basis%
\[
\frac{d\beta}{d\nu}\left(  x\right)  =\frac{2^{\sigma-2}}{\pi\Gamma\left(
\sigma\right)  }\left\vert \Gamma\left(  \frac{\sigma}{2}+\iota\frac{x}%
{2}\right)  \right\vert ^{2}\text{ \ for }\sigma>0\text{ and\ }x\in\mathbb{R}%
\]
and Fourier transform $L\left(  \iota t\right)  =\left(  \cosh t\right)
^{-\sigma}$ as shown on page 28 of \cite{Letac:1990}. So,
\[
L\left(  \theta\right)  =2^{\sigma}\left(  e^{-\iota\theta}+e^{-\iota\theta
}\right)  ^{-\sigma}=\left(  \cos\theta\right)  ^{-\sigma}\text{ \ for
\ }\left\vert \theta\right\vert <2^{-1}\pi,
\]
$\mu\left(  \theta\right)  =\sigma\tan\theta$ and $V\left(  \mu\right)
=\sigma\left(  1+\frac{\mu^{2}}{\sigma^{2}}\right)  $ for $\mu\in\mathbb{R}$.
From Theorems 2 and 3 of \cite{Morris:1982}, we see that $\tilde{c}_{n}^{\ast
}=\int x^{n}e^{\theta x}\beta\left(  dx\right)  $ for $n\geq3$ is a polynomial
of degree $n$ in $\mu$ with at least one non-zero term of order between $1$
and $n-1$. So, Hyperbolic Cosine family is unlikely to have a separable moment sequence.
\end{example}

\section{On uniform consistency in frequency domain for Construction II and
III}

\label{secUniformFrequence}

For Construction I applied to location-shift families with RL type CFs, we
have proved%
\[
\Pr\left(  \sup\nolimits_{\boldsymbol{\mu}\in\mathcal{B}_{m}\left(
\rho\right)  }\left\vert \pi_{1,m}^{-1}\sup\nolimits_{t\in\left[  0,\tau
_{m}\right]  }\hat{\varphi}_{m}\left(  t,\mathbf{z}\right)  -1\right\vert
\rightarrow0\right)  \rightarrow1,
\]
for which the estimator is also consistent uniformly in $t\in\left[
0,\tau_{m}\right]  $ for a positive, increasing sequence $\tau_{m}%
\rightarrow\infty$. This is referred to as \textquotedblleft uniform
consistency in frequency domain\textquotedblright. Even though theoretically
it provides much flexibility in choosing different sequences of values for $t$
when estimating $\pi_{1,m}$, it does not have much practical value since $t$
needs to be large for $\varphi_{m}\left(  t,\boldsymbol{\mu}\right)  $ to
converge to $\pi_{1,m}$ fast so that $\hat{\varphi}_{m}\left(  t,\mathbf{z}%
\right)  $ can accurately estimate $\pi_{1,m}$.

Such uniform consistency is a consequence of the uniform boundedness and
global Lipschitz property of the transform on $\left\{  z_{i}\right\}
_{i=1}^{m}$ that is used to construct $\left\{  K\left(  t,z_{i};\mu
_{0}\right)  \right\}  _{i=1}^{m}$. In contrast, for Constructions II and III,
the corresponding transform is not necessarily uniformly bounded or globally
Lipschitz (see the comparison below), and uniform consistency in frequency
domain is hard to achieve. In fact, it is very challenging to derive good
concentration inequalities for sums of transformed independent random
variables where the transform is neither bounded nor globally Lipschitz. For
progress along this line when the transform is a polynomial, we refer to
readers to \cite{Kim:2000}, \cite{Vu:2002} and \cite{Schudy:2012}.

Now we present the comparison. For location-shift families with RL type CFs,
recall%
\[
K\left(  t,x;\mu_{0}\right)  =\int_{\left[  -1,1\right]  }\frac{\omega\left(
s\right)  w\left(  ts,x\right)  }{r_{\mu_{0}}\left(  ts\right)  }ds,
\]
where $w\left(  y,x\right)  =\cos\left(  yx-h_{\mu_{0}}\left(  y\right)
\right)  $ and $\left\Vert \partial_{y}h_{\mu_{0}}\right\Vert _{\infty}%
=C_{\mu_{0}}<\infty$ is assumed. So,%
\begin{equation}
\left\Vert w\right\Vert _{\infty}<\infty\text{ \ and \ }\left\Vert
\partial_{y}w\left(  \cdot,x\right)  \right\Vert _{\infty}\leq\tilde{C}%
_{0}\left\vert x\right\vert +\tilde{C} \label{eq17a}%
\end{equation}
for finite, positive constants $\tilde{C}_{0}$ and $\tilde{C}$ that do not
depend on $x$. This, together with the location-shift property, implies
uniform consistency in frequency domain for $\hat{\varphi}_{m}\left(
t,\mathbf{z}\right)  $ for an admissible $\omega$. Let us examine
Constructions II and III. First, consider Construction II. For (\ref{II-c}),
i.e.,%
\[
K\left(  t,x;\theta_{0}\right)  =H\left(  \eta_{0}\right)  \int_{\left[
-1,1\right]  }w\left(  ts,x\right)  \omega\left(  s\right)  ds\text{ \ with
\ }\eta_{0}=e^{\theta_{0}},
\]
where%
\[
w\left(  y,x\right)  =\frac{y^{x}\cos\left(  2^{-1}\pi x-y\eta_{0}\right)
}{c_{x}x!}\text{ \ \ for }y\geq0\text{ \ and }x\in\mathbb{N}\text{,}%
\]
we have $\left\Vert w\right\Vert _{\infty}=\infty$, and $\left\Vert
\partial_{y}w\left(  \cdot,x\right)  \right\Vert _{\infty}\leq\tilde{C}%
_{0}\left\vert x\right\vert +\tilde{C}$ does not hold. Secondly, consider
Construction III. Recall \autoref{EgGamma} for Gamma family, i.e.,%
\[
K\left(  t,x;\mu_{0}\right)  =\Gamma\left(  \sigma\right)  \int_{\left[
-1,1\right]  }w\left(  ts,x\right)  \omega\left(  s\right)  ds,
\]
where%
\[
w\left(  y,x\right)  =\sum_{n=0}^{\infty}\frac{\left(  -yx\right)  ^{n}%
\cos\left(  2^{-1}\pi n+y\xi\left(  \theta_{0}\right)  \right)  }%
{n!\Gamma\left(  \sigma+n\right)  }\text{ \ for }y\geq0\text{\ and }x>0.
\]
Decompose $w\left(  y,x\right)  $ into the sum of four series%
\[
S_{l^{\prime}}\left(  x,y\right)  =\sum_{l=0}^{\infty}\frac{\left(
-yx\right)  ^{4l+l^{\prime}}\cos\left(  2^{-1}\pi\left(  4l+l^{\prime}\right)
+y\xi\left(  \theta_{0}\right)  \right)  }{\left(  4l+l^{\prime}\right)
!\Gamma\left(  \sigma+4l+l^{\prime}\right)  }\text{ \ for }l^{\prime}%
\in\left\{  0,1,2,3\right\}  .
\]
Then the summands in $S_{l^{\prime}}\left(  x,y\right)  $ for each $l^{\prime
}$ has a fixed sign uniformly in $x$ and $l$. Further, there exists a sequence
of $y\rightarrow\infty$ such that $\left\vert \cos\left(  \frac{\pi}{2}%
n+y\xi\left(  \theta_{0}\right)  \right)  \right\vert $ is positive uniformly
in $n$. Thus, there exists a sequence of $x$ such that $\left\Vert
w\right\Vert _{\infty}=\infty$.

\section{Proofs related to Construction I}

\label{ProofsI}

\subsection{Proof of \autoref{ThmSeparable}}

First of all, $K\left(  t,x;\mu_{0}\right)  $ defined by (\ref{eq2}) is the
real part of%
\[
K^{\dag}\left(  t,x;\mu_{0}\right)  =\int_{\left[  -1,1\right]  }\frac
{\omega\left(  s\right)  \exp\left(  \iota tsx\right)  }{\hat{F}_{\mu_{0}%
}\left(  ts\right)  }ds
\]
and $\psi\left(  t,\mu;\mu_{0}\right)  $ by (\ref{eq13}) the real part of%
\begin{equation}
\psi^{\dag}\left(  t,\mu;\mu_{0}\right)  =\int_{\left[  -1,1\right]  }%
\omega\left(  s\right)  \frac{\hat{F}_{\mu}\left(  ts\right)  }{\hat{F}%
_{\mu_{0}}\left(  ts\right)  }ds=\int_{\left[  -1,1\right]  }\omega\left(
s\right)  \frac{r_{\mu}\left(  ts\right)  }{r_{\mu_{0}}\left(  ts\right)
}e^{\iota\left(  h_{\mu}\left(  ts\right)  -h_{\mu_{0}}\left(  ts\right)
\right)  }ds. \label{eq8a}%
\end{equation}
With the boundedness of $\omega$, the uniform continuity of $r_{\mu}$ for each
$\mu\in U$, (\ref{eq6b}) and (\ref{eq6}), we can apply Fubini theorem to
obtain%
\[
\psi^{\dag}\left(  t,\mu;\mu_{0}\right)  =\int K^{\dag}\left(  t,x;\mu
_{0}\right)  dF_{\mu}\left(  x\right)  =\int_{\left[  -1,1\right]  }%
\frac{\omega\left(  s\right)  }{\hat{F}_{\mu_{0}}\left(  ts\right)  }%
\nu\left(  ds\right)  \int\exp\left(  \iota tsx\right)  dF_{\mu}\left(
x\right)  .
\]
This justifies (\ref{eq13}). If $\mu=\mu_{0}$, then $\frac{\hat{F}_{\mu}}%
{\hat{F}_{\mu_{0}}}\equiv1$, and (\ref{eq8a}) yields $\psi^{\dag}\left(
t,\mu_{0};\mu_{0}\right)  =1$ since $\omega$ is a density on $\left[
-1,1\right]  $, which justifies the first part of the second claim.

Finally, let $q_{1}\left(  y\right)  =\frac{r_{\mu}\left(  y\right)  }%
{r_{\mu_{0}}\left(  y\right)  }$ for $y\in\mathbb{R}$. Then (\ref{eq6b}) and
(\ref{eq6}) imply%
\[
\sup_{\left(  s,t\right)  \in\left[  -1,1\right]  \times\mathbb{R}}%
q_{1}\left(  st\right)  \leq C<\infty.
\]
Since (\ref{eq6d}) holds, Theorem 3 of \cite{Costin:2016} implies%
\begin{equation}
\lim_{t\rightarrow\infty}\int_{\left[  -1,1\right]  }\omega\left(  s\right)
\frac{\hat{F}_{\mu}\left(  ts\right)  }{\hat{F}_{\mu_{0}}\left(  ts\right)
}ds=0\text{ \ for \ }\mu\neq\mu_{0}, \label{eq8b}%
\end{equation}
which justifies the third claim.

\subsection{Proof of \autoref{MethodICor}}

It suffices to show (\ref{eq8b}). First of all, both $q_{1}\left(  ts\right)
=\frac{r_{\mu}\left(  ts\right)  }{r_{\mu_{0}}\left(  ts\right)  }$ and
$\omega\left(  s\right)  q_{1}\left(  ts\right)  $ belong to $L^{1}\left(
\left[  -1,1\right]  \right)  $ uniformly in $t\in\mathbb{R}$. For any
$\epsilon>0$, there exists a step function $q_{1,\epsilon}$ on $\mathbb{R}$
with compact support $A_{0}$ such that%
\[
q_{1,\epsilon}\left(  y\right)  =\sum_{j=1}^{n_{2,\epsilon}}a_{2j}1_{A_{2j}%
}\left(  y\right)  \text{ \ and \ }\int_{\mathbb{R}}\left\vert q_{1,\epsilon
}\left(  y\right)  -q_{1}\left(  y\right)  \right\vert dy<\epsilon,
\]
where $n_{2,\epsilon}\in\mathbb{N}$ is finite and the sets $\left\{
A_{2j}\right\}  _{j=1}^{n_{2,\epsilon}}$ are disjoint and $\bigcup
\nolimits_{j=1}^{n_{2,\epsilon}}A_{2j}\subseteq A_{0}$. Now consider $t$ with
$\left\vert t\right\vert \geq1$. Then, the boundedness of $\omega$ and
$\frac{r_{\mu}}{r_{\mu_{0}}}$ implies%
\begin{equation}
\int_{\left[  -1,1\right]  }\left\vert \omega\left(  s\right)  q_{1}\left(
ts\right)  -\omega\left(  s\right)  q_{1,\epsilon}\left(  ts\right)
\right\vert ds\leq C\left\vert t\right\vert ^{-1}\int_{\mathbb{R}}\left\vert
q_{1,\epsilon}\left(  y\right)  -q_{1}\left(  y\right)  \right\vert
dy\leq2C\epsilon. \label{eq10b}%
\end{equation}
Let $\tau\left(  ts\right)  =h_{\mu}\left(  ts\right)  -h_{\mu_{0}}\left(
ts\right)  $ and%
\[
a_{t}^{\ast}=\max_{1\leq j\leq n_{2,\epsilon}}\nu\left(  \left\{  s\in\left[
-1,1\right]  :ts\in A_{2j}\right\}  \right)  .
\]
Since the sets $\left\{  A_{2,j}\right\}  _{j=1}^{n_{2,\epsilon}}$ are
uniformly bounded, $\lim_{\left\vert t\right\vert \rightarrow\infty}%
a_{t}^{\ast}=0$ and
\begin{equation}
\int_{\left[  -1,1\right]  }\left\vert \omega\left(  s\right)  q_{1,\epsilon
}\left(  ts\right)  \exp\left(  \iota\tau\left(  ts\right)  \right)
\right\vert ds\leq Cn_{2,\epsilon}a_{t}^{\ast}\rightarrow0\text{ \ as
}\left\vert t\right\vert \rightarrow\infty. \label{eq10c}%
\end{equation}
Combining (\ref{eq10b}) and (\ref{eq10c}) gives (\ref{eq8b}), which justifies
the claim.

\subsection{Proof of \autoref{lmLoc}}

Since $\mathcal{F}$ is a location-shift family, if $z$ has CDF $F_{\mu}$ with
$\mu\in U$, then there exists some $\mu_{0}\in U$ such that $z=\mu^{\prime
}+z^{\prime}$, where $z^{\prime}$ has CDF $F_{\mu_{0}}$ and $\mu^{\prime}%
=\mu-\mu_{0}$. So,%
\[
\hat{F}_{\mu}\left(  t\right)  =\mathbb{E}\left[  \exp\left(  \iota tz\right)
\right]  =\mathbb{E}\left[  \exp\left(  \iota t\left(  \mu^{\prime}+z^{\prime
}\right)  \right)  \right]  =\hat{F}_{\mu_{0}}\left(  t\right)  \exp\left(
\iota t\mu^{\prime}\right)
\]
for all $t$. In particular, in the representation $\hat{F}_{\mu}=r_{\mu
}e^{\iota h_{\mu}}$, the modulus $r_{\mu}$ does not depend on $\mu$ and
$h_{\mu}\left(  t\right)  =t\mu^{\prime}$. If $\hat{F}_{\mu_{0}}\left(
t\right)  \neq0$ for all $t\in\mathbb{R}$, then (\ref{eq6d}) holds and
$\mathcal{\hat{F}}$ is of RL type.

\subsection{Proof of \autoref{CorLocationShift}}

When $\mathcal{F}$ is a location-shift family,%
\[
\int_{A}dF_{\mu}\left(  x\right)  =\int_{A-\left(  \mu-\mu_{0}\right)
}dF_{\mu_{0}}\left(  y\right)
\]
for each $A\subseteq\mathbb{R}$ measurable with respect to $F_{\mu_{0}}$.
Therefore,%
\[
\int K\left(  t,x;\mu_{0}\right)  dF_{\mu}\left(  x\right)  =\int K\left(
t,y+\left(  \mu-\mu_{0}\right)  ;\mu_{0}\right)  dF_{\mu_{0}}\left(  y\right)
\]
and the first identity in (\ref{eq7a}) holds. Further, $\frac{r_{\mu}}%
{r_{\mu_{0}}}\equiv1$ for all $\mu\in U$. So, (\ref{eq13}) reduces to the
second identity in (\ref{eq7a}). Finally, since the cosine function is even on
$\mathbb{R}$, it suffices to consider $t$ and $\mu$ such that $t\left(
\mu-\mu_{0}\right)  >0$ in the representation (\ref{eq7a}). When $\omega$ is
good, the proof of the third claim of Lemma 7.1 of \cite{Jin:2006b} remains
valid, which implies $1\geq\psi\left(  t,\mu;\mu_{0}\right)  \geq0$ for all
$\mu$ and $t$.

\subsection{Proof of \autoref{ThmConsistency}}

Recall
\[
K\left(  t,x;\mu_{0}\right)  =\int_{\left[  -1,1\right]  }\frac{\omega\left(
s\right)  \cos\left(  tsx-h_{\mu_{0}}\left(  ts\right)  \right)  }{r_{\mu_{0}%
}\left(  ts\right)  }ds.
\]
Set $w_{i}\left(  y\right)  =\cos\left(  yz_{i}-h_{\mu_{0}}\left(  y\right)
\right)  $ for each $i$ and $y\in\mathbb{R}$ and define%
\begin{equation}
S_{m}\left(  y\right)  =\frac{1}{m}\sum_{i=1}^{m}\left(  w_{i}\left(
y\right)  -\mathbb{E}\left[  w_{i}\left(  y\right)  \right]  \right)  .
\label{eq2e}%
\end{equation}
Then%
\begin{align*}
\hat{\varphi}_{m}\left(  t,\mathbf{z}\right)  -\varphi_{m}\left(
t,\boldsymbol{\mu}\right)   &  =\dfrac{1}{m}\sum_{i=1}^{m}\left(  K\left(
t,z_{i};\mu_{0}\right)  -\mathbb{E}\left[  K\left(  t,z_{i};\mu_{0}\right)
\right]  \right) \\
&  =\int_{\left[  -1,1\right]  }\frac{\omega\left(  s\right)  }{r_{\mu_{0}%
}\left(  ts\right)  }S_{m}\left(  ts\right)  ds.
\end{align*}
Since $\left\vert w_{i}\left(  ts\right)  \right\vert \leq1$ uniformly in
$\left(  t,s,z_{i},i\right)  $ and $\left\{  z_{i}\right\}  _{i=1}^{m}$ are
independent, (\ref{eq2g}) holds. Further, Hoeffding inequality of
\cite{Hoeffding:1963} implies%
\begin{equation}
\Pr\left(  \left\vert S_{m}\left(  ts\right)  \right\vert \geq\frac{\lambda
}{\sqrt{m}}\right)  \leq2\exp\left(  -2^{-1}\lambda^{2}\right)  \text{ for any
}\lambda>0 \label{eq2c}%
\end{equation}
uniformly in $\left(  t,s,m\right)  \in\mathbb{R}\times\left[  -1,1\right]
\times\mathbb{N}_{+}$. Recall $a\left(  t;\mu_{0}\right)  =\int_{\left[
-1,1\right]  }\frac{ds}{r_{\mu_{0}}\left(  ts\right)  }$ for $t\in\mathbb{R}$.
Therefore,%
\[
\left\vert \hat{\varphi}_{m}\left(  t,\mathbf{z}\right)  -\varphi_{m}\left(
t,\boldsymbol{\mu}\right)  \right\vert \leq\frac{\lambda\left\Vert
\omega\right\Vert _{\infty}}{\sqrt{m}}\int_{\left[  -1,1\right]  }\frac
{1}{r_{\mu_{0}}\left(  ts\right)  }ds=\frac{\lambda\left\Vert \omega
\right\Vert _{\infty}}{\sqrt{m}}a\left(  t;\mu_{0}\right)
\]
with probability $1-2\exp\left(  -2^{-1}\lambda^{2}\right)  $, i.e.,
(\ref{eq2f}) holds.

Consider the second claim. With probability at least $1-2\exp\left(
-2^{-1}\lambda_{m}^{2}\right)  $,%
\begin{align}
\left\vert \frac{\hat{\varphi}_{m}\left(  t_{m},\mathbf{z}\right)  }{\pi
_{1,m}}-1\right\vert  &  \leq\left\vert \frac{\hat{\varphi}_{m}\left(
t_{m},\mathbf{z}\right)  -\varphi_{m}\left(  t_{m},\boldsymbol{\mu}\right)
}{\pi_{1,m}}\right\vert +\left\vert \frac{\varphi_{m}\left(  t_{m}%
,\boldsymbol{\mu}\right)  }{\pi_{1,m}}-1\right\vert \nonumber\\
&  \leq\frac{1}{\pi_{1,m}}\frac{\lambda_{m}\left\Vert \omega\right\Vert
_{\infty}a\left(  t_{m};\mu_{0}\right)  }{\sqrt{m}}+\frac{1}{m\pi_{1,m}}%
\sum_{j\in I_{1,m}}\left\vert \psi\left(  t_{m},\mu_{j};\mu_{0}\right)
\right\vert \nonumber\\
&  \leq\frac{1}{\pi_{1,m}}\frac{\lambda_{m}\left\Vert \omega\right\Vert
_{\infty}a\left(  t_{m};\mu_{0}\right)  }{\sqrt{m}}+\sup_{\left(  t,\left\vert
\mu\right\vert \right)  \in\lbrack t_{m},\infty)\times\lbrack u_{m},\infty
)}\left\vert \psi\left(  t,\mu;\mu_{0}\right)  \right\vert . \label{eq2i}%
\end{align}
However, assumptions (\ref{eq4}), (\ref{eq2h}) and $\left\{  \mu_{i}:i\in
I_{1,m}\right\}  \subseteq\lbrack u_{m},\infty)$ imply that both $\exp\left(
-2^{-1}\lambda_{m}^{2}\right)  $ and the upper bound in (\ref{eq2i}) converge
to $0$. So, (\ref{eq4a}) holds.

\subsection{Proof of \autoref{CorConI}}

Recall \autoref{ThmConsistency} and its proof. First of all, $\lim
_{m\rightarrow\infty}\sup\left\{  \psi\left(  t_{m},\mu_{i};\mu_{0}\right)
:i\in I_{1,m}\right\}  =0$ when $t_{m}\rightarrow\infty$ and $\min_{i\in
I_{1,m}}\left\vert \mu_{i}-\mu_{0}\right\vert \geq\left(  t_{m}\right)
^{-1}\log\log m$. Therefore, it suffices to show%
\begin{equation}
\lim_{m\rightarrow\infty}\frac{\lambda_{m}a\left(  t_{m};\mu_{0}\right)  }%
{\pi_{1,m}\sqrt{m}}=0\ \text{and}\ \lim_{m\rightarrow\infty}\exp\left(
-2^{-1}\lambda_{m}^{2}\right)  =0, \label{eqAx2}%
\end{equation}
where%
\[
a\left(  t;\mu_{0}\right)  =\int_{\left[  -1,1\right]  }\frac{ds}{r_{\mu_{0}%
}\left(  ts\right)  }\text{ }\ \text{for }t\in\mathbb{R}.
\]

Consider $t>0$. For Gaussian family, $r_{\mu_{0}}^{-1}\left(  t\right)
=\exp\left(  2^{-1}t^{2}\sigma^{2}\right)  $. Setting $t_{m}=\sigma^{-1}%
\sqrt{2\gamma\log m}$ for any $\gamma\in\left(  0,0.5\right]  $ gives%
\begin{equation}
\int_{\left[  0,1\right]  }r_{\mu_{0}}^{-1}\left(  ts\right)  ds=\int_{\left[
0,1\right]  }\exp\left(  2^{-1}t_{m}^{2}s^{2}\sigma^{2}\right)  ds\leq
\frac{m^{\gamma}}{\gamma\log m}\left(  1+o\left(  1\right)  \right)  ,
\label{eqB5}%
\end{equation}
where the last inequality follows from the proof of Theorem 4 of
\cite{Jin:2008}. Set $\lambda_{m}=o\left(  t_{m}\right)  $ with $\lambda
_{m}\rightarrow\infty$. Then, (\ref{eqAx2}) holds for all $\pi_{1,m}\geq
Cm^{\gamma-0.5}$.

For Laplace family, $r_{\mu_{0}}^{-1}\left(  t\right)  =1+\sigma^{2}t^{2}$ and%
\begin{equation}
\int_{\left[  0,1\right]  }r_{\mu_{0}}^{-1}\left(  ts\right)  ds\leq
\int_{\left[  0,1\right]  }\left(  1+\sigma^{2}t^{2}\right)  ds=1+\frac
{\sigma^{2}t^{3}}{3}. \label{eqB1}%
\end{equation}
Set $t_{m}=\log m$ and $\lambda_{m}=O\left(  t_{m}\right)  $ with $\lambda
_{m}\rightarrow\infty$. Then, (\ref{eqAx2}) holds for all $\pi_{1,m}\geq
Cm^{-\gamma}$ with $0\leq\gamma<1/2$.

For Hyperbolic Secant family, $r_{\mu}^{-1}\left(  t\right)  =\sigma
\cosh\left(  t\sigma^{-1}\right)  \sim2^{-1}\sigma\exp\left(  \sigma
^{-1}t\right)  $ as $t\rightarrow\infty$ and%
\begin{equation}
\int_{\left[  0,1\right]  }r_{\mu}^{-1}\left(  ts\right)  ds\leq C\sigma
\int_{\left[  0,1\right]  }\exp\left(  \sigma^{-1}ts\right)  ds\leq
C\sigma^{2}\frac{\exp\left(  \sigma^{-1}t\right)  }{t}. \label{eqB2}%
\end{equation}
Set $t_{m}=\sigma\gamma\log m$ with $0<\gamma\leq1/2$ and $\lambda
_{m}=o\left(  t_{m}\right)  $ with $\lambda_{m}\rightarrow\infty$. Then,
(\ref{eqAx2}) holds for all $\pi_{1,m}\geq Cm^{\gamma-0.5}$.

For Logistic family, $r_{\mu}^{-1}\left(  t\right)  =\left(  \pi\sigma
t\right)  ^{-1}\sinh\left(  \pi\sigma t\right)  \sim\left(  2\pi\sigma
t\right)  ^{-1}e^{\pi\sigma t}$ as $t\rightarrow\infty$ and%
\begin{equation}
\int_{\left[  0,1\right]  }r_{\mu}^{-1}\left(  t\right)  ds\leq\left(
\sigma\pi t\right)  ^{-1}\int_{\left[  0,1\right]  }\exp\left(  \sigma\pi
ts\right)  ds=\frac{\exp\left(  \sigma\pi t\right)  }{\sigma^{2}\pi^{2}t^{2}}.
\label{eqB3}%
\end{equation}
Set $t_{m}=\left(  \sigma\pi\right)  ^{-1}\gamma\log m$ with $0<\gamma\leq1/2$
and $\lambda_{m}=o\left(  t_{m}\right)  $ with $\lambda_{m}\rightarrow\infty$.
Then, (\ref{eqAx2}) holds for all $\pi_{1,m}\geq Cm^{\gamma-0.5}$.

For Cauchy family, $r_{\mu}^{-1}\left(  t\right)  =\exp\left(  \sigma
\left\vert t\right\vert \right)  $ and%
\begin{equation}
\int_{\left[  0,1\right]  }r_{\mu_{0}}^{-1}\left(  ts\right)  ds\leq
\int_{\left[  0,1\right]  }\exp\left(  \sigma\left\vert t\right\vert s\right)
ds\leq\frac{\exp\left(  \sigma\left\vert t\right\vert \right)  }%
{\sigma\left\vert t\right\vert }. \label{eqB4}%
\end{equation}
Set $t_{m}=\sigma^{-1}\gamma\log m$ with $0<\gamma\leq1/2$ and $\lambda
_{m}=o\left(  t_{m}\right)  $ with $\lambda_{m}\rightarrow\infty$. Then,
(\ref{eqAx2}) holds for all $\pi_{1,m}\geq Cm^{\gamma-0.5}$.

In each case above, $C>0$ can be any constant for which $\pi_{1,m}\in\left(
0,1\right]  $ as $\gamma$ varies in its designated range.

\subsection{Proof of \autoref{ThmLocationShift}}

The strategy of proof adapts that for Lemma 7.2 of \cite{Jin:2006b} for
Gaussian family, which can be regarded as an application of the
\textquotedblleft chaining method\textquotedblright\ proposed by
\cite{Talagrand:1996}. Since $r_{\mu_{0}}$ has no real zeros and $\mathcal{F}$
is a location-shift family, $r_{\mu}$ has no real zeros for each $\mu\neq
\mu_{0}$ and $h_{\mu}\left(  t\right)  $ is well-defined and continuous in $t$
on $\mathbb{R}$ for each $\mu\in U$. Therefore, $\frac{d}{dy}h_{\mu_{0}%
}\left(  y\right)  $ can be defined.

Recall $w_{i}\left(  y\right)  =\cos\left(  yz_{i}-h_{\mu_{0}}\left(
y\right)  \right)  $ and $S_{m}\left(  y\right)  $ defined by (\ref{eq2e}).
Let $\hat{s}_{m}\left(  y\right)  =\frac{1}{m}\sum_{i=1}^{m}w_{i}\left(
y\right)  $ and $s_{m}\left(  y\right)  =\mathbb{E}\left[  \hat{s}_{m}\left(
y\right)  \right]  $. For the rest of the proof, we will first assume the
existence of the positive constants $\gamma^{\prime}$, $\gamma^{\prime\prime}%
$, $q$, $\vartheta$ and the non-negative constant $\vartheta^{\prime}$ and
then determine them at the end of the proof. Let $\gamma_{m}=\gamma^{\prime
}\log m$. The rest of the proof is divided into three parts.

\textbf{Part I}: to show the assertion \textquotedblleft if%
\begin{equation}
\lim_{m\rightarrow\infty}\frac{m^{\vartheta}\log\gamma_{m}}{R_{m}\left(
\rho\right)  \sqrt{m}\sqrt{2q\gamma_{m}}}=\infty\label{eq11b}%
\end{equation}
where $R_{m}\left(  \rho\right)  =2\mathbb{E}\left[  \left\vert X_{1}%
\right\vert \right]  +2\rho+2C_{\mu_{0}}$ and $X_{1}$ has CDF $F_{\mu_{0}}$,
then, for all large $m$,%
\begin{equation}
\sup_{\boldsymbol{\mu}\in\mathcal{B}_{m}\left(  \rho\right)  }\sup
_{y\in\left[  0,\gamma_{m}\right]  }\left\vert \hat{s}_{m}\left(  y\right)
-\mathbb{E}\left[  \hat{s}_{m}\left(  y\right)  \right]  \right\vert \leq
\frac{\sqrt{2q\gamma_{m}}}{\sqrt{m}} \label{eq11j}%
\end{equation}
holds with probability at least $1-p_{m}\left(  \vartheta,q,h_{\mu_{0}}%
,\gamma_{m}\right)  $, where%
\begin{equation}
p_{m}\left(  \vartheta,q,h_{\mu_{0}},\gamma_{m}\right)  =2m^{\vartheta}%
\gamma_{m}^{2}\exp\left(  -q\gamma_{m}\right)  +4A_{\mu_{0}}q\gamma
_{m}m^{-2\vartheta}\left(  \log\gamma_{m}\right)  ^{-2} \label{eq11e}%
\end{equation}
and $A_{\mu_{0}}$ is the variance of $\left\vert X_{1}\right\vert
$\textquotedblright.

Define the closed interval $G_{m}=\left[  0,\gamma_{m}\right]  $. Let
$\mathcal{P}=\left\{  y_{1},\ldots,y_{l_{\ast}}\right\}  $ for some $l_{\ast
}\in\mathbb{N}_{+}$ with $y_{j}<y_{j+1}$ be a partition of $G_{m}$ with norm
$\Delta=\max_{1\leq j\leq l_{\ast}-1}\left\vert y_{j+1}-y_{j}\right\vert $
such that $\Delta=m^{-\vartheta}$. For each $y\in G_{m}$, pick $y_{j}%
\in\mathcal{P}$ that is the closest to $y$. By Lagrange mean value theorem,%
\begin{align*}
\left\vert \hat{s}_{m}\left(  y\right)  -s_{m}\left(  y\right)  \right\vert
&  \leq\left\vert \hat{s}_{m}\left(  y_{i}\right)  -s_{m}\left(  y_{i}\right)
\right\vert +\left\vert \left(  \hat{s}_{m}\left(  y\right)  -\hat{s}%
_{m}\left(  y_{i}\right)  \right)  -\left(  s_{m}\left(  y\right)
-s_{m}\left(  y_{i}\right)  \right)  \right\vert \\
&  \leq\left\vert \hat{s}_{m}\left(  y_{i}\right)  -s_{m}\left(  y_{i}\right)
\right\vert +\Delta\sup_{y\in\mathbb{R}}\left\vert \partial_{y}\left(  \hat
{s}_{m}\left(  y\right)  -s_{m}\left(  y\right)  \right)  \right\vert ,
\end{align*}
where $\partial_{\cdot}$ denotes the derivative with respect to the subscript.
So,%
\begin{equation}
B_{0}=\Pr\left(  \sup_{\boldsymbol{\mu}\in\mathcal{B}_{m}\left(  \rho\right)
}\max_{y\in G_{m}}\left\vert S_{m}\left(  y\right)  \right\vert \geq
\frac{\sqrt{2q\gamma_{m}}}{\sqrt{m}}\right)  \leq B_{1}+B_{2}, \label{eq10h}%
\end{equation}
where%
\[
B_{1}=\Pr\left(  \sup_{\boldsymbol{\mu}\in\mathcal{B}_{m}\left(  \rho\right)
}\max_{1\leq i\leq l_{\ast}}\left\vert \hat{s}_{m}\left(  y_{i}\right)
-s_{m}\left(  y_{i}\right)  \right\vert \geq\frac{\sqrt{2q\gamma_{m}}-\left(
2q\gamma_{m}\right)  ^{-1/2}\log\gamma_{m}}{\sqrt{m}}\right)
\]
and%
\[
B_{2}=\Pr\left(  \sup_{\boldsymbol{\mu}\in\mathcal{B}_{m}\left(  \rho\right)
}\sup_{y\in\mathbb{R}}\left\vert \partial_{y}\left(  \hat{s}_{m}\left(
y\right)  -s_{m}\left(  y\right)  \right)  \right\vert \geq\frac{\Delta
^{-1}\left(  2q\gamma_{m}\right)  ^{-1/2}\log\gamma_{m}}{\sqrt{m}}\right)  .
\]
Applying to $B_{1}$ the union bound and Hoeffding inequality (\ref{eq2c})
gives%
\begin{equation}
B_{1}\leq2l_{\ast}\exp\left(  -q\gamma_{m}+\log\gamma_{m}\right)  \exp\left(
-\frac{\left(  \log\gamma_{m}\right)  ^{2}}{4q\gamma_{m}}\right)
\leq2m^{\vartheta}\gamma_{m}^{2}\exp\left(  -q\gamma_{m}\right)  .
\label{eq10g}%
\end{equation}

On the other hand, $\partial_{y}w_{i}\left(  y\right)  =-\left(
z_{i}-\partial_{y}h_{\mu_{0}}\left(  y\right)  \right)  \sin\left(
yz_{i}-h_{\mu_{0}}\left(  y\right)  \right)  $, and%
\[
\partial_{y}\mathbb{E}\left[  w_{i}\left(  y\right)  \right]  =\mathbb{E}%
\left[  \partial_{y}w_{i}\left(  y\right)  \right]  =-\mathbb{E}\left[
\left(  z_{i}-\partial_{y}h_{\mu_{0}}\left(  y\right)  \right)  \sin\left(
yz_{i}-h_{\mu_{0}}\left(  y\right)  \right)  \right]
\]
holds since $\int\left\vert x\right\vert ^{2}dF_{\mu}\left(  x\right)
<\infty$ for each $\mu\in U$ and $\sup_{y\in\mathbb{R}}\left\vert \frac{d}%
{dy}h_{\mu_{0}}\left(  y\right)  \right\vert =C_{\mu_{0}}<\infty$. So,%
\begin{equation}
\sup_{y\in\mathbb{R}}\left\vert \partial_{y}\left(  \hat{s}_{m}\left(
y\right)  -s_{m}\left(  y\right)  \right)  \right\vert \leq\frac{1}{m}%
\sum_{i=1}^{m}\left\vert z_{i}\right\vert +2C_{\mu_{0}}+\frac{1}{m}\sum
_{i=1}^{m}\mathbb{E}\left[  \left\vert z_{i}\right\vert \right]
.\label{eq10d}%
\end{equation}
Since $\mathcal{F}$ is a location-shift family, there are independent and
identically distributed (i.i.d.) $\left\{  X_{i}\right\}  _{i=1}^{m}$ with
common CDF $F_{\mu_{0}}$ such that $z_{i}=\left(  \mu_{i}-\mu_{0}\right)
+X_{i}$ for $1\leq i\leq m$. Therefore, the upper bound in (\ref{eq10d})
satisfies%
\begin{align*}
&  \frac{1}{m}\sum_{i=1}^{m}\left\vert z_{i}\right\vert +2C_{\mu_{0}}+\frac
{1}{m}\sum_{i=1}^{m}\mathbb{E}\left[  \left\vert z_{i}\right\vert \right]  \\
&  \leq\frac{1}{m}\sum_{i=1}^{m}\left\vert X_{i}\right\vert +\frac{2}{m}%
\sum_{i=1}^{m}\left\vert \mu_{i}-\mu_{0}\right\vert +2C_{\mu_{0}}+\frac{1}%
{m}\sum_{i=1}^{m}\mathbb{E}\left[  \left\vert X_{i}\right\vert \right]  \\
&  \leq\frac{1}{m}\sum_{i=1}^{m}\left(  \left\vert X_{i}\right\vert
-\mathbb{E}\left[  \left\vert X_{i}\right\vert \right]  \right)  +R_{m}\left(
\rho\right)  ,
\end{align*}
where we recall $R_{m}\left(  \rho\right)  =2\mathbb{E}\left[  \left\vert
X_{1}\right\vert \right]  +2\rho+2C_{\mu_{0}}$. Namely,%
\[
\sup_{y\in\mathbb{R}}\left\vert \partial_{y}\left(  \hat{s}_{m}\left(
y\right)  -s_{m}\left(  y\right)  \right)  \right\vert \leq\frac{1}{m}%
\sum_{i=1}^{m}\left(  \left\vert X_{i}\right\vert -\mathbb{E}\left[
\left\vert X_{i}\right\vert \right]  \right)  +R_{m}\left(  \rho\right)
\]
and%
\[
\sup_{\boldsymbol{\mu}\in\mathcal{B}_{m}\left(  \rho\right)  }\sup
_{y\in\mathbb{R}}\left\vert \partial_{y}\left(  \hat{s}_{m}\left(  y\right)
-s_{m}\left(  y\right)  \right)  \right\vert \leq\frac{1}{m}\sum_{i=1}%
^{m}\left(  \left\vert X_{i}\right\vert -\mathbb{E}\left[  \left\vert
X_{i}\right\vert \right]  \right)  +R_{m}\left(  \rho\right)  .
\]
When (\ref{eq11b}) holds, Chebyshev inequality implies%
\begin{align*}
B_{2,1} &  =\Pr\left(  \frac{1}{m}\sum_{i=1}^{m}\left(  \left\vert
X_{i}\right\vert -\mathbb{E}\left[  \left\vert X_{i}\right\vert \right]
\right)  \geq\frac{\Delta^{-1}\left(  2q\gamma_{m}\right)  ^{-1/2}\log
\gamma_{m}}{\sqrt{m}}-R_{m}\left(  \rho\right)  \right)  \\
&  \leq4A_{\mu_{0}}q\gamma_{m}m^{-2\vartheta}\left(  \log\gamma_{m}\right)
^{-2}%
\end{align*}
for all $m$ large enough, where $A_{\mu_{0}}$ is the variance of $\left\vert
X_{1}\right\vert $. Thus, for all $m$ large enough,%
\[
B_{2}\leq B_{2,1}\leq4A_{\mu_{0}}q\gamma_{m}m^{-2\vartheta}\left(  \log
\gamma_{m}\right)  ^{-2}.
\]
This, together with (\ref{eq10g}) and (\ref{eq10h}) and the continuity of
$\hat{s}_{m}\left(  y\right)  -s_{m}\left(  y\right)  $ in $y$, implies%
\[
B_{0}=\Pr\left(  \sup_{\boldsymbol{\mu}\in\mathcal{B}_{m}\left(  \rho\right)
}\sup_{y\in G_{m}}\left\vert \hat{s}_{m}\left(  y\right)  -s_{m}\left(
y\right)  \right\vert \geq\frac{\sqrt{2q\gamma_{m}}}{\sqrt{m}}\right)  \leq
p_{m}\left(  \vartheta,q,h_{\mu_{0}},\gamma_{m}\right)
\]
for all $m$ large enough. This justifies the assertion.

\textbf{Part II}: to show the uniform bound on $\left\vert \hat{\varphi}%
_{m}\left(  t,\mathbf{z}\right)  -\varphi_{m}\left(  t,\boldsymbol{\mu
}\right)  \right\vert $. Pick a positive sequence $\left\{  \tau_{m}%
:m\geq1\right\}  $ such that $\tau_{m}\leq\gamma_{m}$ for all large $m$ and
$\tau_{m}\rightarrow\infty$. Then, \textbf{Part I} implies that, with
probability at least $1-p_{m}\left(  \vartheta,q,h_{\mu_{0}},\gamma
_{m}\right)  $,
\begin{align*}
&  \sup_{\boldsymbol{\mu}\in\mathcal{B}_{m}\left(  \rho\right)  }\sup
_{t\in\left[  0,\tau_{m}\right]  }\left\vert \hat{\varphi}_{m}\left(
t,\mathbf{z}\right)  -\varphi_{m}\left(  t,\boldsymbol{\mu}\right)
\right\vert \\
&  \leq\sup_{\boldsymbol{\mu}\in\mathcal{B}_{m}\left(  \rho\right)  }%
\sup_{t\in\left[  0,\tau_{m}\right]  }\int_{\left[  -1,1\right]  }%
\omega\left(  s\right)  \frac{\sup_{t\in G_{m}}\left\vert S_{m}\left(
ts\right)  \right\vert }{r_{\mu_{0}}\left(  ts\right)  }ds\leq\Upsilon\left(
q,\tau_{m},\gamma_{m},r_{\mu_{0}}\right)
\end{align*}
for all sufficiently large $m$, where%
\[
\Upsilon\left(  q,\tau_{m},\gamma_{m},r_{\mu_{0}}\right)  =\frac{2\left\Vert
\omega\right\Vert _{\infty}\sqrt{2q\gamma_{m}}}{\sqrt{m}}\sup_{t\in\left[
0,\tau_{m}\right]  }\int_{\left[  0,1\right]  }\frac{ds}{r_{\mu_{0}}\left(
ts\right)  }.
\]

\textbf{Part III}: to determine the constants $\gamma^{\prime}$,
$\gamma^{\prime\prime}$, $q$, $\vartheta$ and $\vartheta^{\prime}$ and a
uniform consistency class. Set $\gamma^{\prime}$, $\vartheta$ and $q$ such
that $q\gamma^{\prime}>\vartheta>2^{-1}$ and $0\leq\vartheta^{\prime
}<\vartheta-1/2$. Then $p_{m}\left(  \vartheta,q,h_{\mu_{0}},\gamma
_{m}\right)  \rightarrow0$, $B_{0}\rightarrow0$ and $m^{\vartheta-1/2}%
\gamma_{m}^{-1/2}\log\gamma_{m}\rightarrow\infty$ as $m\rightarrow\infty$. If
additionally $R_{m}\left(  \rho\right)  =O\left(  m^{\vartheta^{\prime}%
}\right)  $ and $u_{m}\geq\left(  \gamma^{\prime\prime}\tau_{m}\right)
^{-1}\log\log m$. Then $\tau_{m}u_{m}\rightarrow\infty$ as $m\rightarrow
\infty$ and (\ref{eq11b}) holds.

Recall%
\[
\psi\left(  t_{m},\mu_{i};\mu_{0}\right)  =\int_{\left[  -1,1\right]  }%
\omega\left(  s\right)  \cos\left(  t_{m}s\left(  \mu_{i}-\mu_{0}\right)
\right)  ds.
\]
Since $\gamma_{m}\left\vert \mu_{i}-\mu_{0}\right\vert \geq\gamma
_{m}\left\vert u_{m}-\mu_{0}\right\vert \rightarrow\infty$ uniformly for $i\in
I_{1,m}$,%
\[
\lim_{m\rightarrow\infty}\sup\left\{  \psi\left(  t,\mu;\mu_{0}\right)
:\left(  t,\left\vert \mu\right\vert \right)  \in\lbrack\gamma_{m}%
,\infty)\times\lbrack u_{m},\infty)\right\}  =0.
\]
So, when $\pi_{1,m}^{-1}\Upsilon\left(  q,\tau_{m},\gamma_{m},r_{\mu_{0}%
}\right)  \rightarrow0$, the same reasoning used to prove (\ref{eq4a})
implies
\[
\Pr\left(  \sup\nolimits_{\boldsymbol{\mu}\in\mathcal{B}_{m}\left(
\rho\right)  }\left\vert \pi_{1,m}^{-1}\sup\nolimits_{t\in\left[  0,\tau
_{m}\right]  }\hat{\varphi}_{m}\left(  t,\mathbf{z}\right)  -1\right\vert
\rightarrow0\right)  \rightarrow1.
\]
In other words, as claimed,%
\[
\mathcal{Q}_{m}\left(  \boldsymbol{\mu},t;\mathcal{F}\right)  =\left\{
\begin{array}
[c]{c}%
q\gamma^{\prime}>\vartheta>2^{-1},\gamma^{\prime}>0,\gamma^{\prime\prime
}>0,0\leq\vartheta^{\prime}<\vartheta-1/2,\\
R_{m}\left(  \rho\right)  =O\left(  m^{\vartheta^{\prime}}\right)  ,\tau
_{m}\leq\gamma_{m},u_{m}\geq\frac{\log\log m}{\gamma^{\prime\prime}\tau_{m}%
},\\
t\in\left[  0,\tau_{m}\right]  ,\lim\limits_{m\rightarrow\infty}\pi_{1,m}%
^{-1}\Upsilon\left(  q,\tau_{m},\gamma_{m},r_{\mu_{0}}\right)  =0
\end{array}
\right\}
\]
is a uniform consistency class.

\subsection{Proof of \autoref{CorInclusion1}}

For the proof, we will refer to the proofs of \autoref{ThmLocationShift} and
\autoref{CorConI}. For Gaussian family, when $q\sigma^{-1}>\vartheta>2^{-1}$,
we can set $\gamma_{m}=\sigma^{-1}\log m$ and $\tau_{m}=\sigma^{-1}%
\sqrt{2\gamma\log m}$. Then (\ref{eqB5}) implies the claimed uniform
consistency class. Further, we see that the fastest speed of convergence is
$\sqrt{\log m}$, achieved when $\liminf_{m\rightarrow\infty}\pi_{1,m}>0$.

For Laplace family, $r_{\mu_{0}}^{-1}\left(  t\right)  =1+\sigma^{2}t^{2}$,%
\[
\sup_{t\in\left[  0,\gamma_{m}\right]  }\int_{\left[  0,1\right]  }r_{\mu_{0}%
}^{-1}\left(  ts\right)  ds=\int_{\left[  0,1\right]  }\left(  1+\sigma
^{2}\gamma_{m}^{2}s^{2}\right)  ds=1+\frac{\sigma^{2}\gamma_{m}^{2}}{3},
\]
and%
\[
\Upsilon\left(  q,\gamma_{m},\gamma_{m},r_{\mu_{0}}\right)  =\frac{2\left\Vert
\omega\right\Vert _{\infty}\sqrt{2q\gamma_{m}}}{\sqrt{m}}\left(
1+\frac{\sigma^{2}\gamma_{m}^{2}}{3}\right)  .
\]
So, setting $\tau_{m}=\gamma_{m}=\log m$ and $0\leq\gamma<1/2$ gives the
claimed uniform consistency class.

For Hyperbolic Secant family, $r_{\mu}^{-1}\left(  t\right)  =\sigma
\cosh\left(  t\sigma^{-1}\right)  \sim2^{-1}\sigma\exp\left(  \sigma
^{-1}t\right)  $ as $t\rightarrow\infty$ and%
\[
\sup_{t\in\left[  0,\gamma_{m}\right]  }\int_{\left[  0,1\right]  }r_{\mu
}^{-1}\left(  ts\right)  ds\leq C\sup_{t\in\left[  0,\gamma_{m}\right]  }%
\int_{\left[  0,1\right]  }\exp\left(  \sigma^{-1}ts\right)  ds\leq
C\exp\left(  \sigma^{-1}\gamma_{m}\right)
\]
and%
\[
\Upsilon\left(  q,\gamma_{m},\gamma_{m},r_{\mu_{0}}\right)  \leq\frac
{C\sqrt{2q\gamma_{m}}}{\sqrt{m}}\exp\left(  \sigma^{-1}\gamma_{m}\right)  .
\]
When $q\sigma>\vartheta>2^{-1}$, we can set $\gamma_{m}=\sigma\log m$ and
$\tau_{m}=\sigma\gamma\log m$ with $0<\gamma<1/2$. This gives the claimed
uniform consistency class.

For Logistic family, $r_{\mu}^{-1}\left(  t\right)  =\left(  \pi\sigma
t\right)  ^{-1}\sinh\left(  \pi\sigma t\right)  \sim\left(  2\pi\sigma
t\right)  ^{-1}e^{\pi\sigma t}$ as $t\rightarrow\infty$. Fix a small
$\varepsilon^{\prime\prime}\in\left(  0,1\right)  $. We can pick a small
$\varepsilon^{\prime}\in\left(  0,1\right)  $ such that $r_{\mu}^{-1}\left(
t\right)  \geq1-\varepsilon^{\prime\prime}$ for all $t\in\left[
0,\varepsilon^{\prime}\right]  $. Then%
\[
\sup_{t\in\left[  0,\gamma_{m}\right]  }\int_{\left[  0,1\right]  }r_{\mu
}^{-1}\left(  ts\right)  ds\leq C\left(  1+\frac{1}{\varepsilon^{\prime}}%
\sup_{t\in\left[  0,\gamma_{m}\right]  }\int_{\left[  0,1\right]  }\exp\left(
\sigma\pi ts\right)  ds\right)  \leq C\exp\left(  \sigma\pi\gamma_{m}\right)
.
\]
So,%
\[
\Upsilon\left(  q,\gamma_{m},\gamma_{m},r_{\mu_{0}}\right)  \leq\frac
{C\sqrt{\gamma_{m}}}{\sqrt{m}}\exp\left(  \sigma\pi\gamma_{m}\right)  .
\]
When $q\left(  \sigma\pi\right)  ^{-1}>\vartheta>2^{-1}$, we can set
$\gamma_{m}=\left(  \sigma\pi\right)  ^{-1}\log m$ and $\tau_{m}=\left(
\sigma\pi\right)  ^{-1}\gamma\log m$ with $0<\gamma<1/2$. This then gives the
claimed uniform consistency class.

\section{Proofs related to Construction II}

\label{ProofsII}

\subsection{Proof of \autoref{ThmDiscreteInfiniteSupport}}

Clearly, $L\left(  \theta\right)  =H\left(  e^{\theta}\right)  $ and
$c_{k}=\frac{H^{\left(  k\right)  }\left(  0\right)  }{k!}$, where $H^{\left(
k\right)  }$ is the $k$th order derivative of $H$ and $H^{\left(  0\right)
}=H$. Let%
\begin{equation}
K^{\dag}\left(  t,x;\theta_{0}\right)  =H\left(  e^{\theta_{0}}\right)
\int_{\left[  -1,1\right]  }\frac{\left(  \iota ts\right)  ^{x}}{\exp\left(
\iota tse^{\theta_{0}}\right)  H^{\left(  x\right)  }\left(  0\right)  }%
\omega\left(  s\right)  ds \label{eq13g}%
\end{equation}
and $\psi^{\dag}\left(  t,\theta;\theta_{0}\right)  =\int K^{\dag}\left(
t,x;\theta_{0}\right)  dG_{\theta}\left(  x\right)  $. Then%
\begin{align*}
\psi^{\dag}\left(  t,\theta;\theta_{0}\right)   &  =\int K^{\dag}\left(
t,x;\theta_{0}\right)  dG_{\theta}\left(  x\right) \\
&  =\frac{H\left(  e^{\theta_{0}}\right)  }{H\left(  e^{\theta}\right)  }%
\int_{\left[  -1,1\right]  }\sum_{k=0}^{\infty}\frac{\left(  \iota ts\right)
^{k}e^{\theta k}c_{k}}{\exp\left(  \iota tse^{\theta_{0}}\right)  H^{\left(
k\right)  }\left(  0\right)  }\omega\left(  s\right)  ds\\
&  =\frac{H\left(  e^{\theta_{0}}\right)  }{H\left(  e^{\theta}\right)  }%
\int_{\left[  -1,1\right]  }\exp\left(  \iota st\left(  e^{\theta}%
-e^{\theta_{0}}\right)  \right)  \omega\left(  s\right)  ds,
\end{align*}
for which $\psi^{\dag}\left(  t,\theta_{0};\theta_{0}\right)  =1$ for any $t$
and $\lim_{t\rightarrow\infty}$ $\psi^{\dag}\left(  t,\theta;\theta
_{0}\right)  =0$ for each $\theta\neq\theta_{0}$ by the RL Lemma. Taking the
real parts of $K^{\dag}$ and $\psi^{\dag}$ yields the claim.

\subsection{Proof of \autoref{LmNEFDiscrete}}

By simple calculations, we obtain the following: (1) $c_{k}k!\equiv1$ for
Poisson family; (2) $\left(  c_{k}k!\right)  ^{-1}=\frac{\left(  n-1\right)
!}{\left(  k+n-1\right)  !}$ for Negative Binomial family with a fixed $n$;
(3) $\left(  c_{k}k!\right)  ^{-1}=\left(  1+k\right)  ^{-\left(  k-1\right)
}$ for Abel family; (4) $\left(  c_{k}k!\right)  ^{-1}=\left(  k+1\right)
!\left(  \left(  2k\right)  !\right)  ^{-1}$ for Tak\'{a}cs family. Therefore,
(\ref{eqUPde}) holds. Fix a $\sigma>0$. Then for Strict Arcsine family,%
\[
c_{k}k!=c_{k}^{\ast}\left(  1\right)  \geq2^{k-2}\left(  \left(  \lfloor
2^{-1}k\rfloor-1\right)  !\right)  ^{2},
\]
and for Large Arcsine family,
\[
c_{k}k!=\frac{k+1}{c_{k}^{\ast}\left(  1+k\right)  }\geq\left(  1+k\right)
^{\lfloor2^{-1}k\rfloor-1}2^{k-2}\left(  \left(  \lfloor2^{-1}k\rfloor
-1\right)  !\right)  ^{2}.
\]
So, (\ref{eqUPde}) does not hold for these two families.

Now we show the third claim. Since $H\left(  z\right)  =\sum_{k=0}^{\infty
}c_{k}z^{k}$ has a positive radius $R_{H}$ of convergence, there exists
$R_{H}>\tilde{r}>0$ such that%
\[
H^{\left(  k\right)  }\left(  0\right)  =\frac{k!}{2\pi\iota}\int_{\left\{
z\in\mathbb{C}:\left\vert z\right\vert =\tilde{r}\right\}  }\frac{H\left(
z\right)  }{z^{k+1}}dz\text{ \ \ \ for all }k\in\mathbb{N}.
\]
However, $H$ has all positive coefficients. Therefore, $\max_{\left\{
z\in\mathbb{C}:\left\vert z\right\vert =\tilde{r}\right\}  }\left\vert
H\left(  z\right)  \right\vert $ is achieved when $z=\tilde{r}$, and%
\[
\left\vert H^{\left(  k\right)  }\left(  0\right)  \right\vert \leq\frac
{k!}{2\pi}2\pi\tilde{r}\frac{\sup_{\left\vert z\right\vert =\tilde{r}%
}\left\vert H\left(  z\right)  \right\vert }{\tilde{r}^{k+1}}=H\left(
\tilde{r}\right)  \frac{k!}{\tilde{r}^{k}}.
\]
Observing that $H^{\left(  k\right)  }\left(  0\right)  $ is real and
$H^{\left(  k\right)  }\left(  0\right)  =c_{k}k!$ for $k\in\mathbb{N}$ gives
(\ref{eqLBder}).

\subsection{Proof of \autoref{II-concentration}}

First, we prove the following lemma.

\begin{lemma}
If $Z$ has CDF $G_{\theta}$ with GF $H$, (\ref{eqUPde}) holds and $\eta>0$,
then
\begin{equation}
\mathbb{E}\left[  t^{2Z}\left(  H^{\left(  Z\right)  }\left(  0\right)
\right)  ^{-2}\right]  \leq\frac{C}{L\left(  \theta\right)  }\frac{\exp\left(
2t\sqrt{\eta}\right)  }{\sqrt{t\sqrt{\eta}}} \label{eq15a}%
\end{equation}
for positive and sufficiently large $t$.
\end{lemma}

\begin{proof}
Recall $H\left(  z\right)  =\sum_{k=0}^{\infty}c_{k}z^{k}$ and $H^{\left(
k\right)  }\left(  0\right)  =c_{k}k!$ for $k\in\mathbb{N}$. Let $\chi
_{Z}\left(  t\right)  =\mathbb{E}\left[  t^{2Z}\left(  H^{\left(  Z\right)
}\left(  0\right)  \right)  ^{-2}\right]  $. Since (\ref{eqUPde}) holds,%
\[
\chi_{Z}\left(  t\right)  =\frac{1}{L\left(  \theta\right)  }\sum
_{k=0}^{\infty}\frac{t^{2k}c_{k}\eta^{k}}{\left(  c_{k}k!\right)  ^{2}}%
\leq\frac{C}{L\left(  \theta\right)  }B_{\mathrm{II}}\left(  2t\sqrt{\eta
}\right)  ,
\]
where $B_{\mathrm{II}}\left(  x\right)  =\sum_{k=0}^{\infty}\frac{\left(
x/2\right)  ^{2k}}{\left(  k!\right)  ^{2}}$ for $x>0$. So, it suffices to
bound $B_{\mathrm{II}}\left(  x\right)  $. For $\sigma>0$ and $y\in\mathbb{C}%
$, let%
\begin{equation}
J_{\sigma}\left(  y\right)  =\sum_{n=0}^{\infty}\frac{\left(  -1\right)  ^{n}%
}{n!\Gamma\left(  n+\sigma+1\right)  }\left(  \frac{y}{2}\right)  ^{2n+\sigma}
\label{BesselA}%
\end{equation}
Then $J_{\sigma}$ is the Bessel function of the first kind of order $\sigma$;
see definition (1.17.1) in Chapter 1 of \cite{Szego:1975}, and $B_{\mathrm{II}%
}\left(  x\right)  =J_{0}\left(  -\iota x\right)  $. By identity (1.71.8) in
Chapter 1 of \cite{Szego:1975} that was derived on page 368 of
\cite{Whittaker:1940}, we have%
\begin{equation}
J_{\sigma}\left(  y\right)  =\sqrt{2}\left(  \pi y\right)  ^{-1/2}\cos\left(
y-c_{0}\right)  \left(  1+O\left(  y^{-2}\right)  \right)  \label{BesselB}%
\end{equation}
\ as $y\rightarrow\infty$ whenever\ $\left\vert \arg y\right\vert <\pi$, where
$c_{0}=2^{-1}\sigma\pi-4^{-1}\pi$. So,%
\begin{equation}
B_{\mathrm{II}}\left(  2t\sqrt{\eta}\right)  =\left\vert J_{0}\left(  -2\iota
t\sqrt{\eta}\right)  \right\vert =C\left(  t\sqrt{\eta}\right)  ^{-1/2}%
\exp\left(  2t\sqrt{\eta}\right)  \label{BesselC}%
\end{equation}
as $0<t\eta\rightarrow\infty$, where we have used the identity $\left\vert
\cos z\right\vert ^{2}=\cosh^{2}\left(  \Im\left(  z\right)  \right)
-\sin^{2}\left(  \Re\left(  z\right)  \right)  $. The bound given by
(\ref{BesselC}) is tight up to a multiple of a positive constant, which can be
seen from inequality (5) of \cite{Gronwall:1932}. On the other hand,
$B_{\mathrm{II}}\left(  t\eta\right)  =O\left(  1\right)  $\ when
$t\eta=O\left(  1\right)  $. Thus, when $\eta>0$ and (\ref{eqUPde}) holds,
(\ref{eq15a}) holds for all positive and sufficiently large $t$.
\end{proof}

Now we show the first claim of the theorem. Recall $\eta_{i}=e^{\theta_{i}}$
when $z_{i}$ has CDF $G_{\theta_{i}}$ for $1\leq i\leq m$. Let $V_{m}\left(
\hat{\varphi}\right)  =\mathbb{V}\left[  \hat{\varphi}_{m}\left(
t,\mathbf{z}\right)  -\varphi_{m}\left(  t,\boldsymbol{\theta}\right)
\right]  $ and $\tilde{p}_{m}\left(  \lambda\right)  =\Pr\left(  \left\vert
\hat{\varphi}_{m}\left(  t,\mathbf{z}\right)  -\varphi_{m}\left(
t,\boldsymbol{\theta}\right)  \right\vert \geq\lambda\right)  $. Let%
\begin{equation}
w\left(  t,x\right)  =\frac{t^{x}\cos\left(  \frac{\pi}{2}x-t\eta_{0}\right)
}{c_{x}x!}\text{ }\ \ \text{for \ }t\geq0\text{ and }x\in\mathbb{N}
\label{eq12h}%
\end{equation}
and
\[
S_{m}\left(  t\right)  =\frac{1}{m}\sum_{i=1}^{m}\left(  w\left(
t,z_{i}\right)  -\mathbb{E}\left[  w\left(  t,z_{i}\right)  \right]  \right)
\text{ \ for \ }t\geq0.\text{ }%
\]
Then (\ref{II-c}) is equivalent to $K\left(  t,x;\theta_{0}\right)  =H\left(
\eta_{0}\right)  \int_{\left[  -1,1\right]  }w\left(  ts,x\right)
\omega\left(  s\right)  ds$ and%
\[
\hat{\varphi}_{m}\left(  t,\mathbf{z}\right)  -\varphi_{m}\left(
t,\boldsymbol{\theta}\right)  =H\left(  \eta_{0}\right)  \int_{\left[
-1,1\right]  }S_{m}\left(  ts\right)  \omega\left(  s\right)  ds.
\]
Since $\mathbb{E}\left[  w^{2}\left(  t,Z\right)  \right]  \leq\chi_{Z}\left(
t\right)  $, inequality (\ref{eq15a}) implies%
\[
V_{m}\left(  \hat{\varphi}\right)  \leq\frac{C}{m^{2}}\sum_{i=1}^{m}\frac
{\exp\left(  2t\sqrt{\eta_{i}}\right)  }{L\left(  \theta_{i}\right)  \left(
t\sqrt{\eta_{i}}\right)  ^{1/2}}\leq\frac{C}{m}\frac{\exp\left(  2t\left\Vert
\boldsymbol{\eta}\right\Vert _{\infty}^{1/2}\right)  }{\sqrt{t}\phi_{m}\left(
L,\boldsymbol{\theta}\right)  }.
\]
So,%
\[
\Pr\left(  \left\vert S_{m}\left(  t\right)  \right\vert \geq\lambda\right)
\leq\frac{CV_{m}\left(  \hat{\varphi}\right)  }{\lambda^{2}}\leq\frac
{C}{\lambda^{2}m}\frac{\exp\left(  2t\left\Vert \boldsymbol{\eta}\right\Vert
_{\infty}^{1/2}\right)  }{\sqrt{t}\phi_{m}\left(  L,\boldsymbol{\theta
}\right)  }.
\]
From%
\[
\int_{\left[  -1,1\right]  }S_{m}\left(  ts\right)  \omega\left(  s\right)
ds=\frac{1}{t}\int_{\left[  -t,t\right]  }S_{m}\left(  y\right)  \omega\left(
yt^{-1}\right)  dy
\]
for $t>0$, we have%
\[
\tilde{p}_{m}\left(  \lambda\right)  \leq\Pr\left(  \left\vert S_{m}\left(
t\right)  \right\vert \geq\frac{\lambda}{2H\left(  \eta_{0}\right)  \left\Vert
\omega\right\Vert _{\infty}}\right)  \leq\frac{CH^{2}\left(  \eta_{0}\right)
\left\Vert \omega\right\Vert _{\infty}^{2}}{\lambda^{2}m}\frac{\exp\left(
2t\left\Vert \boldsymbol{\eta}\right\Vert _{\infty}^{1/2}\right)  }{\sqrt
{t}\phi_{m}\left(  L,\boldsymbol{\theta}\right)  }.
\]

Finally, we show the second claim of the theorem. For Poisson family,
$c_{k}k!=1$ for all $k\in\mathbb{N}$. So, from (\ref{eq12h}) we obtain%
\[
\mathbb{E}\left[  w^{2}\left(  t,Z\right)  \right]  \leq\frac{H^{2}\left(
\eta_{0}\right)  \exp\left(  t^{2}\sqrt{\eta}\right)  }{L\left(
\theta\right)  }\text{ \ and \ \ }V_{m}\left(  \hat{\varphi}\right)  \leq
\frac{C}{m}\frac{\exp\left(  t^{2}\left\Vert \boldsymbol{\eta}\right\Vert
_{\infty}^{1/2}\right)  }{\min_{1\leq i\leq m}L\left(  \theta_{i}\right)  }%
\]
and%
\[
\tilde{p}_{m}\left(  \lambda\right)  \leq\Pr\left(  \left\vert S_{m}\left(
t\right)  \right\vert \geq\frac{\lambda}{2H\left(  \eta_{0}\right)  \left\Vert
\omega\right\Vert _{\infty}}\right)  \leq\frac{C}{\lambda^{2}m}\frac
{\exp\left(  t^{2}\left\Vert \boldsymbol{\eta}\right\Vert _{\infty}%
^{1/2}\right)  }{\min_{1\leq i\leq m}L\left(  \theta_{i}\right)  }%
\]
for positive and sufficiently large $t$.

\subsection{Proof of \autoref{II-consistency}}

Obviously, $\phi_{m}\left(  L,\boldsymbol{\theta}\right)  $ is positive and
finite when $\left\Vert \boldsymbol{\theta}\right\Vert _{\infty}\leq\rho$.
First, consider the case when (\ref{eqUPde}) holds. Then (\ref{eq15e4})
implies%
\[
\Pr\left(  \left\vert \hat{\varphi}_{m}\left(  t,\mathbf{z}\right)
-\varphi_{m}\left(  t,\boldsymbol{\theta}\right)  \right\vert \geq
\lambda\right)  \leq\frac{C}{\lambda^{2}m}\frac{\exp\left(  2t\left\Vert
\boldsymbol{\eta}\right\Vert _{\infty}^{1/2}\right)  }{\sqrt{t}}.
\]
We can set $t=2^{-1}\left\Vert \boldsymbol{\eta}\right\Vert _{\infty}%
^{-1/2}\gamma\log m$ for $\gamma\in\left(  0,1\right]  $, which induces%
\[
\Pr\left(  \left\vert \hat{\varphi}_{m}\left(  t,\mathbf{z}\right)
-\varphi_{m}\left(  t,\boldsymbol{\theta}\right)  \right\vert \geq
\lambda\right)  \leq\frac{C}{\lambda^{2}m^{1-\gamma}}\frac{1}{\sqrt{\gamma\log
m}}.
\]
Let $\varepsilon>0$ be any finite constant. If $\pi_{1,m}\geq Cm^{\left(
\gamma-1\right)  /2}$, then
\[
\Pr\left(  \frac{\left\vert \hat{\varphi}_{m}\left(  t,\mathbf{z}\right)
-\varphi_{m}\left(  t,\boldsymbol{\theta}\right)  \right\vert }{\pi_{1,m}}%
\geq\varepsilon\right)  \leq\frac{C\varepsilon^{-2}}{\sqrt{\gamma\log m}%
}\rightarrow0\text{ \ as \ }m\rightarrow\infty.
\]
Moreover,
\[
\psi\left(  t,\theta;\theta_{0}\right)  =\frac{H\left(  \eta_{0}\right)
}{H\left(  \eta\right)  }\int_{\left[  -1,1\right]  }\cos\left(  st\left(
\eta-\eta_{0}\right)  \right)  \omega\left(  s\right)  ds\rightarrow0\text{
\ as \ }m\rightarrow\infty
\]
whenever $\lim_{m\rightarrow\infty}t\min_{1\leq i\leq m}\left\vert \eta
_{0}-\eta_{i}\right\vert =\infty$.

Secondly, we deal with Poisson family. Clearly, $L_{\min}^{\left(  m\right)
}=\min_{1\leq i\leq m}L\left(  \theta_{i}\right)  $ is positive and finite
when $\left\Vert \boldsymbol{\theta}\right\Vert _{\infty}\leq\rho$. So,
inequality (\ref{eqII3}) implies%
\[
\Pr\left(  \left\vert \hat{\varphi}_{m}\left(  t,\mathbf{z}\right)
-\varphi_{m}\left(  t,\boldsymbol{\theta}\right)  \right\vert \geq
\lambda\right)  \leq\frac{C}{\lambda^{2}m}\exp\left(  t^{2}\left\Vert
\boldsymbol{\eta}\right\Vert _{\infty}^{1/2}\right)  .
\]
So, we can set $t=\sqrt{\left\Vert \boldsymbol{\eta}\right\Vert _{\infty
}^{-1/2}\gamma\log m}$ for a fixed $\gamma\in\left(  0,1\right)  $, which
induces%
\[
\Pr\left(  \left\vert \hat{\varphi}_{m}\left(  t,\mathbf{z}\right)
-\varphi_{m}\left(  t,\boldsymbol{\theta}\right)  \right\vert \geq
\lambda\right)  \leq\frac{C}{\lambda^{2}m^{1-\gamma}}.
\]
If $\pi_{1,m}\geq Cm^{\left(  \gamma^{\prime}-1\right)  /2}$ for any
$\gamma^{\prime}>\gamma$, then%
\[
\Pr\left(  \frac{\left\vert \hat{\varphi}_{m}\left(  t,\mathbf{z}\right)
-\varphi_{m}\left(  t,\boldsymbol{\theta}\right)  \right\vert }{\pi_{1,m}}%
\geq\varepsilon\right)  \leq\frac{C\varepsilon^{-2}}{m^{\gamma^{\prime}%
-\gamma}}\rightarrow0\text{ \ as \ }m\rightarrow\infty.
\]
This completes the proof.

\section{Proofs related to Construction III}

\label{ProofsIII}

\subsection{Proof of \autoref{ThmConstructionMoments}}

Let%
\begin{equation}
K^{\dag}\left(  t,x;\mu_{0}\right)  =\int_{\left[  -1,1\right]  }\exp\left(
\iota ts\xi\left(  \theta_{0}\right)  \right)  \sum_{n=0}^{\infty}%
\frac{\left(  -\iota tsx\right)  ^{n}}{\tilde{a}_{n}n!}\omega\left(  s\right)
ds. \label{eq1i}%
\end{equation}
Since $\left\{  \tilde{c}_{n}\left(  \theta\right)  \right\}  _{n\geq1}$ is
separable at $\theta_{0}$, then%
\begin{align}
\psi^{\dag}\left(  t,\mu;\mu_{0}\right)   &  =\frac{1}{\zeta\left(  \theta
_{0}\right)  }\int K^{\dag}\left(  t,x;\theta_{0}\right)  dG_{\theta}\left(
x\right) \nonumber\\
&  =\frac{1}{\zeta\left(  \theta_{0}\right)  }\int dG_{\theta}\left(
x\right)  \int_{\left[  -1,1\right]  }\exp\left(  \iota ts\xi\left(
\theta_{0}\right)  \right)  \sum_{n=0}^{\infty}\frac{\left(  -\iota
tsx\right)  ^{n}}{\tilde{a}_{n}n!}\omega\left(  s\right)  ds\nonumber\\
&  =\frac{1}{\zeta\left(  \theta_{0}\right)  }\int_{\left[  -1,1\right]  }%
\exp\left(  \iota ts\xi\left(  \theta_{0}\right)  \right)  \omega\left(
s\right)  ds\sum_{n=0}^{\infty}\frac{\left(  -\iota ts\right)  ^{n}}{\tilde
{a}_{n}n!}\tilde{c}_{n}\left(  \theta\right) \nonumber\\
&  =\frac{\zeta\left(  \theta\right)  }{\zeta\left(  \theta_{0}\right)  }%
\int_{\left[  -1,1\right]  }\exp\left(  \iota ts\xi\left(  \theta_{0}\right)
\right)  \omega\left(  s\right)  ds\sum_{n=0}^{\infty}\frac{\left(  -\iota
ts\right)  ^{n}}{n!}\xi^{n}\left(  \theta\right) \nonumber\\
&  =\frac{\zeta\left(  \theta\right)  }{\zeta\left(  \theta_{0}\right)  }%
\int_{\left[  -1,1\right]  }\exp\left(  \iota ts\left(  \xi\left(  \theta
_{0}\right)  -\xi\left(  \theta\right)  \right)  \right)  \omega\left(
s\right)  ds. \label{eq1j}%
\end{align}
Further, $\psi^{\dag}\left(  t,\mu;\mu_{0}\right)  =1$ when $\mu=\mu_{0}$ for
all $t$, and the RL Lemma implies that $\lim_{t\rightarrow\infty}$ $\psi
^{\dag}\left(  t,\mu;\mu_{0}\right)  =0$ for each $\theta\neq\theta_{0}$.
Taking the real parts of $K^{\dag}$ and $\psi^{\dag}$ gives the claim.

\subsection{Proof of \autoref{ConcentrationIII}}

First, we prove the following lemma.

\begin{lemma}
For a fixed $\sigma>0$, let%
\[
\tilde{w}\left(  z,x\right)  =\sum_{n=0}^{\infty}\frac{\left(  zx\right)
^{n}}{n!\Gamma\left(  \sigma+n\right)  }\text{ \ for }z,x>0.
\]
If $Z$ has CDF $G_{\theta}$ from the Gamma family with scale parameter
$\sigma$, then%
\begin{equation}
\mathbb{E}\left[  \tilde{w}^{2}\left(  z,Z\right)  \right]  \leq C\left(
\frac{z}{1-\theta}\right)  ^{3/4-\sigma}\exp\left(  \frac{4z}{1-\theta
}\right)  \label{eq15}%
\end{equation}
for positive and sufficiently large $z$.
\end{lemma}

\begin{proof}
Recall the Bessel function $J_{\sigma}$ defined by (\ref{BesselA}) and the
asymptotic bound (\ref{BesselB}). Then,%
\[
\tilde{w}\left(  z,x\right)  =\left(  \sqrt{zx}\right)  ^{1-\sigma}\left\vert
\iota^{1-\sigma}J_{\sigma-1}\left(  \iota2\sqrt{zx}\right)  \right\vert
=\left(  zx\right)  ^{\frac{1}{4}-\frac{\sigma}{2}}\exp\left(  2\sqrt
{zx}\right)  \left(  1+O\left(  \left(  zx\right)  ^{-1}\right)  \right)
\]
when $zx\rightarrow\infty$. Let $A_{1,z}=\left\{  x\in\left(  0,\infty\right)
:zx=O\left(  1\right)  \right\}  $. Then, on the set $A_{1,z}$, $f_{\theta
}\left(  x\right)  =O\left(  x^{\sigma-1}\right)  $ and $\tilde{w}\left(
z,x\right)  \leq Ce^{zx}=O\left(  1\right)  $ when $\theta<1$. Therefore,%
\begin{equation}
\int_{A_{1,z}}\tilde{w}^{2}\left(  z,x\right)  dG_{\theta}\left(  x\right)
\leq C\left(  1-\theta\right)  ^{\sigma}\int_{A_{1,z}}x^{\sigma-1}dx\leq
C\left(  1-\theta\right)  ^{\sigma}z^{-\sigma}. \label{eq16f}%
\end{equation}
On the other hand, let $A_{2,z}=\left\{  x\in\left(  0,\infty\right)
:\lim_{z\rightarrow\infty}zx=\infty\right\}  $. Then%
\begin{align}
\int_{A_{2,z}}\tilde{w}^{2}\left(  z,x\right)  dG_{\theta}\left(  x\right)
&  \leq C\int_{A_{2,z}}\left(  zx\right)  ^{\frac{1}{2}-\sigma}\exp\left(
4\sqrt{zx}\right)  dG_{\theta}\left(  x\right) \nonumber\\
&  \leq\int_{A_{2,z}}\left(  zx\right)  ^{\frac{1}{2}-\sigma}\sum
_{n=0}^{\infty}\frac{\left(  4\sqrt{zx}\right)  ^{n}}{n!}dG_{\theta}\left(
x\right)  =z^{\frac{1}{2}-\sigma}B_{\mathrm{III}}\left(  z\right)  ,
\label{eq16a}%
\end{align}
where%
\[
B_{\mathrm{III}}\left(  z\right)  =\sum_{n=0}^{\infty}\frac{4^{n}z^{n/2}}%
{n!}c_{n/2}^{\dagger}\text{ \ \ and \ }c_{n/2}^{\dagger}=\int x^{2^{-1}\left(
n+1\right)  -\sigma}dG_{\theta}\left(  x\right)
\]
and $c_{n/2}^{\dagger}$ is referred to as a \textquotedblleft
half-moment\textquotedblright.

However,%
\[
c_{n/2}^{\dagger}=\frac{\left(  1-\theta\right)  ^{\sigma}}{\Gamma\left(
\sigma\right)  }\int_{0}^{\infty}x^{2^{-1}\left(  n+1\right)  -\sigma
}e^{\theta x}x^{\sigma-1}e^{-x}dx=\frac{\Gamma\left(  2^{-1}n+2^{-1}\right)
}{\Gamma\left(  \sigma\right)  }\frac{\left(  1-\theta\right)  ^{\sigma-1/2}%
}{\left(  1-\theta\right)  ^{n/2}},
\]
and by Stirling formula,%
\begin{align*}
\frac{\Gamma\left(  2^{-1}n+2^{-1}\right)  }{n!}  &  \leq C\frac{\sqrt
{\pi\left(  n-1\right)  }\left(  \frac{n-1}{2}\right)  ^{\frac{n-1}{2}}%
}{e^{\frac{n-1}{2}}\sqrt{2\pi n}\left(  \frac{n}{e}\right)  ^{n}}\leq
Ce^{\frac{n}{2}}2^{-\frac{n}{2}}\frac{\left(  n-1\right)  ^{n/2}}{n^{n/2}%
}\frac{\left(  n-1\right)  ^{-1/2}}{n^{n/2}}\\
&  \leq Ce^{\frac{n}{2}}2^{-\frac{n}{2}}\frac{\left(  n-1\right)  ^{-1/2}%
}{n^{n/2}}\leq C2^{-\frac{n}{2}}\frac{n^{-1/4}}{\sqrt{n!}}.
\end{align*}
Therefore,%
\begin{equation}
B_{\mathrm{III}}\left(  z\right)  \leq C\left(  1-\theta\right)  ^{\sigma
-1/2}\sum_{n=0}^{\infty}\frac{4^{n}z^{n/2}2^{-\frac{n}{2}}}{\left(
1-\theta\right)  ^{n/2}}\frac{1}{\sqrt{n!}}=C\left(  1-\theta\right)
^{\sigma-1/2}Q^{\ast}\left(  \frac{8z}{1-\theta}\right)  , \label{eq16}%
\end{equation}
where $Q^{\ast}\left(  z\right)  =\sum_{n=0}^{\infty}\frac{z^{n/2}}{\sqrt{n!}%
}$. By definition (8.01) and identity (8.07) in Chapter 8 of \cite{Olver:1974}%
,%
\begin{equation}
Q^{\ast}\left(  z\right)  =\sqrt{2}\left(  2\pi z\right)  ^{1/4}\exp\left(
2^{-1}z\right)  \left(  1+O\left(  z^{-1}\right)  \right)  . \label{eq16c}%
\end{equation}
Combining (\ref{eq16a}), (\ref{eq16}) and (\ref{eq16c}) gives%
\[
\int_{A_{2,z}}\tilde{w}^{2}\left(  z,x\right)  dG_{\theta}\left(  x\right)
\leq C\left(  1-\theta\right)  ^{\sigma-1/2}z^{\frac{1}{2}-\sigma}\left(
\frac{z}{1-\theta}\right)  ^{1/4}\exp\left(  \frac{4z}{1-\theta}\right)
\]
for all positive and sufficiently large $z$. Recall (\ref{eq16f}). Thus, when
$1-\theta>0$, $\sigma>0$ and $z$ is positive and sufficiently large,
\begin{align*}
\mathbb{E}\left[  \tilde{w}^{2}\left(  z,Z\right)  \right]   &  \leq
\int_{A_{1,z}}\tilde{w}^{2}\left(  z,x\right)  dG_{\theta}\left(  x\right)
+\int_{A_{2,z}}\tilde{w}^{2}\left(  z,x\right)  dG_{\theta}\left(  x\right) \\
&  \leq C\left(  z^{-\sigma}+\left(  \frac{z}{1-\theta}\right)  ^{3/4-\sigma
}\exp\left(  \frac{4z}{1-\theta}\right)  \right) \\
&  \leq C\left(  \frac{z}{1-\theta}\right)  ^{3/4-\sigma}\exp\left(  \frac
{4z}{1-\theta}\right)  .
\end{align*}
Thus, (\ref{eq15}) holds.
\end{proof}

Now we show the theorem. Define%
\[
w\left(  t,x\right)  =\Gamma\left(  \sigma\right)  \sum_{n=0}^{\infty}%
\frac{\left(  -tx\right)  ^{n}\cos\left(  2^{-1}\pi n+t\xi\left(  \theta
_{0}\right)  \right)  }{n!\Gamma\left(  n+\sigma\right)  }\text{ \ for }%
t\geq0\text{ and }x>0.
\]
Set $S_{m}\left(  t\right)  =m^{-1}\sum_{i=1}^{m}\left(  w\left(
t,z_{i}\right)  -\mathbb{E}\left[  w\left(  t,z_{i}\right)  \right]  \right)
$. Then%
\[
K\left(  t,x;\theta_{0}\right)  =\int_{\left[  -1,1\right]  }w\left(
ts,x\right)  \omega\left(  s\right)  ds.
\]
Recall $V_{m}\left(  \hat{\varphi}\right)  =\mathbb{V}\left[  \hat{\varphi
}_{m}\left(  t,\mathbf{z}\right)  -\varphi_{m}\left(  t,\boldsymbol{\theta
}\right)  \right]  $. Since $\left\vert w\left(  t,x\right)  \right\vert
\leq\Gamma\left(  \sigma\right)  \tilde{w}\left(  t,x\right)  $ uniformly in
$\left(  t,x\right)  $, (\ref{eq15}) implies, for positive and sufficiently
large $t$,%
\begin{align*}
V_{m}\left(  \hat{\varphi}\right)   &  \leq m^{-2}\Gamma^{2}\left(
\sigma\right)  \sum_{i=1}^{m}\mathbb{E}\left[  \tilde{w}^{2}\left(
t,x\right)  \right] \\
&  \leq\frac{C}{m^{2}}\sum_{i=1}^{m}\left(  \frac{t}{1-\theta_{i}}\right)
^{3/4-\sigma}\exp\left(  \frac{4t}{1-\theta_{i}}\right)  \leq\frac{1}%
{m}V_{\mathrm{III}}^{\left(  m\right)  },
\end{align*}
where%
\[
V_{\mathrm{III}}^{\left(  m\right)  }=\frac{C}{m}\exp\left(  \frac{4t}%
{u_{3,m}}\right)  \sum_{i=1}^{m}\left(  \frac{t}{1-\theta_{i}}\right)
^{3/4-\sigma}.
\]
and $u_{3,m}=\min_{1\leq i\leq m}\left\{  1-\theta_{i}\right\}  $. Recall
$\tilde{p}_{m}\left(  \lambda\right)  =\Pr\left(  \left\vert \hat{\varphi}%
_{m}\left(  t,\mathbf{z}\right)  -\varphi_{m}\left(  t,\boldsymbol{\theta
}\right)  \right\vert \geq\lambda\right)  $. So, $\Pr\left(  \left\vert
S_{m}\left(  t\right)  \right\vert \geq\lambda\right)  \leq\lambda
^{-2}V_{\mathrm{III}}^{\left(  m\right)  }$ and%
\[
\tilde{p}_{m}\left(  \lambda\right)  \leq\Pr\left(  \left\vert S_{m}\left(
t\right)  \right\vert \geq\frac{\zeta\left(  \theta_{0}\right)  \lambda
}{2\left\Vert \omega\right\Vert _{\infty}}\right)  \leq\frac{\left\Vert
\omega\right\Vert _{\infty}^{2}}{\zeta^{2}\left(  \theta_{0}\right)  }\frac
{1}{\lambda^{2}}V_{\mathrm{III}}^{\left(  m\right)  },
\]
completing the proof.

\subsection{Proof of \autoref{III-uniformConsistent}}

Recall%
\[
V_{\mathrm{III}}^{\left(  m\right)  }=\frac{C}{m}\exp\left(  \frac{4t}%
{u_{3,m}}\right)  \sum_{i=1}^{m}\left(  \frac{t}{1-\theta_{i}}\right)
^{3/4-\sigma}%
\]
and (\ref{eq18}), i.e.,%
\[
\Pr\left(  \left\vert \hat{\varphi}_{m}\left(  t,\mathbf{z}\right)
-\varphi_{m}\left(  t,\boldsymbol{\mu}\right)  \right\vert \geq\lambda\right)
\leq\frac{1}{\lambda^{2}}V_{\mathrm{III}}^{\left(  m\right)  },
\]
where $u_{3,m}=\min_{1\leq i\leq m}\left\{  1-\theta_{i}\right\}  $ and
$\theta_{i}<1$, we divide the rest of the proof into two cases: $\sigma>3/4$
or $\sigma\leq3/4$. If $\sigma>3/4$ and $\left\Vert \boldsymbol{\theta
}\right\Vert_{\infty} \leq\rho$, then%
\[
V_{\mathrm{III}}^{\left(  m\right)  }\leq C t^{3/4-\sigma}\exp\left(
\frac{4t}{u_{3,m}}\right)  \left(  \max_{1\leq i\leq m}\left(  1-\theta
_{i}\right)  \right)  ^{\sigma-3/4}\leq C t^{3/4-\sigma}\exp\left(
\frac{4t}{u_{3,m}}\right)  .
\]
So, we can set $t=4^{-1}u_{3,m}\gamma\log m$ for any fixed $\gamma\in\left(
0,1\right]  $ to obtain%
\begin{equation}
\Pr\left(  \left\vert \hat{\varphi}_{m}\left(  t,\mathbf{z}\right)
-\varphi_{m}\left(  t,\boldsymbol{\mu}\right)  \right\vert \geq\lambda\right)
\leq\frac{C}{m^{1-\gamma}\lambda^{2}}\left(  4^{-1}u_{3,m}\gamma\log m\right)
^{3/4-\sigma}, \label{eq18c}%
\end{equation}
which implies%
\[
\Pr\left(  \frac{\left\vert \hat{\varphi}_{m}\left(  t,\mathbf{z}\right)
-\varphi_{m}\left(  t,\boldsymbol{\mu}\right)  \right\vert }{\pi_{1,m}}%
\geq\varepsilon\right)  \leq C\varepsilon^{-2}\left(  4^{-1}u_{3,m}\gamma\log
m\right)  ^{3/4-\sigma}\rightarrow0\text{ as }m\rightarrow\infty
\]
for any fixed $\varepsilon>0$ whenever $\pi_{1,m}\geq Cm^{\left(
\gamma-1\right)  /2}$ and $u_{3,m}\gamma\log m\rightarrow\infty$.

In contrast, if $\sigma\leq3/4$, then%
\[
V_{\mathrm{III}}^{\left(  m\right)  }\leq C \left(  \frac{t}{u_{3,m}%
}\right)  ^{3/4-\sigma}\exp\left(  \frac{4t}{u_{3,m}}\right)  .
\]
So, we can still set $t=4^{-1}u_{3,m}\gamma\log m$ for any fixed $\gamma
\in\left(  0,1\right)  $, which implies%
\[
\Pr\left(  \frac{\left\vert \hat{\varphi}_{m}\left(  t,\mathbf{z}\right)
-\varphi_{m}\left(  t,\boldsymbol{\mu}\right)  \right\vert }{\pi_{1,m}}%
\geq\varepsilon\right)  \leq\frac{C\varepsilon^{-2}\left(  4^{-1}\gamma\log
m\right)  ^{3/4-\sigma}}{m^{\gamma^{\prime}-\gamma}}\rightarrow0\text{ as
}m\rightarrow\infty
\]
whenever $\pi_{1,m}\geq Cm^{\left(  \gamma^{\prime}-1\right)  /2}$ for any
$\gamma^{\prime}>\gamma$. Recall
\[
\psi\left(  t,\mu;\mu_{0}\right)  =\int_{\left[  -1,1\right]  }\cos\left(
ts\left(  \xi\left(  \theta_{0}\right)  -\xi\left(  \theta\right)  \right)
\right)  \omega\left(  s\right)  ds,
\]
which converges to $0$ as $m\rightarrow\infty$ when $\lim_{m\rightarrow\infty
}t\min_{i\in I_{1,m}}\left\vert \xi\left(  \theta_{0}\right)  -\xi\left(
\theta_{i}\right)  \right\vert =\infty$. Noticing $u_{3,m}=\min_{1\leq i\leq
}\xi^{-1}\left(  \theta_{i}\right)  $, we have shown the claim.

\newpage
\section{Additional simulation results}
\label{secSimAdd}

This section presents the performances of the hybrid estimator ``Jin'' induced by the estimator of \cite{Jin:2008}, the ``MR'' estimator of \cite{Meinshausen:2006}, and the proposed estimator ``New'' when they are applied to each of the five families, i.e., Cauchy, Laplace, Poisson, Negative Binomial and Gamma families.

\begin{figure}[H]
\centering
\includegraphics[height=0.8\textheight]{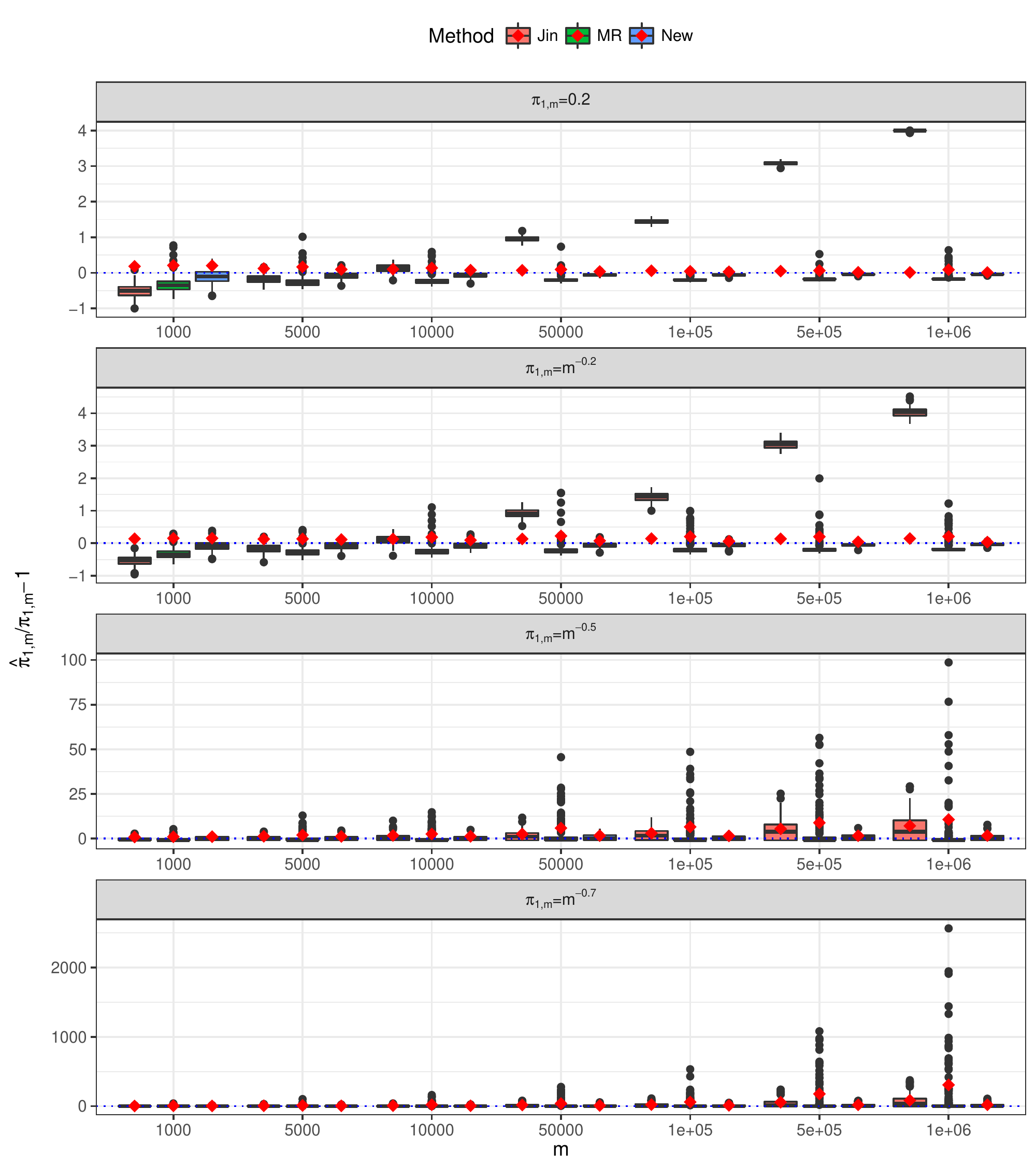}
\caption[Cauchy]{Boxplot of the excess $\tilde{\delta}_{m}=\hat{\pi}_{1,m}\pi_{1,m}^{-1}-1$ of an estimator $\hat{\pi}_{1,m}$ when it is applied to Cauchy family. The boxplot for the proposed estimator is the right one in each triple of boxplots for each $m$. The thick horizontal line and the diamond in each boxplot are respectively the mean and standard deviation of $\tilde{\delta}_{m}$, and the dotted horizontal line in each panel corresponding to a setting of $\pi_{1,m}$ is the reference for $\tilde{\delta}_{m}=0$.}
\label{figCauchy}%
\end{figure}

\begin{figure}[H]
\centering
\includegraphics[height=0.8\textheight]{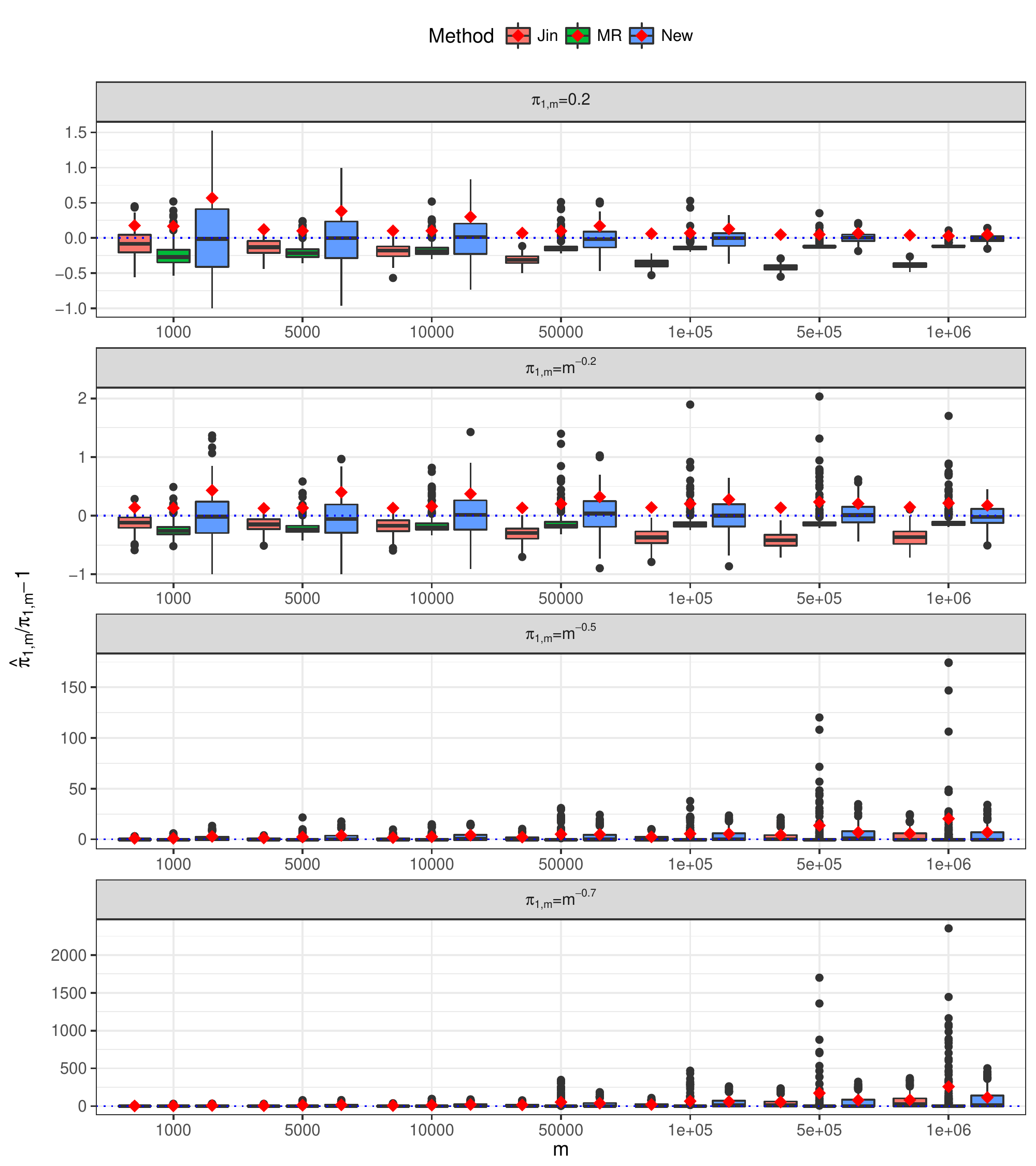}
\caption[Laplace]{Boxplot of the excess $\tilde{\delta}_{m}=\hat{\pi}_{1,m}\pi_{1,m}^{-1}-1$ of an estimator $\hat{\pi}_{1,m}$ when it is applied to Laplace family. The boxplot for the proposed estimator is the right one in each triple of boxplots for each $m$. The thick horizontal line and the diamond in each boxplot are respectively the mean and standard deviation of $\tilde{\delta}_{m}$, and the dotted horizontal line in each panel corresponding to a setting of $\pi_{m}$ is the reference for $\tilde{\delta}_{m}=0$.}
\label{figLaplace}%
\end{figure}

\begin{figure}[H]
\centering
\includegraphics[height=0.8\textheight]{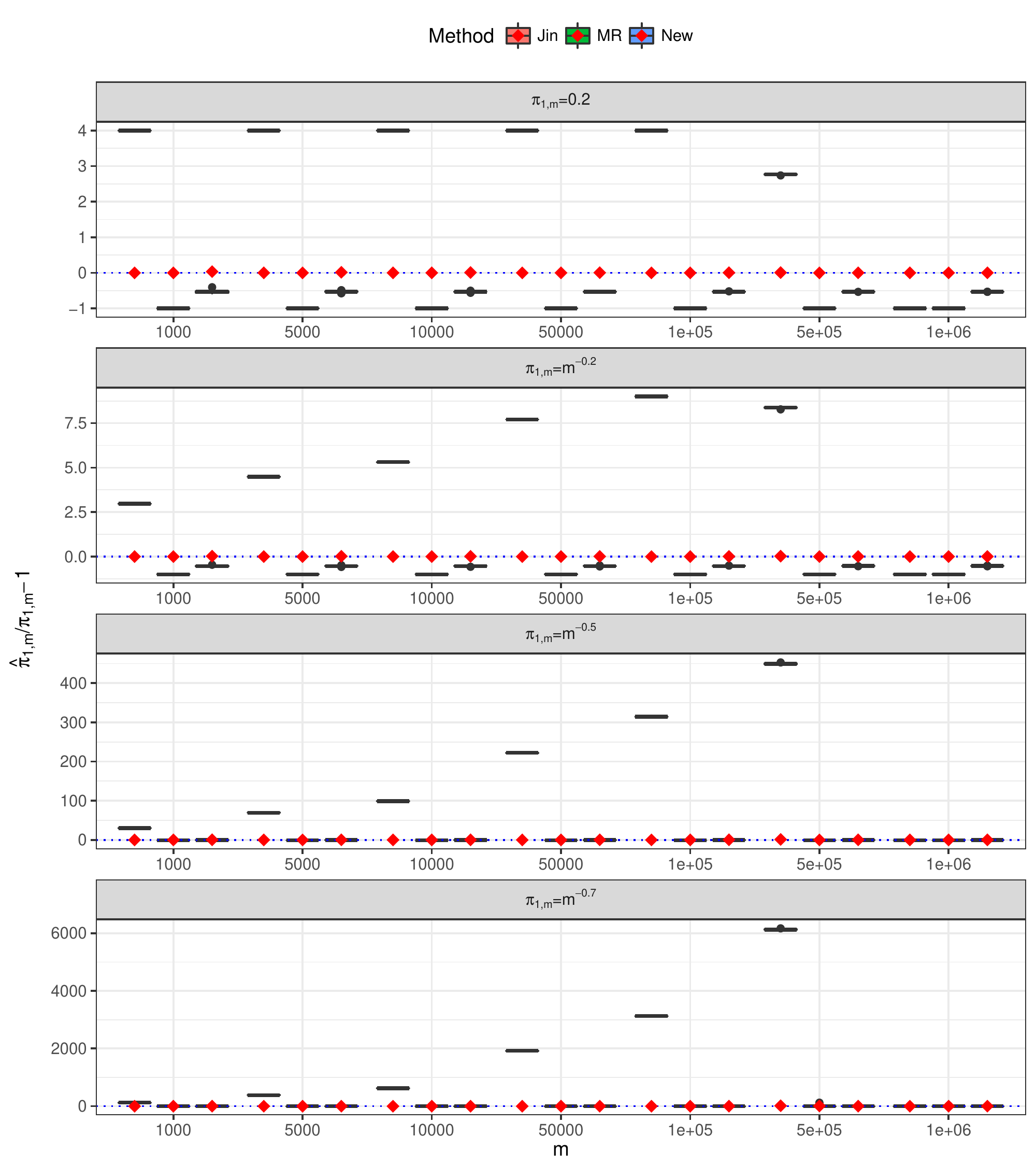}
\caption[NegativeBinomial]{Boxplot of the excess $\tilde{\delta}_{m}=\hat{\pi}_{1,m}\pi_{1,m}^{-1}-1$ of an estimator $\hat{\pi}_{1,m}$ when it is applied to Negative Binomial family. The boxplot for the proposed estimator is the right one in each triple of boxplots for each $m$. The thick horizontal line and the diamond in each boxplot are respectively the mean and standard deviation of $\tilde{\delta}_{m}$, and the dotted horizontal line in each panel corresponding to a setting of $\pi_{1,m}$ is the reference for $\tilde{\delta}_{m}=0$.}
\label{figNegBin}%
\end{figure}

\begin{figure}[H]
\centering
\includegraphics[height=0.8\textheight]{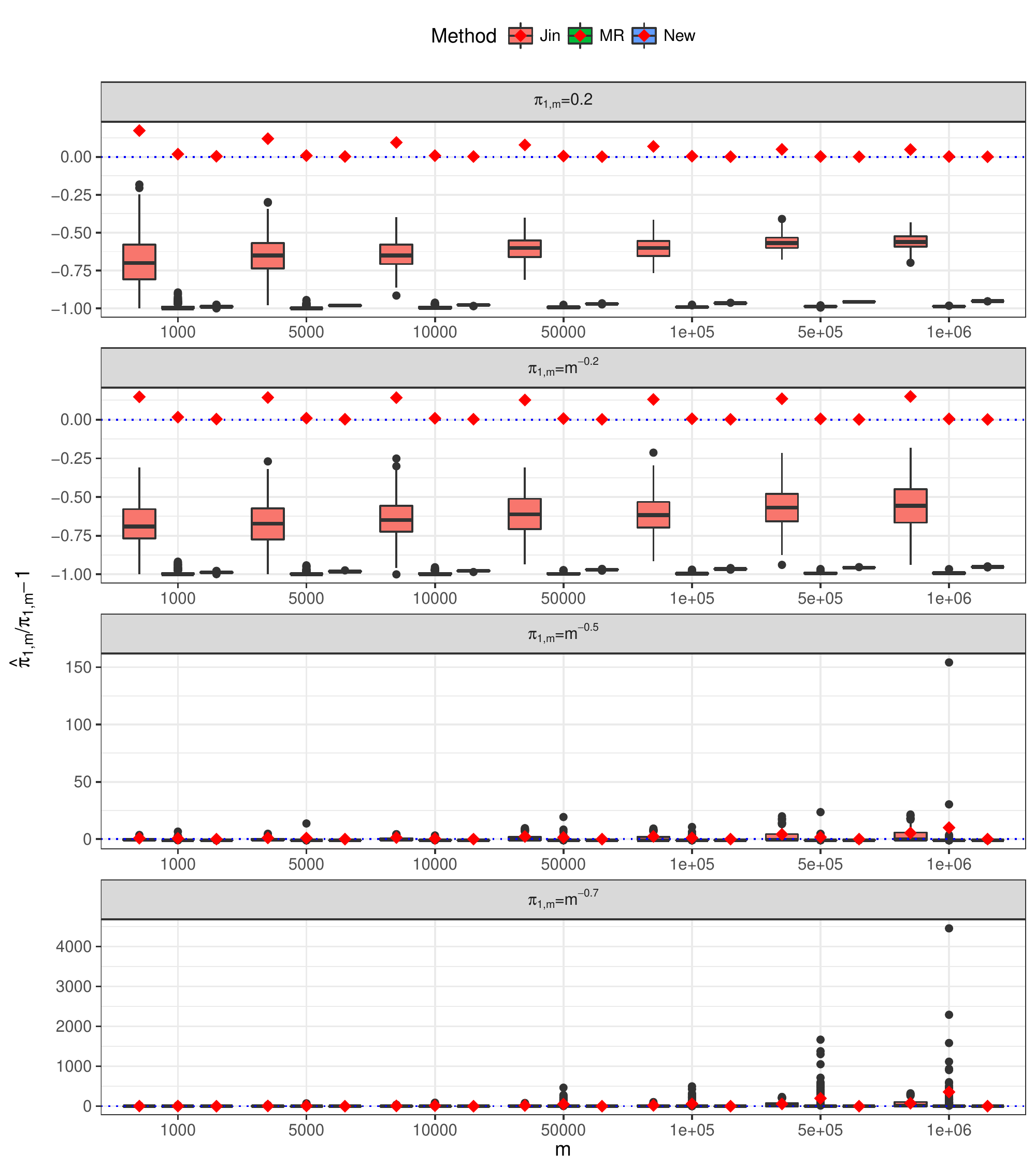}
\caption[Chi-square]{Boxplot of the excess $\tilde{\delta}_{m}=\hat{\pi}_{1,m}\pi_{1,m}^{-1}-1$ of an estimator $\hat{\pi}_{1,m}$ when it is applied to central Chi-square family. The boxplot for the proposed estimator is the right one in each triple of boxplots for each $m$. The thick horizontal line and the diamond in each boxplot are respectively the mean and standard deviation of $\tilde{\delta}_{m}$, and the dotted horizontal line in each panel corresponding to a setting of $\pi_{1,m}$ is the reference for $\tilde{\delta}_{m}=0$.}
\label{figChisq}
\end{figure}

\begin{figure}[H]
\centering
\includegraphics[height=0.8\textheight]{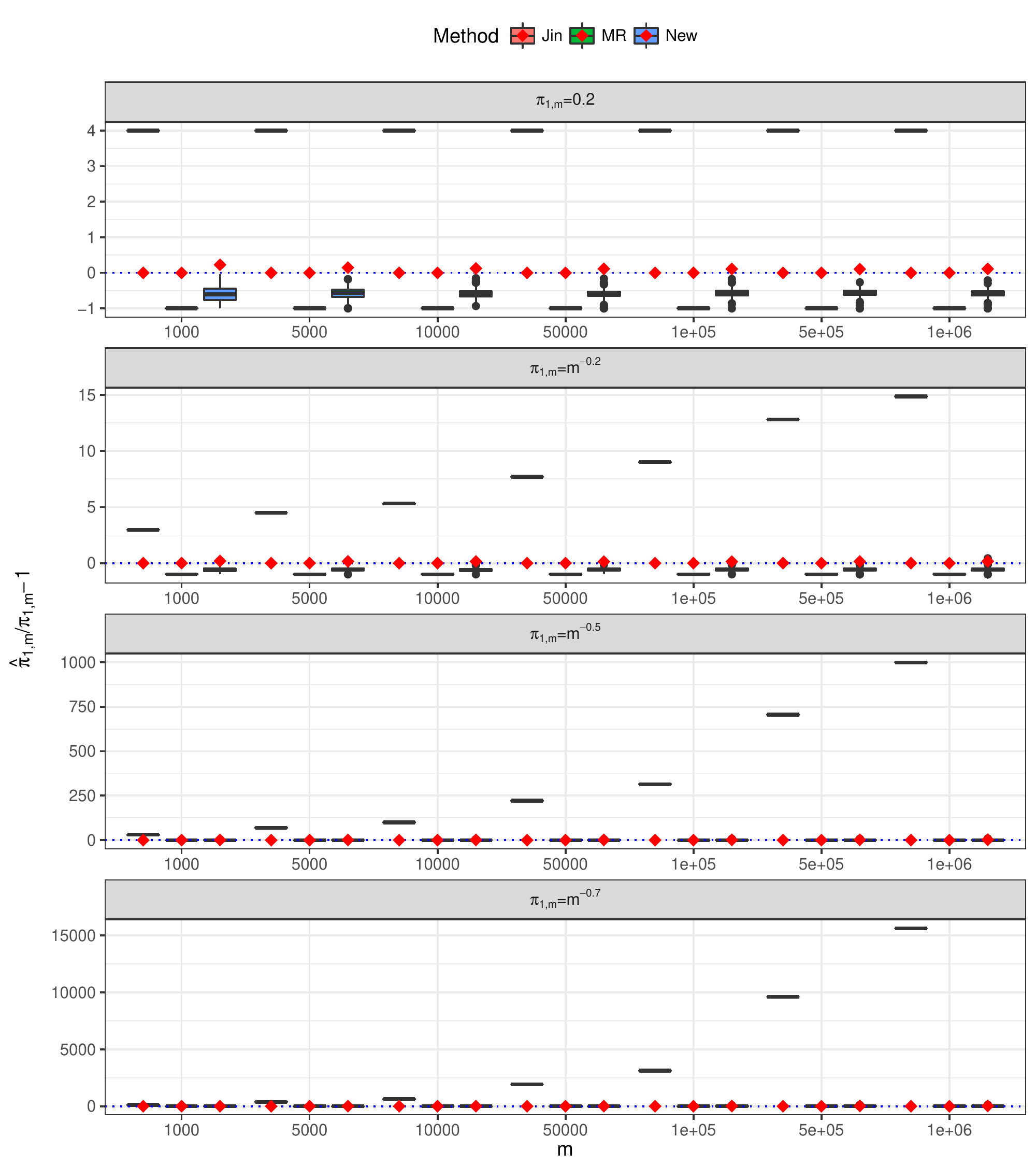}
\caption[Poisson]{Boxplot of the excess $\tilde{\delta}_{m}=\hat{\pi}_{1,m}\pi_{1,m}^{-1}-1$ of an estimator $\hat{\pi}_{1,m}$ when it is applied to Poisson family. The boxplot for the proposed estimator is the right one in each triple of boxplots for each $m$. The thick horizontal line and the diamond in each boxplot are respectively the mean and standard deviation of $\tilde{\delta}_{m}$, and the dotted horizontal line in each panel corresponding to a setting of $\pi_{1,m}$ is the reference for $\tilde{\delta}_{m}=0$.}
\label{figPoisson}%
\end{figure}

\end{document}